\documentclass[12pt, a4paper, reqno]{amsart}
\usepackage{amsmath}
\usepackage{amsfonts}
\usepackage{amssymb}
\usepackage{color}
\usepackage{comment,graphicx}
\usepackage[T1]{fontenc}
\usepackage{url}
\usepackage{ae}
\usepackage{mathtools}
\usepackage{esint}
\usepackage{latexsym}
\usepackage{amscd}

\def\phi{\varphi}
\def\rho{\varrho}
\def\epsilon{\varepsilon}

\numberwithin{equation}{section}
\theoremstyle{plain}
\newtheorem{thm}[equation]{Theorem}
\newtheorem{lem}[equation]{Lemma}
\newtheorem{prop}[equation]{Proposition}
\newtheorem{cor}[equation]{Corollary}
\theoremstyle{definition}
\newtheorem{defn}[equation]{Definition}

\theoremstyle{remark}
\newtheorem{rem}[equation]{Remark}
\newtheorem{ex}[equation]{Example}

\renewcommand{\leq}{\leqslant}
\renewcommand{\geq}{\geqslant}

\DeclareMathOperator{\supp}{supp}

\begin{document}
\title[Weigthed  Herz-type Besov and Triebel-Lizorkin spaces]{ Herz-type Besov and Triebel-Lizorkin spaces with
general weights}
\author[D. Drihem]{Douadi Drihem}
\address{Douadi Drihem\\
M'sila University\\
Department of Mathematics\\
Laboratory of Functional Analysis and Geometry of Spaces\\
M'sila 28000, Algeria.}
\email{douadidr@yahoo.fr, douadi.drihem@univ-msila.dz}
\thanks{ }
\date{\today }

\begin{abstract}
In this paper, we introduce a new family of function spaces of generalized smoothness that encompasses the classical Besov and Triebel-Lizorkin spaces, their weighted counterparts, as well as Herz-type Besov and Triebel-Lizorkin spaces. We provide a characterization of these spaces via the $\varphi$-transform in the sense of Frazier and Jawerth, and establish the corresponding Sobolev embeddings. Furthermore, we obtain smooth atomic and molecular decompositions for these spaces. Our approach relies on several new maximal inequalities.
\end{abstract}

\maketitle

\section{Introduction}

The spaces of generalized smoothness are have been introduced by several
authors. We refer, for instance, to Bownik \cite{M07}, Cobos and Fernandez 
\cite{CF88}, Goldman \cite{Go79} and \cite{Go83}, and Kalyabin \cite{Ka83};
see also\ Besov \cite{B03} and \cite{B05}, and Kalyabin and Lizorkin \cite%
{Kl87}. More general Besov spaces with variable smoothness were explicitly
studied by Ansorena and Blasco \cite{AB95} and \cite{AB96}, including
characterizations by differences and atomic decomposition. The wavelet
decomposition (with respect to a compactly supported wavelet basis of
Daubechies type) of nonhomogeneous Besov spaces of generalized smoothness
was achieved by Almeida \cite{Al}. The theory of these spaces had a
remarkable development in part due to its usefulness in applications. For
instance, they appear in the study of trace spaces on fractals, see Edmunds
and Triebel \cite{ET96} and \cite{ET99}, were they introduced the spaces $%
B_{p,q}^{s,\Psi }$, where $\Psi $ is a so-called admissible function,
typically of log-type near $0$. For a complete treatment of these spaces we
refer the reader the work of Moura \cite{Mo01}. Further results on Besov
spaces of variable smoothness are given in \cite{AA16}.

More general function spaces of generalized smoothness can be found in
Caetano and Leopold \cite{CL}, Farkas and Leopold \cite{FL06}, and reference
therein.

 Dominguez and Tikhonov \cite{DT} gave a treatment of function
spaces with logarithmic smoothness (Besov, Sobolev, Triebel-Lizorkin),
including various new characterizations for Besov norms in terms of
different, sharp estimates for Besov norms of derivatives and potential
operators (Riesz and Bessel potentials) in terms of norms of functions
themselves and sharp embeddings between the Besov spaces defined by
differences and by Fourier-analytical decompositions as well as between
Besov and Sobolev/Triebel-Lizorkin spaces.

Tyulenev introduced in \cite{Ty15}, \cite{Ty-N-L} and \cite{Ty-151} a new
family of Besov spaces of variable smoothness which cover many classes of
Besov spaces, where the norm on these spaces was defined\ with the help of
classical differences. 

Based on this weighted framework, the author in \cite{D20} and \cite{D20.1} introduced Besov and Triebel-Lizorkin spaces of variable smoothness. Several fundamental properties of these spaces were established, including their $\varphi$-transform characterization in the sense of Frazier and Jawerth, as well as the corresponding Sobolev  embeddings. Furthermore, smooth atomic and molecular decompositions for these spaces were obtained. Further properties of such spaces are given in  \cite{D7},\cite{D20.2}, \cite{DRu25} and \cite{DFA26}.

Based on Herz spaces, Besov and Triebel-Lizorkin spaces, the authors in \cite%
{XuYang05} and \cite{Xu05} introduced a new family of function spaces called
Herz-type Besov spaces $\dot{K}_{p}^{\alpha ,q}B_{\beta }^{s}$ and Herz-type
Triebel-Lizorkin spaces $\dot{K}_{p}^{\alpha ,q}F_{\beta }^{s}$. These
spaces unify and generalize the classical Besov spaces and Triebel-Lizorkin
spaces. Several basic properties were established, such as the Fourier
analytical characterisation, lifting properties and embeddings properties.
When $\alpha =0$ and$\ p=q$ the spaces $\dot{K}_{p}^{0,p}B_{\beta }^{s}$ and 
$\dot{K}_{p}^{0,p}F_{\beta }^{s}$ coincide with the usual function spaces $%
B_{p,\beta }^{s}$ and $F_{p,\beta }^{s}$, respectively. The embeddings\ in $%
\dot{K}_{p}^{\alpha ,q}B_{\beta }^{s}$ and $\dot{K}_{p}^{\alpha ,q}F_{\beta
}^{s}$ such us Sobolev, Franke and Jewarth, extend and improve Sobolev,
Franke and Jewarth embeddings of Besov and Triebel-Lizorkin spaces; see \cite%
{Drihem1.13}, \cite{Drihem2.13},\cite{Drmana} and \cite{drihem2016jawerth}.

The interest in Herz-type Besov-Triebel-Lizorkin spaces comes not only from
theoretical reasons but also from their applications to several classical
problems in analysis. In \cite{LuYang97}, Lu and Yang introduced the
Herz-type Sobolev and Bessel potential spaces. They gave some applications
to partial differential equations. In \cite{Dr22.Banach} the author studied
the Cauchy problem for the semilinear parabolic equations%
\begin{equation*}
\partial _{t}u-\Delta u=G(u)
\end{equation*}%
with initial data in Herz-type Triebel-Lizorkin spaces and under some
suitable conditions on $G$.

In this paper, we study Herz-type Besov and  Triebel-Lizorkin spaces of generalised smoothness. We present the $\phi$-transform characterization of these spaces in the sense of
Frazier and Jawerth and we prove their Sobolev embeddings. Also, we establish the
smooth atomic and molecular decomposition of these function spaces. To do these
we need a generalization of some maximal inequality to the case of general weights.

The paper is organized into seven sections. In Section 2, we present some preliminaries and recall basic facts about the Muckenhoupt classes and Tyulenev's weighted class. We also establish some key technical lemmas that will be needed in the proofs of our main results.

In Section 3, we present a new version of the weighted vector-valued maximal inequality of Fefferman and Stein in the Herz setting. In Section 4, we establish several fundamental properties of the function spaces under consideration, including their $\varphi$-transform characterization. In Section 5, we extend the well-known Sobolev embeddings to these function spaces. Section 6 is devoted to the smooth atomic and molecular decompositions of these spaces.

Finally, in Section 7, we consider the inhomogeneous versions of the spaces under consideration and outline the corresponding analogous results..

\section{Basic tools}

 We begin by describing some notation and conventions we will need throughout the paper.

\subsection{Notation and conventions}

We denote by $\mathbb{R}^{n}$ the $n$-dimensional real Euclidean space, $%
\mathbb{N}$ the collection of all natural numbers and $\mathbb{N}_{0}=%
\mathbb{N}\cup \{0\}$. The letter $\mathbb{Z}$ stands for the set of all
integer numbers.\ The expression $f\lesssim g$ means that $f\leq c\,g$ for
some independent constant $c$ (and non-negative functions $f$ and $g$), and $%
f\approx g$ means $f\lesssim g\lesssim f$.\vskip5pt

By $\supp f$, we denote the support of the function $f$, that is, the closure of the set of points where $f$ is nonzero. If $E\subset {\mathbb{R}^{n}}$ is a measurable set, then $%
|E|$ stands for the (Lebesgue) measure of $E$ and $\chi _{E}$ denotes its
characteristic function. By $c$ we denote generic positive constants, which
may have different values at different occurrences. \vskip5pt

A weight is a nonnegative locally integrable function on $\mathbb{R}^{n}$
that takes values in $(0,\infty )$ almost everywhere. For measurable set $%
E\subset \mathbb{R}^{n}$ and a weight $\gamma $, $\gamma (E)$ denotes 
\begin{equation*}
\int_{E}\gamma (x)dx.
\end{equation*}%
Given a measurable set $E\subset \mathbb{R}^{n}$ and $0<p\leq \infty $, we
denote by $L_{p}(E)$ the space of all functions $f:E\rightarrow \mathbb{C}$
equipped with the quasi-norm 
\begin{equation*}
\big\|f\big\|_{L_{p}(E)}=\Big(\int_{E}\left\vert f(x)\right\vert ^{p}dx\Big)%
^{1/p}<\infty ,
\end{equation*}%
with $0<p<\infty $ and%
\begin{equation*}
\big\|f\big\|_{L_{\infty }(E)}=\underset{x\in E}{\text{ess-sup}}\left\vert
f(x)\right\vert <\infty .
\end{equation*}%
If $E=\mathbb{R}^{n}$ then we put $\big\|f\big\|_{L_{p}(\mathbb{R}^{n})}$=$\big\|f\big\|_{p}$. 
For a function $f$ in $L_{1}^{\mathrm{loc}}(\mathbb{R}^{n})$, we set 
\begin{equation*}
M_{A}(f)=\frac{1}{|A|}\int_{A}\left\vert f(x)\right\vert dx
\end{equation*}%
for any $A\subset \mathbb{R}^{n}$. Furthermore, we put%
\begin{equation*}
M_{A,p}(f)=\Big(\frac{1}{|A|}\int_{A}\left\vert f(x)\right\vert ^{p}dx\Big)%
^{1/p},
\end{equation*}%
with $0<p<\infty $. If $1\leq p\leq \infty $ and $1/p+1/p^{\prime }=1$, then $p^{\prime }$ is
called the conjugate exponent of $p$.

Let $0<p,q\leq \infty $. The space $L_{p}(\ell _{q})$ is defined to be the
set of all sequences $\{f_{k}\}_{k\in \mathbb{Z}}$ of functions such that%
\begin{equation*}
\big\|\{f_{k}\}_{k\in \mathbb{Z}}\big\|_{L_{p}(\ell _{q})}=\Big\|\Big(\sum_{k=-\infty }^{\infty
}|f_{k}|^{q}\Big)^{1/q}\Big\|_{p}<\infty
\end{equation*}%
with the usual modifications if $q=\infty $ and if $\{t_{k}\}_{k\in \mathbb{Z}}$ is a sequence
of functions then%
\begin{equation*}
\big\|\{f_{k}\}_{k\in \mathbb{Z}}\big\|_{L_{p}(\ell _{q},\{t_{k}\}_{k\in \mathbb{Z}})}=\big\|\{t_{k}f_{k}\}_{k\in \mathbb{Z}}\big\|_{L_{p}(\ell _{q})}.
\end{equation*}
In what follows, $Q$ will denote a cube in the space $\mathbb{R}^{n}$\ with
sides parallel to the coordinate axes and $l(Q)$\ will denote the side
length of the cube $Q$. For all cubes $Q$ and $r>0$, let $rQ$ be the cube
concentric with $Q$ having side length $rl(Q)$. For $v\in \mathbb{Z}$ and $%
m\in \mathbb{Z}^{n}$, denote by $Q_{v,m}$ the dyadic cube,%
\begin{equation*}
Q_{v,m}=2^{-v}([0,1)^{n}+m).
\end{equation*}%
For the collection of all such cubes we use $\mathcal{Q}=\{Q_{v,m}:v\in 
\mathbb{Z},m\in \mathbb{Z}^{n}\}$. For each cube $Q$, we denote by $x_{v,m}$
the lower left-corner $2^{-v}m$ of $Q=Q_{v,m}$.

\subsection{Muckenhoupt weights}

The purpose of this subsection is to review some known properties of\
Muckenhoupt class.

\begin{defn}
Let $1<p<\infty $. We say that a weight $\gamma $ belongs to the Muckenhoupt
class $A_{p}(\mathbb{R}^{n})$ if there exists a constant $C>0$ such that for
every cube $Q$ the following inequality holds 
\begin{equation}
M_{Q}(\gamma )M_{Q,p^{\prime }/p}(\gamma ^{-1})\leq C.  \label{Ap-constant}
\end{equation}
\end{defn}

The smallest constant $C$ for which $\mathrm{\eqref{Ap-constant}}$ holds is
denoted by $A_{p}(\gamma )$. As an example, we can take%
\begin{equation*}
\gamma (x)=|x|^{\alpha },\quad \alpha \in \mathbb{R}.
\end{equation*}%
Then $\gamma \in A_{p}(\mathbb{R}^{n})$, $1<p<\infty $, if and only if $%
-n<\alpha <n(p-1)$.

For $p=1$ we rewrite the above definition in the following way.

\begin{defn}
We say that a weight $\gamma $ belongs to the Muckenhoupt class $A_{1}( 
\mathbb{R}^{n})$ if there exists a constant $C>0$ such that for every cube $%
Q $\ and for a.e.\ $y\in Q$ the following inequality holds 
\begin{equation}
M_{Q}(\gamma )\leq C\gamma (y).  \label{A1-constant}
\end{equation}
\end{defn}

The smallest constant $C$ for which $\mathrm{\eqref{A1-constant}}$ holds is
denoted by $A_{1}(\gamma )$. The above classes have been first studied by
Muckenhoupt\ \cite{Mu72} and use to characterize the boundedness of the
Hardy-Littlewood maximal function on weighted Lebesgue spaces, see the monographs 
\cite{GR85} and \cite{L. Graf14}\ for a complete account on the theory of
Muckenhoupt weights.

We recall a few basic properties of the class of $A_{p}(\mathbb{R}^{n})$
weights, see \cite[Chapter 7]{Du01}, \cite[Chapter 7]{L. Graf14} and \cite[%
Chapter 5]{St93}.

\begin{lem}
\label{Ap-Property}Let $1\leq p<\infty $.\newline
$\mathrm{(i)}$ If $\gamma \in A_{p}(\mathbb{R}^{n})$, then for any $1\leq
p<q $, $\gamma \in A_{q}(\mathbb{R}^{n})$.\newline
$\mathrm{(ii)}$ Let $1<p<\infty $. $\gamma \in A_{p}(\mathbb{R}^{n})$ if and
only if $\gamma ^{1-p^{\prime }}\in A_{p^{\prime }}(\mathbb{R}^{n})$.\newline
$\mathrm{(iii)}$ Let $\gamma \in A_{p}(\mathbb{R}^{n})$. There is $C>0$ such
that for any cube $Q$ and a measurable subset $E\subset Q$ 
\begin{equation*}
\Big(\frac{|E|}{|Q|}\Big)^{p-1}M_{Q}(\gamma )\leq CM_{E}(\gamma ).
\end{equation*}
$\mathrm{(iv)}$ Suppose that $\gamma \in A_{p}(\mathbb{R}^{n})$ for some $%
1<p<\infty $. Then there exists a $1<p_{1}<p<\infty $ such that $\gamma \in
A_{p_{1}}(\mathbb{R}^{n})$.\newline
$\mathrm{(v)}$ If $\gamma \in A_{p}(\mathbb{R}^{n})$, then for any $%
0<\varepsilon \leq 1$, $\gamma ^{\varepsilon }\in A_{p}(\mathbb{R}^{n})$.%
\newline
$\mathrm{(vi)}$ Let $1\leq p<\infty $ and $\gamma \in A_{p}(\mathbb{R}^{n})$%
. Then there exist $\delta \in (0,1)$ and $C>0$ depending only on $n$, $p$,
and $A_{p}(\gamma )$ such that for any cube $Q$ and any measurable subset $S$
of $Q$ we have 
\begin{equation*}
\frac{M_{S}(\gamma )}{M_{Q}(\gamma )}\leq C\Big(\frac{|S|}{|Q|}\Big)%
^{\delta-1}.
\end{equation*}
\end{lem}

The following theorem gives a useful property of $A_{p}(\mathbb{R}^{n})$
weights (reverse H\"{o}lder inequality), see \cite[Chapter 7]{L. Graf14} or 
\cite[Chapter 1]{LuDiYa07}.

\begin{thm}
\label{reverse Holder inequality}Let $1\leq p<\infty \ $and $\gamma \in
A_{p}(\mathbb{R}^{n})$. Then there exist  constants $C>0$ and $\varepsilon
>0$, depending only on $p$ and the $A_{p}(\mathbb{R}^{n})$-constant of $%
\gamma $, such that for every cube $Q$, 
\begin{equation*}
M_{Q,1+\varepsilon }(\gamma )\leq CM_{Q}(\gamma ).
\end{equation*}
\end{thm}

\subsection{The weight class $\dot{X}_{\protect\upsilon ,\protect\sigma ,p}$}

Let $0<p\leq \infty $. A weight sequence $\{t_{k}\}_{k\in \mathbb{Z}}$ is called $p$%
-admissible if $t_{k}\in L_{p}^{\mathrm{loc}}(\mathbb{R}^{n})$ for all $k\in 
\mathbb{Z}$. We mention here that 
\begin{equation*}
\int_{E}t_{k}^{p}(x)dx<c(k)
\end{equation*}%
for any\ compact set $E\subset \mathbb{R}^{n}$. For a $p$-admissible weight
sequence $\{t_{k}\}_{k\in \mathbb{Z}}$\ we set%
\begin{equation*}
t_{k,m}=\big\|t_{k}|\big\|_{L_{p}(Q_{k,m})},\quad k\in \mathbb{Z},m\in \mathbb{Z%
}^{n}.
\end{equation*}

Tyulenev\ \cite{Ty14} introduced the following new weighted class\ and used
it to study Besov spaces of variable smoothness.

\begin{defn}
\label{Tyulenev-class}Let $\upsilon _{1}$, $\upsilon_{2}\in \mathbb{R}$, $%
p,\sigma _{1}$, $\sigma _{2}$ $\in (0,+\infty ]$, $\upsilon =(\upsilon
_{1},\upsilon _{2})$ and let $\sigma =(\sigma _{1},\sigma _{2})$. We let $\dot{%
		X}_{\upsilon ,\sigma ,p}=\dot{X}_{\upsilon ,\sigma ,p}(\mathbb{R}^{n})$ denote
the set of $p$-admissible weight sequences $\{t_{k}\}_{k\in \mathbb{Z}}$ satisfying the
following conditions. There exist numbers $C_{1},C_{2}>0$ such that for any $%
k\leq j$\ and every cube $Q,$ 
\begin{eqnarray}
M_{Q,p}(t_{k})M_{Q,\sigma _{1}}(t_{j}^{-1}) &\leq &C_{1}2^{\alpha _{1}(k-j)},
\label{Asum1} \\
M_{Q,p}^{-1}(t_{k})M_{Q,\sigma _{2}}(t_{j}) &\leq &C_{2}2^{\alpha _{2}(j-k)}, 
\label{Asum2}
\end{eqnarray}
\end{defn}
the modifications of \eqref{Asum1} or \eqref{Asum2} for $\sigma _{1}=\infty$ or $\sigma _{2}=\infty$ are clear. The constants $C_{1},C_{2}>0$ are independent of both the indexes $k$ and $j$%
.

\begin{rem}
$\mathrm{(i)}$\ We would like to mention that if $\{t_{k}\}_{k\in \mathbb{Z}}$ satisfies $%
\mathrm{\eqref{Asum1}}$ with $\sigma _{1}=r\left( p/r\right) ^{\prime }$ and 
$0<r<p<\infty $, then $t_{k}^{p}\in A_{p/r}(\mathbb{R}^{n})$ for any $k\in 
\mathbb{Z}$.\newline
$\mathrm{(ii)}$ We say that $t_{k}\in A_{p}(\mathbb{R}^{n})$,\ $k\in \mathbb{%
\ Z}$, $1\leq p<\infty ,$ have the same Muckenhoupt constant if 
\begin{equation*}
A_{p}(t_{k})=c,\quad k\in \mathbb{Z},
\end{equation*}%
where $c$ is independent of $k$.\newline
$\mathrm{(iii)}$ Definition \ref{Tyulenev-class} is different from the one
used in \cite[Definition 2.1]{Ty14} and Definition 2.7 in \cite{Ty15},
because we used the boundedness of the maximal function on weighted Lebesgue
spaces.
\end{rem}

\begin{ex}
\label{Example1}Let $0<r<p<\infty $, a weight $\omega ^{p}\in A_{\frac{p}{r}
}(\mathbb{R}^{n})$ and 
\begin{equation*}
\{s_{k}\}_{k\in \mathbb{Z}}=\{2^{ks}\omega ^{p}(2^{-k}\cdot )\}_{k\in \mathbb{Z}},\quad s\in 
\mathbb{R}.
\end{equation*}
Clearly, $\{s_{k}\}_{k\in \mathbb{Z}}$ lies in $\dot{X}_{\upsilon ,\sigma ,p}$
for $\upsilon _{1}=\upsilon _{2}=s$, $\sigma =(r(\frac{p}{r})^{\prime },p)$.
\end{ex}

\begin{rem}
\label{Tyulenev-class-properties}Let $0<\theta \leq p<\infty $. Let $\upsilon
_{1}$, $\upsilon _{2}\in \mathbb{R}$, $\sigma _{1},\sigma _{2}\in (0,+\infty )$
, $\sigma _{2}\geq p$, $\upsilon =(\upsilon _{1},\upsilon _{2})$ and let $\sigma
=(\sigma _{1}=\theta \left( \frac{p}{\theta }\right) ^{\prime },\sigma _{2})$
. Let a $p$-admissible weight sequence $\{t_{k}\}_{k\in \mathbb{Z}}\in \dot{X}_{\upsilon ,\sigma
,p}$. Then 
\begin{equation*}
\upsilon _{2}\geq \upsilon _{1},
\end{equation*}
see \cite{D20}.
\end{rem}


Further notation will be properly introduced whenever needed.

\subsection{Auxiliary results}

In this subsection we present some results which are useful for us. Let
recall the vector-valued maximal inequality of Fefferman and Stein \cite%
{FeSt71}. As usual, we put%
\begin{equation*}
\mathcal{M}(f)(x)=\sup_{Q}M_{Q}(f),\quad f\in L_{1}^{\mathrm{loc}}(\mathbb{R}%
^{n}),
\end{equation*}%
where the supremum\ is taken over all cubes with sides parallel to the axis
and $x\in Q$. Also we set 
\begin{equation*}
\mathcal{M}_{\sigma }(f)=\sup_{Q}M_{Q,\sigma }(f),\quad 0<\sigma <\infty .
\end{equation*}%
Observe that $\mathcal{M}_{\sigma }(f)$ can be rewritten as%
\begin{equation*}
\mathcal{M}_{\sigma }(f)=(\mathcal{M}(|f|^{\sigma }))^{1/\sigma },\quad
0<\sigma <\infty .
\end{equation*}

\begin{thm}
\label{Maximal}Let $1<p\leq \infty $. Then 
\begin{equation*}
\big\|\mathcal{M}(f)\big\|_{p}\lesssim \big\|f\big\|_{p}
\end{equation*}
holds for all $f\in L_{p}(\mathbb{R}^{n})$.
\end{thm}

For the proof see \cite[Chapter 7]{L. Graf14}. Now, we state the
vector-valued maximal inequality of Fefferman and Stein \cite{FeSt71}.

\begin{thm}
\label{FS-inequality}Let $0<p<\infty ,0<\beta \leq \infty $ and $0<\nu
<\min (p,\beta)$. Then 
\begin{equation}
\Big\|\Big(\sum\limits_{j=-\infty }^{\infty }\big(\mathcal{M}_{\nu
}(f_{j})\big)^{\beta }\Big)^{1/\beta }\Big\|_{p}\lesssim \Big\|\Big(%
\sum\limits_{j=-\infty }^{\infty }|f_{j}|^{\beta }\Big)%
^{1/\beta }\Big\|_{p}  \label{Fe-St71}
\end{equation}%
holds for all sequence of functions $\{f_{k}\}_{k\in \mathbb{Z}}\in L_{p}(\ell _{\beta})$.
\end{thm}

We shall require the following theorem, that is the Fefferman-Stein's
inequality, see \cite{FeSt71}.

\begin{lem}
	\label{FS-lemma}Let $1<p<\infty $. Given a non-negative real valued
	functions $f$ and $g$. We have 
	\begin{equation*}
	\int_{\mathbb{R}^{n}}\left( \mathcal{M}(f)(x)\right) ^{p}g(x)dx\leq c\int_{ 
		\mathbb{R}^{n}}\left( f(x)\right) ^{p}\mathcal{M}(g)(x)dx,
	\end{equation*}
	with $c$ independent of $f$ and $g$.
\end{lem}

For the proof see \cite[Chapter 7]{L. Graf14}.

We recall the following Hardy type inequality.

\begin{lem}
\label{lem:lq-inequality}Let $0<a<1$ and $0<q\leq \infty $. Let $%
\{\varepsilon _{k}\}_{k\in \mathbb{Z}}$ be a sequence of positive real
numbers, such that 
\begin{equation*}
\big\|\{\varepsilon _{k}\}_{k\in \mathbb{Z}}\big\|_{\ell ^{q}}=I<\infty .
\end{equation*}%
Then the sequences 
\begin{equation*}
\big\{\delta _{k}:\delta _{k}=\sum_{j=-\infty }^{k}a^{k-j}\varepsilon _{j}%
\big\}_{k\in \mathbb{Z}}\quad \text{and}\quad \big\{\eta _{k}:\eta
_{k}=\sum_{j=k}^{\infty }a^{j-k}\varepsilon _{j}\big\}_{k\in \mathbb{Z}}
\end{equation*}%
belong to $\ell ^{q}$, and 
\begin{equation*}
\big\|\{\delta _{k}\}_{k\in \mathbb{Z}}\big\|_{\ell ^{q}}+\big\|\{\eta
_{k}\}_{k\in \mathbb{Z}}\big\|_{\ell ^{q}}\leq C\text{\ }I,
\end{equation*}%
with 
\begin{equation*}
C=2\Big(\frac{1}{1-a^{\min (1,q)}}\Big)^{1/(\min(1,q)}.
\end{equation*}
\end{lem}
\section{Maximal\ inequalities}

The aim of this section is to present some weighted vector-valued inequalities for the Hardy-Littlewood maximal function on Herz spaces. First, we present some fundamental properties of Herz spaces. Let $%
k\in \mathbb{Z}$. For convenience, we set 
\begin{equation*}
B_{k}=B(0,2^{k}), \qquad R_{k}=B_{k}\setminus B_{k-1}.
\end{equation*}
We denote by $\chi_{k}$ the characteristic function of the set $R_{k}$.

\begin{defn}
	\label{def-inh-Herz} Let $\alpha \in \mathbb{R}$ and $0< p,q\leq \infty$. The homogeneous Herz space $\dot{K}%
	_{p}^{\alpha,q}(\mathbb{R}^{n})$ is defined as the set of all functions $%
	f\in L_{p}^{\mathrm{loc}}(\mathbb{R}^{n}\setminus\{0\})$ such that 
	\begin{equation*}
	\|f\|_{\dot{K}_{p}^{\alpha,q}(\mathbb{R}^{n})} = \left(
	\sum_{k=-\infty}^{\infty} 2^{k\alpha q} \|f\chi_{k}\|_{p}^{q} \right)^{1/q}
	<\infty,
	\end{equation*}
	with the usual modifications when $p=\infty$ and/or $q=\infty$.

\end{defn}

\begin{rem}
	The spaces $\dot{K}_{p}^{\alpha,q}(%
	\mathbb{R}^{n})$ are quasi-Banach spaces. If $\alpha=0$ and $0< p=q\leq \infty$%
	, then  $\dot{K}_{p}^{0,p}(\mathbb{R}^{n})$ coincides with the Lebesgue space $L^{p}(\mathbb{R}^{n})$. Moreover, 
	\begin{equation*}
	\dot{K}_{p}^{\alpha,p}(\mathbb{R}^{n}) = L^{p}(\mathbb{R}^{n},|\cdot|^{%
		\alpha p}),
	\end{equation*}
	that is, the weighted Lebesgue space endowed with the norm 
	\begin{equation*}
	\|f\|_{L^{p}(\mathbb{R}^{n},|\cdot|^{\alpha p})} = \left( \int_{\mathbb{R}%
		^{n}} |f(x)|^{p}|x|^{\alpha p}\,dx \right)^{1/p}.
	\end{equation*}
	Let $0<p\leq \infty$, $0<q_{1}\leq q_{2}\leq \infty$, and $\alpha\in\mathbb{R}$.
	Then the following continuous embeddings hold:
	\begin{equation}
	\dot{K}_{p}^{\alpha,q_{1}}(\mathbb{R}^{n})
	\hookrightarrow
	\dot{K}_{p}^{\alpha,q_{2}}(\mathbb{R}^{n}).
	\label{herz-emb}
	\end{equation}

A detailed discussion of the properties of these spaces may be found in the
monograph \cite{LYH08}, the papers \cite{LuYang1.95}, \cite{LuYang2.95}, 
\cite{Hernandez1998}, and the references therein.

\end{rem}
\begin{defn}
	Let $\alpha \in \mathbb{R}$ and $0< p,q,\beta\leq \infty$. The space $\dot{K}%
	_{p}^{\alpha,q}(\ell _{\beta})$ is defined to be the
	set of all sequences $\{f_{k}\}_{k\in \mathbb{Z}}$ of functions such that%
	\begin{equation*}
	\big\|\{f_{k}\}_{k\in \mathbb{Z}}\big\|_{\dot{K}%
		_{p}^{\alpha,q}(\ell _{\beta})}=\Big\|\Big(\sum_{k=-\infty }^{\infty
	}|f_{k}|^{\beta}\Big)^{1/\beta}\Big\|_{\dot{K}%
		_{p}^{\alpha,q}}<\infty
	\end{equation*}%
	with the usual modifications if $\beta=\infty $.
	\end{defn}

Consider a mapping $T$ that maps measurable functions on $\mathbb{R}^n$ to
measurable functions on $\mathbb{R}^n$. We say that $T$ is sublinear if it
satisfies the following conditions:

\begin{enumerate}

\item If $f=f_1+f_2$ and $T(f_1)$ and $T(f_2)$ are defined, then $T(f)$ is
also defined and 
\begin{equation*}
|T(f_1+f_2)(x)| \leq |T(f_1)(x)|+|T(f_2)(x)|.
\end{equation*}

\item For any scalar $c$, if $T(f)$ is defined, then $T(cf)$ is also defined
and 
\begin{equation*}
|T(cf)(x)|=|c|\,|T(f)(x)|.
\end{equation*}
\end{enumerate}

	We recall the following interpolation result from Lemma 2.2 of
\cite{AJ80}.

\begin{lem}
	\label{Calderon-Zygmund}
	Let $1\leq p_0,p_1<\infty$, $M_0,M_1>0$, and
	$1<\beta_0,\beta_1<\infty$.
	Suppose that $T$ is a sublinear operator satisfying
	\begin{equation*}
	\left\|\{T(f_k)\}_{k\in\mathbb Z}\right\|_
	{L_{p_i}(\ell_{\beta_i})}
	\leq M_i
	\left\|\{f_k\}_{k\in\mathbb Z}\right\|_
	{L_{p_i}(\ell_{\beta_i})},
	\qquad i=0,1,
	\end{equation*}
	for every $\{f_k\}_{k\in\mathbb Z}\in
	L_{p_i}(\ell_{\beta_i})$, $i=0,1$.
	Then $T$ extends to a bounded operator on
	$L_p(\ell_\beta)$ such that
	\begin{equation*}
	\left\|\{T(f_k)\}_{k\in\mathbb Z}\right\|_
	{L_p(\ell_\beta)}
	\leq M_0^{1-\theta}M_1^\theta
	\left\|\{f_k\}_{k\in\mathbb Z}\right\|_
	{L_p(\ell_\beta)},
	\label{CZ19561}
	\end{equation*}
	for every $\{f_k\}_{k\in\mathbb Z}\in L_p(\ell_\beta)$,
	where
	\begin{equation*}
	\frac{1}{p}=\frac{1-\theta}{p_0}+\frac{\theta}{p_1},
	\qquad
	\frac{1}{\beta}=\frac{1-\theta}{\beta_0}
	+\frac{\theta}{\beta_1},
	\qquad 0<\theta<1.
	\end{equation*}
\end{lem}

Next, we present a weighted version of the preceding lemma.

\begin{lem}
	\label{Calderon-Zygmund1}
	Let $w_k$, $k\in\mathbb Z$, be measurable functions on
	$\mathbb R^n$ such that $w_k>0$ almost everywhere.
	Let $1\leq p_0,p_1<\infty$, $M_0,M_1>0$, and
	$1<\beta_0,\beta_1<\infty$.
	Suppose that $T$ is a sublinear operator satisfying
	\begin{equation*}
	\left\|\{w_kT(f_k)\}_{k\in\mathbb Z}\right\|_
	{L_{p_i}(\ell_{\beta_i})}
	\leq M_i
	\left\|\{w_kf_k\}_{k\in\mathbb Z}\right\|_
	{L_{p_i}(\ell_{\beta_i})},
	\qquad i=0,1,
	\end{equation*}
	for every sequence $\{f_k\}_{k\in\mathbb Z}$ such that
	$\{w_kf_k\}_{k\in\mathbb Z}\in
	L_{p_i}(\ell_{\beta_i})$, $i=0,1$.
	Then $T$ satisfies
	\begin{equation*}
	\left\|\{w_kT(f_k)\}_{k\in\mathbb Z}\right\|_
	{L_p(\ell_\beta)}
	\leq M_0^{1-\theta}M_1^\theta
	\left\|\{w_kf_k\}_{k\in\mathbb Z}\right\|_
	{L_p(\ell_\beta)},
	\label{CZ1957}
	\end{equation*}
	for every sequence $\{f_k\}_{k\in\mathbb Z}$ such that
	$\{w_kf_k\}_{k\in\mathbb Z}\in L_p(\ell_\beta)$, where
	\begin{equation*}
	\frac{1}{p}=\frac{1-\theta}{p_0}+\frac{\theta}{p_1},
	\qquad
	\frac{1}{\beta}=\frac{1-\theta}{\beta_0}
	+\frac{\theta}{\beta_1},
	\qquad 0<\theta<1.
	\end{equation*}
\end{lem}

\begin{proof}
	For each $k\in\mathbb Z$, define
	\begin{equation*}
	H(f_k)=w_kT(w_k^{-1}f_k).
	\end{equation*}
	Since $w_k>0$ almost everywhere, this definition is meaningful.
	Moreover, $H$ is a sublinear operator.
	
	Let $\{g_k\}_{k\in\mathbb Z}\in
	L_{p_i}(\ell_{\beta_i})$, $i=0,1$.
	By the hypotheses on $T$, applied to the sequence
	$\{w_k^{-1}g_k\}_{k\in\mathbb Z}$, we obtain
	\begin{align*}
	\big\|\{H(g_k)\}_{k\in\mathbb Z}\big\|_
	{L_{p_i}(\ell_{\beta_i})}
	&=
	\big\|\{w_kT(w_k^{-1}g_k)\}_{k\in\mathbb Z}\big\|_
	{L_{p_i}(\ell_{\beta_i})}\\
	&\leq M_i
	\big\|\{w_kw_k^{-1}g_k\}_{k\in\mathbb Z}\big\|_
	{L_{p_i}(\ell_{\beta_i})}\\
	&=M_i
	\big\|\{g_k\}_{k\in\mathbb Z}\big\|_
	{L_{p_i}(\ell_{\beta_i})},
	\qquad i=0,1.
	\end{align*}
	
	Consequently, applying Lemma \ref{Calderon-Zygmund} to $H$,
	we obtain
	\begin{equation*}
	\big\|\{H(g_k)\}_{k\in\mathbb Z}\big\|_
	{L_p(\ell_\beta)}
	\leq M_0^{1-\theta}M_1^\theta
	\big\|\{g_k\}_{k\in\mathbb Z}\big\|_
	{L_p(\ell_\beta)}.
	\end{equation*}
	
	Now, let $\{f_k\}_{k\in\mathbb Z}$ be such that
	$\{w_kf_k\}_{k\in\mathbb Z}\in L_p(\ell_\beta)$.
	Taking $g_k=w_kf_k$, we obtain
	\begin{align*}
	\big\|\{w_kT(f_k)\}_{k\in\mathbb Z}\big\|_
	{L_p(\ell_\beta)}
	&=
	\big\|\{H(w_kf_k)\}_{k\in\mathbb Z}\big\|_
	{L_p(\ell_\beta)}\\
	&\leq M_0^{1-\theta}M_1^\theta
	\big\|\{w_kf_k\}_{k\in\mathbb Z}\big\|_
	{L_p(\ell_\beta)}.
	\end{align*}
	This completes the proof.
\end{proof}

\begin{rem}
	The proof of Lemma \ref{Calderon-Zygmund1} provides an alternative
	proof of Lemma 2.13 in \cite{D20.1}.
\end{rem}
\begin{defn}
	Let $A$ be a Banach space, $0<q\leq\infty$, and
	$\alpha\in\mathbb{R}$. The homogeneous Herz-type space
	$\dot{K}^{\alpha,q}(A)$ is defined by
	\[
	\dot{K}^{\alpha,q}(A)
	=
	\left\{
	f:\ f\chi_k\in A \text{ for all } k\in\mathbb{Z},
	\ \|f\|_{\dot{K}^{\alpha,q}(A)}<\infty
	\right\},
	\]
	where
	\[
	\|f\|_{\dot{K}^{\alpha,q}(A)}
	=
	\left(
	\sum_{k\in\mathbb{Z}}
	2^{k\alpha q}\|f\chi_k\|_A^q
	\right)^{1/q},
	\]
	with the usual modification when $q=\infty$.
\end{defn}
We now turn to the theory of real interpolation. We recall the definition of $K$-functional. Let $(A_{0},A_{1})$ be a compatible couple of normed spaces.
For any $f\in A_{0}+A_{1}$, the $K$-functional is defined as
\begin{equation*}
K(t,f;A_{0},A_{1})=\inf_{f=f_{0}+f_{1}}\big(\big\|f_{0}\big\|_{A_{0}}+t\big\|%
f_{1}\big\|_{A_{1}}\big),\quad t>0,
\end{equation*}%
is the $K$-functional.
 \begin{defn}
	\label{real1}Let $0<\theta <1$ and $1\leq q\leq\infty $. Let $(A_{0},A_{1})$\
	be a compatible couple of normed spaces. The space $(A_{0},A_{1})_{\theta ,q}$ consists of
	all $f$ in $A_{0}+A_{1}$ for which the functional%
	\begin{equation*}
	\big\|f\big\|_{(A_{0},A_{1})_{\theta ,q}}=\Big(\int_{0}^{\infty }t^{-\theta
		q}(K(t,f))^{q}\frac{dt}{t}\Big)^{1/q},
	\end{equation*}%
	is finite, with the usual modification when $q=\infty$.
\end{defn}
For a comprehensive treatment of the theory of real interpolation, the reader is referred to \cite{BL76} and \cite{T78}.

 The following interpolation theorem plays a central role in this section; see \cite{Dr-Biskra}.

\begin{thm}
	\label{InterpolationHerz}
	Let $A$ be a Banach ideal space, $0<\theta<1$,
	$1\leq q_0,q_1,q\leq\infty$, and
	$\alpha_0,\alpha_1\in\mathbb{R}$. If
	$\alpha_0\neq\alpha_1$, then
	\[
	\big(
	\dot{K}^{\alpha_0,q_0}(A),
	\dot{K}^{\alpha_1,q_1}(A)
	\big)_{\theta,q}
	=
	\dot{K}^{\alpha,q}(A),
	\]
	where
	\[
	\alpha=(1-\theta)\alpha_0+\theta\alpha_1.
	\]
	The equality is understood up to equivalence of norms.
\end{thm}

We shall say, as usual, that an operator $T$ mapping a vector
space $A$ into $L^0(\Omega,\mu)$, the space of $\mu$-measurable
functions on a measure space $(\Omega,\mu)$, is
$C$-subadditive if $T(f+g)$ is defined whenever $Tf$ and $Tg$
are defined, and if
\[
|T(f+g)(t)|
\leq C\bigl(|Tf(t)|+|Tg(t)|\bigr)
\quad \text{a.e.},
\]
where $C$ is a positive constant independent of $a$ and $b$.
If, moreover,
\[
|T(\lambda f)(t)|
=|\lambda|\,|Tf(t)|
\quad \text{a.e.}
\]
for any scalar $\lambda$, then $T$ is called
$C$-sublinear (sublinear when $C=1$). We need the following result; see Theorem 6 in \cite{Maligranda}.
\begin{thm}
	\label{Lech}
	Let $1\leq q\leq \infty$. Let $B=(B_0,B_1)$ be a couple of Banach lattices on
	a measure space $(\Omega,\mu)$. Let
	\[
	T:A_0+A_1\longrightarrow L^0(\Omega,\mu)
	\]
	satisfy the following conditions:

		\item[(i)] $T$ is $C$-subadditive.
		\item[(ii)] For $i=0,1$, if $f\in A_i$, then
		\[
		Tf\in B_i
		\quad\text{and}\quad
		\|Tf\|_{B_i}\leq M_i\|f\|_{A_i}.
		\]
	
	Then $T$ maps $A_{\theta,q}$ into $B_{\theta,q}$ and
	\[
	\|Tf\|_{B_{\theta,q}}
	\leq C M_0^{1-\theta}M_1^\theta
	\|f\|_{A_{\theta,q}}.
	\]
\end{thm}

We shall also need the following result; see \cite{D20.1} and \cite{DRu25}.

\begin{lem}
\label{key-estimate2}Let $1<\theta \leq p<\infty $\ and $1<\beta \leq\infty $.
Let $\{t_{j}\}_{j\in \mathbb{Z}}$\ be a $p$-admissible\ weight\ sequence\ such that $%
t_{j}^{p}\in A_{\frac{p}{\theta }}(\mathbb{R}^{n})$, $j\in \mathbb{Z}$.
Assume that $t_{j}^{p}$,\ $j\in \mathbb{Z}$ have the same Muckenhoupt
constant. Then 
\begin{equation*}
\Big\|\Big(\sum\limits_{j=-\infty }^{\infty }t_{j}^{\beta }\big(\mathcal{M}%
(f_{j})\big)^{\beta }\Big)^{1/\beta }\Big\|_{p}\lesssim \Big\|\Big(%
\sum\limits_{j=-\infty }^{\infty }t_{j}^{\beta }|f_{j}|^{\beta }\Big)%
^{1/\beta }\Big\|_{p}
\end{equation*}%
holds for all sequences of functions $\{t_{j}f_{j}\}_{j\in \mathbb{Z}}\in L_{p}(\ell _{\beta
}) $, with the usual modification when $\beta=\infty$. In particular	\begin{equation}
\big\|t_{j}\mathcal{M}(f_{j})\big\|_{p}\leq c\big\| %
t_{j}f_{j}\big\|_{p}  \label{key-est1herz}
\end{equation}
for all sequences of functions\ $f_{j}\in L_{p}(\mathbb{R}^{n},t_{j})$, $%
j\in \mathbb{Z}$, where $c>0$ is independent of $j$.
\end{lem}
\begin{rem}
	\label{r-estimatescopy(2)} 
	$\newline
	\mathrm{(i)}$ Let $1<\theta \leq p<\infty $. Let $\{t_{j}\}_{j\in \mathbb{Z}}$\ be a $p$-admissible\ weight\ sequence\ such that $%
	t_{j}^{p}\in A_{\frac{p}{\theta }}(\mathbb{R}^{n})$, $j\in \mathbb{Z}$.
	Assume that $t_{j}^{p}$,\ $j\in \mathbb{Z}$ have the same Muckenhoupt
	constant. The property \eqref{key-est1herz} can be generalized
	in the following way.\newline
	$\mathrm{\bullet }$ If $\{t_{j}^{p}\}_{j\in \mathbb{Z}}$ satisfies $\mathrm{\ %
		\eqref{Asum1}}$, then 
	\begin{equation*}
	\big\|t_{k}\mathcal{M}(f_{j})\big\|_{p}\leq c\text{ }%
	2^{\alpha _{1}(k-j)}\big\|t_{j}f_{j}\big\|_{p}
	\end{equation*}%
	holds for all sequences of functions $f_{j}\in L_{p}(\mathbb{R}^{n},t_{j})$, 
	$j\in \mathbb{Z}$ and $j\geq k$, where $c>0$ is independent of $k$ and $j$,
	see \cite{D20}. $\newline
	\mathrm{\bullet }$ If $\{t_{j}^{p}\}_{j\in \mathbb{Z}}$ satisfies $\mathrm{\ %
		\eqref{Asum2}}$ with $\sigma _{2}\geq p$, then 
	\begin{equation*}
	\big\|t_{k}\mathcal{M}(f_{j})\big\|_{p}\leq c\text{ }%
	2^{\alpha _{2}(k-j)}\big\|t_{j}f_{j}\big\|_{p}
	\end{equation*}%
	holds for all sequences of functions $f_{j}\in L_{p}(\mathbb{R}^{n},t_{j})$, 
	$j\in \mathbb{Z}$ and $k\geq j$, where $c>0$ is independent of $k$ and $j$,
	see \cite{D20}.
	 $\newline
	\mathrm{(ii)}$\ In Lemma \ref{key-estimate2}, we can assume that $%
	t_{j}^{p}\in A_{p}(\mathbb{R}^{n})$,\ $j\in \mathbb{Z}$, $1<p<\infty $ with 
	\begin{equation*}
	A_{p}(t_{j}^{p})\leq c,\quad k\in \mathbb{Z},
	\end{equation*}%
	which follows by Lemma \ref{Ap-Property}/(iii).
\end{rem}
Next, we present the first main result of this section.

\begin{lem}
\label{key-estimate1}Let $1<\theta \leq p<\infty$, $1\leq q  \leq\infty $ and $1<\beta <\infty $%
. Let $\{t_{j}\}_{j\in \mathbb{Z}}$\ be a $p$-admissible\ weight\ sequence\ such that $%
t_{j}^{p}\in A_{\frac{p}{\theta }}(\mathbb{R}^{n})$, $j\in \mathbb{Z}$.
Assume that $t_{j}^{p}$,\ $j\in \mathbb{Z}$ have the same Muckenhoupt
constant. Let $\delta $ be as in Lemma \ref{Ap-Property}. Let $-\frac{%
n\delta }{p}<\alpha <n-\frac{n}{p}$. Then 
\begin{equation}
\Big\|\Big(\sum\limits_{j=-\infty }^{\infty }t_{j}^{\beta }\big(\mathcal{M}%
(f_{j})\big)^{\beta }\Big)^{1/\beta }\Big\|_{\dot{K}_{p}^{\alpha ,q}}\lesssim %
\Big\|\Big(\sum\limits_{j=-\infty }^{\infty }t_{j}^{\beta }|f_{j}|^{\beta }%
\Big)^{1/\beta }\Big\|_{\dot{K}_{p}^{\alpha ,q}}  \label{key-est}
\end{equation}%
holds for all sequences of functions $\{t_{j}f_{j}\}_{j\in \mathbb{Z}}\in
\dot{K}_{p}^{\alpha ,q}(\ell _{\beta })$.
\end{lem}

\begin{proof}
We will do the proof into two steps.

\textit{Step 1.} We prove $\mathrm{\eqref{key-est}}$ with $1<\beta <p<\infty 
$. Let $\kappa \in \mathbb{N}$ be such that $\beta ^{\kappa }\leq p<\beta
^{\kappa +1}$.

\textit{Substep 1.1. }We consider the case $\beta ^{\kappa }<p<\beta
^{\kappa +1}$. In order to prove we additionally do it into the two steps
Substeps 1.1.1, 1.1.2 and 1.1.3.

\textit{Substep 1.1.1.} We consider the case $1<\beta \leq \theta ^{\frac{1}{%
\kappa }}$. From Lemma \ref{Ap-Property}/(i), $t_{j}^{p}\in A_{\frac{p}{\beta^{\kappa }}
}(\mathbb{R}^{n})$, $j\in \mathbb{Z}$. Observing that

\begin{equation*}
\Big\|\Big(\sum\limits_{j=-\infty }^{\infty }t_{j}^{\beta }\big(\mathcal{M}%
(f_{j})\big)^{\beta }\Big)^{1/\beta }\chi _{k}\Big\|_{p},\quad k\in \mathbb{Z%
}
\end{equation*}
can be estimated from above by 
\begin{equation*}
\sup \sum\limits_{j=-\infty }^{\infty }\int_{R_{k}}t_{j}(x)\mathcal{\ M}%
(f_{j})(x)|g_{j}(x)|dx=\sup \sum\limits_{j=-\infty }^{\infty }T_{j,k},
\end{equation*}%
where the supremum is taken over all sequence of functions $\{g_{j}\chi_{k}\}_{j\in 
\mathbb{Z}}\in L_{p^{\prime }}(\ell _{\beta ^{\prime }})$ with 
\begin{equation}
\big\|\{g_{j}\chi_{k}\}_{j\in \mathbb{Z}}\big\|_{L_{p^{\prime }}(\ell _{\beta
^{\prime }})}\leq 1,\label{dual}
\end{equation}%
where $p^{\prime }$ and $\beta ^{\prime }$ are the conjugate exponent of $p$
and $\beta $, respectively. We split the operator $\mathcal{M}(f_{j}),j\in 
\mathbb{Z}$ into 
\begin{equation}
\mathcal{M}(f_{j})\leq \mathcal{M}(f_{j}\chi _{B_{k-2}})+\mathcal{M}%
(f_{j}\chi _{\widetilde{R}_{k}})+\mathcal{M}(f_{j}\chi _{{\mathbb{R}^{n}}%
\setminus B_{k+2}}),  \label{splitM}
\end{equation}%
where $\widetilde{R}_{k}=\left\{ x\in {\mathbb{R}^{n}}:\,2^{k-2}\leq
|x|<2^{k+2}\right\} $ with $j,k\in \mathbb{Z}$. We have $T_{j,k}\leq
\sum_{m=1}^{3}T_{m,j,k}$, where 
\begin{equation*}
T_{1,j,k}=\int_{R_{k}}t_{j}(x)\mathcal{M}(f_{j}\chi
_{B_{k-2}})(x)|g_{j}(x)|dx,\quad T_{2,j,k}=\int_{R_{k}}t_{j}(x)\mathcal{M}%
(f_{j}\chi _{\widetilde{R}_{k}})(x)|g_{j}(x)|dx
\end{equation*}%
and%
\begin{equation*}
T_{3,j,k}=\int_{R_{k}}t_{j}(x)\mathcal{M}(f_{j}\chi _{{\mathbb{R}^{n}}%
\setminus B_{k+2}})(x)|g_{j}(x)|dx.
\end{equation*}%
Hence it suffices to prove that 
\begin{equation}
\Big\|\Big\{2^{k\alpha }\sup\Big( \sum\limits_{j=-\infty }^{\infty }T_{m,j,k}\Big)\Big\}%
_{k\in \mathbb{Z}}\Big\|_{\ell _{q}}\lesssim \Big\|\Big(\sum\limits_{j=-%
\infty }^{\infty }t_{j}^{\beta }|f_{j}|^{\beta }\Big)^{1/\beta }\Big\|%
_{\dot{K}_{p}^{\alpha ,q}},\quad m\in \{1,2,3\},  \label{main-est}
\end{equation}%
where the supremum is taken over all sequence of functions $\{g_{j}\chi_{k}\}_{j\in 
	\mathbb{Z}}\in L_{p^{\prime }}(\ell _{\beta ^{\prime }})$ such that \eqref{dual} holds.

\textit{Estimate of $\Big\{2^{k\alpha }\sup\Big(\sum\limits_{j=-\infty }^{\infty
} T_{1,j,k}\Big)\Big\}_{k\in \mathbb{Z}}$ in $\ell _{q}$ norm}. For $x\in R_{k}$
and $y\in B_{k-2}$ we have $|x-y|\geq |x|-|y|>2^{k-2}$. Then 
\begin{align*}
2^{k\alpha }\mathcal{M}(f_{j}\chi _{B_{k-2}})(x)& \lesssim 2^{k(\alpha
-n)}\int_{B_{k-2}}|f_{j}(y)|\,dy \\
& \lesssim \sum_{i=-\infty }^{k-2}2^{(k-i)\alpha -kn}\int_{R_{i}}2^{i\alpha
}|f_{j}(y)|\,dy.
\end{align*}%
This leads to%
\begin{equation*}
2^{k\alpha }T_{1,j,k}\lesssim \sum_{i=-\infty }^{k-2}2^{(k-i)\alpha
-kn}\int_{R_{k}}t_{j}(x)|g_{j}(x)|dx\int_{R_{i}}2^{i\alpha }|f_{j}(y)|\,dy.
\end{equation*}%
Applying H\"{o}lder's inequality, we get 
\begin{align}
\int_{R_{k}}t_{j}(x)|g_{j}(x)|dx& \leq \Big(\int_{R_{k}}t_{j}^{\tau }(x)dx%
\Big)^{1/\tau }\Big(\int_{R_{k}}|g_{j}(x)|^{\tau ^{\prime }}dx\Big)^{1/\tau
^{\prime }}  \label{main-est2} \\
& \leq |R_{k}|M_{B_{k},\tau }(t_{j})M_{B_{k},\tau ^{\prime }}(g_{j}\chi _{k})
\notag
\end{align}%
with $\tau =p(1+\varepsilon )$, where $\varepsilon $ as in Theorem \ref%
{reverse Holder inequality}, which is possible since $t_{j}^{p}\in
A_{\frac{p}{\theta }}(\mathbb{R}^{n})$,  $j \in \mathbb{Z}$. Obviously, we have 
\begin{equation*}
M_{B_{k},\tau }(t_{j})=M_{B_{k},p(1+\varepsilon )}(t_{j})\leq c\text{ }%
M_{B_{k},p}(t_{j}),\quad k,j\in \mathbb{Z}.
\end{equation*}%
Since $t_{j}^{p}$,\ $j\in \mathbb{Z}$ have the same Muckenhoupt constant and
from the proof of Theorem 7.2.2 in \cite{L. Graf14} the constant $c$ is
independent of $k$. From Lemma \ref{Ap-Property}/(iii) we obtain%
\begin{equation*}
M_{B_{k},p}(t_{j})\leq C2^{(i-k)(\frac{n}{p}-n)}M_{B_{i},p}(t_{j}),\quad i\leq k-2.
\end{equation*}%
 Therefore, 
\begin{equation*}
\int_{R_{k}}t_{j}(x)|g_{j}(x)|dx\lesssim 2^{(i-k)(\frac{n}{p}%
-n)}|R_{k}|M_{B_{i},p}(t_{j})M_{B_{k},\tau ^{\prime }}(g_{j}\chi _{k}),
\end{equation*}%
By H\"{o}lder's inequality, we obtain 
\begin{equation*}
M_{R_{i}}(f_{j})\leq M_{B_{i},s^{\prime
}}(t_{j}^{-1})M_{B_{i},s}(t_{j}f_{j}\chi _{i}),\quad s=\frac{p}{\beta
^{\kappa }}.
\end{equation*}%
Again, by H\"{o}lder's inequality with the help of the fact that $\beta
^{\kappa }>1$, we obtain%
\begin{equation*}
\int_{B_{i}}t_{j}^{-s^{\prime }}(x)dx=\int_{B_{i}}t_{j}^{-\frac{s^{\prime }}{%
s}s}(x)dx\leq \Big(\int_{B_{i}}t_{j}^{-\frac{s^{\prime }}{s}p}(x)dx\Big)%
^{s/p}|B_{i}|^{1-\frac{s}{p}},
\end{equation*}%
which yields that 
\begin{equation*}
M_{B_{i},s^{\prime }}(t_{j}^{-1})\leq \big(M_{B_{i},\frac{s^{\prime }}{s}%
}(t_{j}^{-p})\big)^{1/p}.
\end{equation*}%
Hence 
\begin{equation*}
\int_{R_{i}}|f_{j}(y)|\,\chi
_{i}(y)dy\int_{R_{k}}t_{j}(x)|g_{j}(x)|dx\lesssim 2^{(i-k)(\frac{n}{p}%
-n)}|R_{k}||R_{i}|M_{B_{i},s}(t_{j}f_{j}\chi _{i})M_{B_{k},\tau ^{\prime
}}(g_{j}\chi _{k}),
\end{equation*}%
because of $t_{j}^{p}\in A_{s}(\mathbb{R}^{n})$,\ $j\in \mathbb{Z}$. The
last expression is bounded by 
\begin{align*}
&c2^{(i-k)(\frac{n}{p}-n)}|R_{k}|\int_{R_{i}}M_{B_{i},s}(t_{j}f_{j}\chi
_{i})M_{B_{k},\tau ^{\prime }}(g_{j}\chi _{k})dx\\
&\lesssim 2^{(i-k)(\frac{n}{p}-n)}|R_{k}|\int_{%
\mathbb{R}^{n}}\mathcal{M}_{s}(t_{j}f_{j})(x)\mathcal{M}_{\tau ^{\prime
}}(g_{j}\chi _{k})(x)dx.
\end{align*}%
Therefore, 
\begin{equation*}
2^{k\alpha }\sum\limits_{j=-\infty }^{\infty }T_{1,j,k}\lesssim
\sum_{i=-\infty }^{k-2}2^{(k-i)(\alpha -\frac{n}{p}+n)}2^{\alpha i}\int_{\mathbb{R}%
^{n}}\sum\limits_{j=-\infty }^{\infty }\mathcal{M}_{s}(t_{j}f_{j}\chi _{i})(x)\mathcal{M}_{\tau ^{\prime }}(g_{j}\chi _{k})(x)dx.
\end{equation*}%
By H\"{o}lder's inequality the term inside the integral is bounded by 
\begin{equation*}
\Big(\sum\limits_{j=-\infty }^{\infty }\big(\mathcal{M}_{s}\big(t_{j}f_{j}\chi _{i}\big)(x)\big)^{\beta}\Big)^{1/\beta}\Big(\sum\limits_{j=-\infty
}^{\infty }\big(\mathcal{M}_{\tau ^{\prime }}(g_{j}\chi _{k})(x)\big)%
^{\beta^{\prime }}\Big)^{1/\beta^{\prime }}
\end{equation*}%
for any $x\in \mathbb{R}^{n}$. Again by H\"{o}lder's inequality $2^{k\alpha
}\sum\limits_{j=-\infty }^{\infty }T_{1,j,k}$ can be estimated by 
\begin{align*}
& c\sum_{i=-\infty }^{k-2}2^{(k-i)(\alpha -\frac{n}{p}+n)}2^{\alpha i}\Big\|%
\Big(\sum\limits_{j=-\infty }^{\infty }\big(\mathcal{M}_{s}\big(%
t_{j}f_{j}\chi _{i}\big)\big)^{\beta }\Big)^{1/\beta }\Big\|_{p} \\ \times
& \Big\|\Big(\sum\limits_{j=-\infty }^{\infty }\big(\mathcal{M}_{\tau
^{\prime }}(g_{j}\chi _{k})\big)^{\beta ^{\prime }}\Big)^{1/\beta ^{\prime }}%
\Big\|_{p^{\prime }} \\
& \lesssim \sum_{i=-\infty }^{k-2}2^{(k-i)(\alpha -\frac{n}{p}+n)}2^{\alpha
i}\Big\|\Big(\sum\limits_{j=-\infty }^{\infty }\big(\mathcal{M}_{s}\big(%
t_{j}f_{j}\chi _{i}\big)\big)^{\beta }\Big)^{1/\beta }\Big\|_{p} \\
& \lesssim \sum_{i=-\infty }^{k-2}2^{(k-i)(\alpha -\frac{n}{p}+n)}2^{\alpha
i}\Big\|\Big(\sum\limits_{j=-\infty }^{\infty }|t_{j}f_{j}\chi _{i}|^{\beta }%
\Big)^{1/\beta }\Big\|_{p},
\end{align*}%
where in the second and third inequality we have used the vector-valued
maximal inequality of Fefferman and Stein $\mathrm{\eqref{Fe-St71}}$,
because of $0<s<\beta <p$, and the fact that 
\begin{align*}
\Big\|\Big(\sum\limits_{j=-\infty }^{\infty }\big(\mathcal{M}_{\tau ^{\prime
}}(g_{j}\chi _{k})\big)^{\beta ^{\prime }}\Big)^{1/\beta ^{\prime }}\Big\|%
_{p^{\prime }}& \lesssim \big\|\{g_{j}\chi _{k}\}_{j\in \mathbb{Z}} \big\|_{L_{p^{\prime }}(\ell
_{\beta^{\prime }})} \\
& \lesssim 1,
\end{align*}%
since $\tau =p(1+\varepsilon )>p>\beta $.\textrm{\ }Applying Lemma \ref%
{lem:lq-inequality} we obtain $\mathrm{\eqref{main-est}}$ in the case of $%
m=1.$

\textit{Estimate of $\Big\{2^{k\alpha }\sup\Big(\sum\limits_{j=-\infty }^{\infty} T_{2,j,k}\Big)\Big\}_{k\in \mathbb{Z}}$ in $\ell _{q}$ norm}. By H\"{o}lder's inequality $%
\sum\limits_{j=-\infty }^{\infty }T_{2,j,k}$ is bounded above by%
\begin{align*}
&\int_{R_{k}}\Big(\sum\limits_{j=-\infty }^{\infty }t_{j}^{\beta }\big(%
\mathcal{M}(f_{j}\chi _{\widetilde{R}_{k}})\big)^{\beta }\Big)^{1/\beta }%
\Big(\sum\limits_{j=-\infty }^{\infty }|g_{j}|^{\beta ^{\prime
}}\Big)^{1/\beta ^{\prime }}dx \\
&\lesssim \Big\|\Big(\sum\limits_{j=-\infty }^{\infty }t_{j}^{\beta }\big(%
\mathcal{M}(f_{j}\chi _{\widetilde{R}_{k}})\big)^{\beta }\Big)^{1/\beta }%
\Big\|_{p}\Big\|\Big(\sum\limits_{j=-\infty }^{\infty }|g_{j}\chi_{k}|%
^{\beta ^{\prime }}\Big)^{1/\beta ^{\prime }} \Big\|_{p^{\prime }} \\
&\lesssim \Big\|\Big(\sum\limits_{j=-\infty }^{\infty }t_{j}^{\beta }\big(%
\mathcal{M}(f_{j}\chi _{\widetilde{R}_{k}})\big)^{\beta }\Big)^{1/\beta }%
\Big\|_{p}.
\end{align*}%
Applying  Lemma \ref{key-estimate2}, we obtain 
\begin{equation*}
\Big\|\Big(\sum\limits_{j=-\infty }^{\infty }t_{j}^{\beta }\big(\mathcal{M}%
(f_{j}\chi _{\widetilde{R}_{k}})\big)^{\beta }\Big)^{1/\beta }\Big\|_{p}
\end{equation*}%
is bounded by%
\begin{equation*}
c\Big\|\Big(\sum\limits_{j=-\infty }^{\infty }t_{j}^{\beta }|f_{j}|^{\beta }\chi
_{\widetilde{R}_{k}}\Big)^{1/\beta }\Big\|_{p},\quad k\in 
\mathbb{Z}.
\end{equation*}%
This leads to%
\begin{equation*}
\sum\limits_{j=-\infty }^{\infty }T_{2,j,k}\lesssim \Big\|\Big(%
\sum\limits_{j=-\infty }^{\infty }t_{j}^{\beta }|f_{j}|^{\beta }%
\Big)^{1/\beta }\chi _{\widetilde{R}_{k}}\Big\|_{p}.
\end{equation*}%
Multiplying by $2^{k\alpha}$ and taking the power $q$ of both sides and summing
over $k\in \mathbb{Z}$ we obtain $\mathrm{\eqref{main-est}}$ in the case of $%
m=2\mathrm{.}$

\textit{Estimate of $\Big\{2^{k\alpha }\sup\Big(\sum\limits_{j=-\infty }^{\infty} T_{3,j,k}\Big)\Big\}_{k\in \mathbb{Z}}$ in $\ell _{q}$ norm}. Noting
that $|x-y|>2^{i-3}$ for $x\in R_{k}$ and $y\in R_{i}$, with $i\geq k+3$.
This implies that 
\begin{equation*}
2^{k{\alpha }}\mathcal{M}(f_{j}\chi _{{\mathbb{R}^{n}}\setminus
B_{k+2}})(x)\lesssim \sum_{i=k+3}^{\infty }2^{(k-i)\alpha
-in}\,\int_{R_{i}}2^{i\alpha }|f_{j}(y)|dy.
\end{equation*}
Recall that 
\begin{equation}
\int_{R_{k}}t_{j}(x)|g_{j}(x)|dx\leq \Big(\int_{R_{k}}t_{j}^{\tau }(x)dx\Big)%
^{1/\tau }\Big(\int_{R_{k}}|g_{j}(x)|^{\tau ^{\prime }}dx\Big)^{1/\tau
^{\prime }},  \label{main-est3}
\end{equation}%
with $\tau =p(1+\varepsilon )$, where $\varepsilon $ as in Theorem \ref%
{reverse Holder inequality}. We choose $\varepsilon $ small enough that \ $%
\alpha >-\frac{n\delta }{p(1+\varepsilon )}$, which is possible since $%
\alpha >-\frac{n\delta }{p}$; see also Theorem 7.2.2 in \cite{L. Graf14}. Since $t_{j}^{\tau }\in A_{p}(\mathbb{R}^{n})$, $j\in \mathbb{Z}$,
Lemma \ref{Ap-Property}/(vi) gives%
\begin{equation*}
\Big(\int_{B_{k}}t_{j}^{\tau }(x)dx\Big)^{1/\tau }\lesssim 2^{\frac{%
(k-i)n\delta }{\tau }}\Big(\int_{B_{i}}t_{j}^{\tau }(x)dx\Big)^{1/\tau },\quad i \geq k+3.
\end{equation*}%
Inserting this estimate in $\mathrm{\eqref{main-est3}}$, we find%
\begin{equation*}
\int_{R_{k}}t_{j}(x)|g_{j}(x)|dx\lesssim 2^{\frac{(k-i)n\delta }{\tau }%
}|R_{i}|M_{B_{i},\tau }(t_{j})M_{B_{i},\tau ^{\prime }}(g_{j}\chi _{k}),
\end{equation*}%
where the implicit constant is independent of $k$ and $i$. Obviously, we
have 
\begin{equation*}
M_{B_{i},\tau }(t_{j})=M_{B_{i},p(1+\varepsilon )}(t_{j})\leq c\text{ }%
M_{B_{i},p}(t_{j}).
\end{equation*}%
Therefore 
\begin{equation*}
\int_{B_{i}}|f_{j}(y)|\,\chi
_{i}(y)dy\int_{R_{k}}t_{j}(x)|g_{j}(x)|dx\lesssim 2^{\frac{(k-i)n\delta }{%
\tau }}|R_{i}|^{2}M_{B_{i},s}(t_{j}f_{j}\chi _{i})M_{B_{i},\tau ^{\prime
}}(g_{j}\chi _{k}).
\end{equation*}%
As in the estimate of $2^{k\alpha }\sup\Big(\sum\limits_{j=-\infty }^{\infty
} T_{1,j,k}\Big),k\in \mathbb{Z}$ we get
that $2^{k\alpha }\sup\Big(\sum\limits_{j=-\infty }^{\infty
} T_{3,j,k}\Big),k\in \mathbb{Z}$ can be
estimated from above by 
\begin{equation*}
c\sum_{i=k+3}^{\infty }2^{(k-i)(\alpha +\frac{n\delta }{p(1+\varepsilon )}%
)}2^{\alpha i}\Big\|\Big(\sum\limits_{j=-\infty }^{\infty }|
t_{j}f_{j}\chi _{i}|^{\beta }\Big)^{1/\beta }\Big\|_{p}.
\end{equation*}%
Applying Lemma \ref{lem:lq-inequality} we obtain $\mathrm{\eqref{main-est}}$
in the case of $m=3\mathrm{.}$

\textit{Substep 1.1.2.} We prove $\mathrm{\eqref{key-est}}$ with $\beta =p$. By
Lemma \ref{Ap-Property}/(iv)\ there exists $1<p_{3}<p<\infty $ such that $%
t_{k}^{p}\in A_{p_{3}}(\mathbb{R}^{n}),k\in \mathbb{Z}$. To study\ this
case, we proceed exactly as in Substep 1.1.1 but with $p_{3}$ in place of $s$

\textit{Substep 1.1.3.} We consider the case $\theta ^{\frac{1}{\kappa }%
}<\beta <p$. Let 
\begin{equation*}
\frac{1}{\beta }=\frac{1-\lambda }{p}+\frac{\lambda }{\theta ^{\frac{1}{%
			\kappa }}},\quad 0<\lambda <1.
\end{equation*}%
Let $h\in \mathbb{R}$ be such that  $0<h<\min\Big(n-\frac{n}{p}-\alpha,\alpha+\dfrac{n\delta}{p}\Big)$. We put $\alpha _{1}=\alpha -h$ and $\alpha _{2}=\alpha+h$. Then

\begin{equation*}
-\frac{n\delta }{p}<\alpha _{i}<n-\frac{n}{p},\quad i\in \{1,2\}\label{alpha4}
\end{equation*}
and%
\begin{equation*}
\alpha =(1-\frac{1}{2} )\alpha _{1}+\frac{1}{2} \alpha _{2}\label{alpha5},\quad \alpha _{1}\neq\alpha _{2}. 
\end{equation*}
Substeps 1.1.1 and 1.1.2 give
\begin{equation*}
\Big\|\Big(\sum\limits_{j=-\infty }^{\infty }t_{j}^{p}\big(\mathcal{M}(f_{j})%
\big)^{p}\Big)^{1/p}\Big\|_{\dot{K}_{p}^{\alpha_{i} ,p}}\lesssim \Big\|\Big(%
\sum\limits_{j=-\infty }^{\infty }t_{j}^{p}\left\vert f_{j}\right\vert ^{p}%
\Big)^{1/p}\Big\|_{\dot{K}_{p}^{\alpha_{i} ,p}},
\end{equation*}%
and 
\begin{equation*}
\Big\|\Big(\sum\limits_{j=-\infty }^{\infty }t_{j}^{\theta ^{\frac{1}{\kappa 
}}}\big(\mathcal{M}(f_{j})\big)^{\theta ^{\frac{1}{\kappa }}}\Big)^{1/\theta
	^{\frac{1}{\kappa }}}\Big\|_{\dot{K}_{p}^{\alpha_{i} ,p}}\lesssim \Big\|\Big(%
\sum\limits_{j=-\infty }^{\infty }t_{j}^{\theta ^{\frac{1}{\kappa }%
}}\left\vert f_{j}\right\vert ^{\theta ^{\frac{1}{\kappa }}}\Big)^{1/\theta
	^{\frac{1}{\kappa }}}\Big\|_{\dot{K}_{p}^{\alpha_{i} ,p}}
\end{equation*}%
for any $i\in\{1,2\}$. Applying Lemma \ref{Calderon-Zygmund1} we obtain 
\begin{equation*}
\Big\|\Big(\sum\limits_{j=-\infty }^{\infty }t_{j}^{\beta}\big(\mathcal{M}(f_{j})%
\big)^{\beta}\Big)^{1/\beta}\Big\|_{\dot{K}_{p}^{\alpha_{i} ,p}}\lesssim \Big\|\Big(%
\sum\limits_{j=-\infty }^{\infty }t_{j}^{\beta}\left\vert f_{j}\right\vert ^{\beta}%
\Big)^{1/\beta}\Big\|_{\dot{K}_{p}^{\alpha_{i} ,p}}
\end{equation*}
for any $i\in\{1,2\}$. By Theorem \ref{Lech}, we get 
\begin{equation*}
\Big\|\{t_{k}\mathcal{M}%
(f_{k})\}_{k\in \mathbb{Z}}\Big\|_{(\dot{K}^{\alpha _{1},p}(L_{p}(\ell _{\beta})),\dot{K}^{\alpha _{2},p}(L_{p}(\ell _{\beta}))_{\frac{1}{2} ,q}}\lesssim \Big\|\{t_{k}
f_{k}\}_{k\in \mathbb{Z}}\Big\|_{(\dot{K}^{\alpha _{1},p}(L_{p}(\ell _{\beta})),\dot{K}^{\alpha _{2},p}(L_{p}(\ell _{\beta}))_{\frac{1}{2} ,q}}.
\end{equation*}%
Using Theorem 	\ref{InterpolationHerz}, we obtain
\begin{equation*}
\Big\|\Big(\sum\limits_{k=-\infty }^{\infty }t_{k}^{\beta }\big(\mathcal{M}%
(f_{k})\big)^{\beta }\Big)^{1/\beta }\Big\|_{\dot{K}_{p}^{\alpha ,q}}\lesssim %
\Big\|\Big(\sum\limits_{k=-\infty }^{\infty }t_{k}^{\beta }\left\vert
f_{k}\right\vert ^{\beta }\Big)^{1/\beta }\Big\|_{\dot{K}_{p}^{\alpha ,q}},
\end{equation*}
which is the desired estimate.

\textit{Substep 1.2. }We have proved $\mathrm{\eqref{key-est}}$ in case of $%
\beta ^{\kappa }<p<\beta ^{\kappa +1}$. Now we study the case $p=\beta
^{\kappa },\kappa \in \mathbb{N}$. First observing that $\kappa \geq 2$. Since $t_{k}^{p}\in
A_{p}(\mathbb{R}^{n}),k\in \mathbb{Z}$, from Theorem \ref{reverse Holder
inequality} there exists a number $\gamma _{1}>0$ such 
\begin{equation*}
M_{Q,1+\gamma _{1}}(t_{k}^{p})\lesssim M_{Q}(t_{k}^{p}),\quad k\in \mathbb{Z}%
,
\end{equation*}%
where the implicit constant is independent of $k$, see, e.g., \cite[Theorem
7.2.5 and Corollary 7.2.6.]{L. Graf14}. In addition by Lemma \ref%
{Ap-Property}/(ii), $t_{k}^{-p^{\prime }}\in A_{p^{\prime }}(\mathbb{R}^{n})$
and again by Theorem \ref{reverse Holder inequality} there exists a number $%
\gamma _{2}>0$ such 
\begin{equation*}
M_{Q,1+\gamma _{2}}(t_{k}^{-p^{\prime }})\lesssim M_{Q}(t_{k}^{-p^{\prime
}}),\quad k\in \mathbb{Z}.
\end{equation*}%
Let 
\begin{equation}
0<\gamma <\min \Big(\frac{p(n-\alpha )-n}{\max (0,\alpha )},\gamma
_{1},\gamma _{2},\frac{p-\frac{p}{\beta} }{\frac{p}{\beta} -1}\Big)\label{reserve}.
\end{equation}%
Then 
\begin{equation*}
M_{Q,1+\gamma }(t_{k}^{p})\big(M_{Q,1+\gamma }(t_{k}^{-p^{\prime }})\big)^{%
\frac{p}{p^{\prime }}}\lesssim 1,\quad k\in \mathbb{Z},
\end{equation*}%
where the implicit constant is independent of $k$. Hence $t_{k}^{p(1+\gamma
)}\in A_{p}(\mathbb{R}^{n}),k\in \mathbb{Z}$. Therefore $t_{k}^{p}\in
A_{p_{1}}(\mathbb{R}^{n}),k\in \mathbb{Z}$, where 
\begin{equation*}
p_{1}=\frac{p+\gamma }{1+\gamma },
\end{equation*}%
see  \cite[Exercise 7.1.3]{L. Graf14}. From \eqref{reserve}
\begin{equation*}
\beta ^{\kappa-1}<p_{1}<p=\beta ^{\kappa}.
\end{equation*}%
In addition by Lemma %
\ref{Ap-Property}/(v), we obtain $t_{k}^{p_{1}}\in A_{p_{1}}(\mathbb{R}^{n}),k\in \mathbb{Z}$ . The same procedure yields that $t_{k}^{p(1+\nu )}\in
A_{p}(\mathbb{R}^{n})\subset A_{p(1+\nu )}(\mathbb{R}^{n}),k\in \mathbb{Z}$
with 
\begin{equation}
0<\nu <\min \Big(\gamma _{1},\gamma _{2},\beta-1,\frac{n\delta}{p\max(0,-\alpha)}-1\Big)\label{v}
\end{equation}%
and then 
\begin{equation*}
\beta ^{\kappa }=p<p(1+\nu )=p_{2}<\beta ^{\kappa +1}.
\end{equation*} 
Observe that%
\begin{equation*}
-\frac{n\delta }{p_{1}}<-\frac{n\delta }{p_{2}}<\alpha <n-\frac{n}{p_{1}}<n-\frac{n}{p_{2}},
\end{equation*}%
by \eqref{reserve}, \eqref{v}  and the fact that $p_{1}<p<p_{2}$%
. Let $0<\sigma <1$ be such that 
\begin{equation*}
\frac{1}{p}=\frac{1-\sigma }{p_{1}}+\frac{\sigma }{p_{2}}.
\end{equation*}%
Let $ \tau >0$  be such that
 \begin{equation*} \tau <\min\Big(\alpha+\frac{n\delta}{p_{2}},\frac{n-\frac{n}{p_{1}}-\alpha}{\frac{1}{\sigma}-1}\Big).
 \end{equation*} 
 We put $\alpha _{1}=\alpha-\tau$ and $\alpha _{2}=\alpha+\frac{(1-\sigma)\tau}{\sigma}$.
 Then 
 \begin{equation*}
 -\frac{n\delta }{p_{1}}<-\frac{n\delta }{p_{2}}<\alpha_{i} <n-\frac{n}{p_{1}}<n-\frac{n}{p_{2}},\quad i=1,2
 \end{equation*}
 and
\begin{equation}
\alpha =(1-\sigma )\alpha _{1}+\sigma \alpha _{2}\label{alpha3},\quad \alpha _{1}\neq\alpha _{2}.
\end{equation}
Substep 1.1 gives 
\begin{equation*}
\Big\|\Big(\sum\limits_{j=-\infty }^{\infty }t_{j}^{\beta }\big(\mathcal{M}%
(f_{j})\big)^{\beta }\Big)^{1/\beta }\Big\|_{\dot{K}_{p_{1}}^{\alpha
_{1},p_{1}}}\lesssim \Big\|\Big(\sum\limits_{j=-\infty }^{\infty
}t_{j}^{\beta }\left\vert f_{j}\right\vert ^{\beta }\Big)^{1/\beta }\Big\|%
_{\dot{K}_{p_{1}}^{\alpha _{1},p_{1}}}
\end{equation*}%
and 
\begin{equation*}
\Big\|\Big(\sum\limits_{j=-\infty }^{\infty }t_{j}^{\beta }\big(\mathcal{M}%
(f_{j})\big)^{\beta }\Big)^{1/\beta }\Big\|_{\dot{K}_{p_{2}}^{\alpha
_{1},p_{2}}}\lesssim \Big\|\Big(\sum\limits_{j=-\infty }^{\infty
}t_{j}^{\beta }\left\vert f_{j}\right\vert ^{\beta }\Big)^{1/\beta }\Big\|%
_{\dot{K}_{p_{2}}^{\alpha _{1},p_{2}}}.
\end{equation*}%
By Lemma 2.13 in \cite{D20.1}, we obtain%
\begin{equation*}
\Big\|\Big(\sum\limits_{j=-\infty }^{\infty }t_{j}^{\beta }\big(\mathcal{M}%
(f_{j})\big)^{\beta }\Big)^{1/\beta }\Big\|_{\dot{K}_{p}^{\alpha _{1},p}}\lesssim %
\Big\|\Big(\sum\limits_{j=-\infty }^{\infty }t_{j}^{\beta }\left\vert
f_{j}\right\vert ^{\beta }\Big)^{1/\beta }\Big\|_{\dot{K}_{p}^{\alpha _{1},p}}.
\end{equation*}%
Similary, we get 
\begin{equation*}
\Big\|\Big(\sum\limits_{k=-\infty }^{\infty }t_{k}^{\beta }\big(\mathcal{M}%
(f_{k})\big)^{\beta }\Big)^{1/\beta }\Big\|_{\dot{K}_{p}^{\alpha _{2},p}}\lesssim %
\Big\|\Big(\sum\limits_{k=-\infty }^{\infty }t_{k}^{\beta }\left\vert
f_{k}\right\vert ^{\beta }\Big)^{1/\beta }\Big\|_{\dot{K}_{p}^{\alpha _{2},p}}.
\end{equation*}%
Consequently
\begin{equation*}
\Big\|\{t_{k}\mathcal{M}%
(f_{k})\}_{k\in \mathbb{Z}}\Big\|_{\dot{K}^{\alpha _{i},p}(L_{p}(\ell _{q}))}\lesssim \Big\|\{t_{k}
f_{k}\}_{k\in \mathbb{Z}}\Big\|_{\dot{K}^{\alpha _{i},p}(L_{p}(\ell _{q}))}
\end{equation*}%
for any $i\in\{1,2\}$. By Theorem \ref{Lech}, we get 
\begin{equation*}
\Big\|\{t_{k}\mathcal{M}%
(f_{k})\}_{k\in \mathbb{Z}}\Big\|_{(\dot{K}^{\alpha _{1},p}(L_{p}(\ell _{q})),\dot{K}^{\alpha _{2},p}(L_{p}(\ell _{q}))_{\sigma ,q}}\lesssim \Big\|\{t_{k}
f_{k}\}_{k\in \mathbb{Z}}\Big\|_{(\dot{K}^{\alpha _{1},p}(L_{p}(\ell _{q})),\dot{K}^{\alpha _{2},p}(L_{p}(\ell _{q}))_{\sigma ,q}}.
\end{equation*}%
Using Theorem 	\ref{InterpolationHerz}, we obtain
\begin{equation*}
\Big\|\Big(\sum\limits_{k=-\infty }^{\infty }t_{k}^{\beta }\big(\mathcal{M}%
(f_{k})\big)^{\beta }\Big)^{1/\beta }\Big\|_{\dot{K}_{p}^{\alpha ,q}}\lesssim %
\Big\|\Big(\sum\limits_{k=-\infty }^{\infty }t_{k}^{\beta }\left\vert
f_{k}\right\vert ^{\beta }\Big)^{1/\beta }\Big\|_{\dot{K}_{p}^{\alpha ,q}}.
\end{equation*}

\textit{Step 2.} We shall prove $\mathrm{\eqref{key-est}}$ with $1<p<\beta<\infty $.

\textit{Substep 2.1}. $-\frac{n\delta }{p}<\alpha \leq0$. 
Let $\varrho>0$ be such that
\begin{equation}
1<\varrho <\min \Big(\frac{pn\delta}{n\delta+p\alpha},q,\beta,\theta \Big),\label{varrho1}
\end{equation}
Let $ \delta_{\{t_{k}^{-\varrho (\frac{p}{\varrho} )^{\prime }}\}_{k\in\mathbb{{Z}}}} $ be obtained by applying  Lemma \ref{Ap-Property}/(vi) to the sequence $\{t_{k}^{-\varrho (\frac{p}{\varrho} )^{\prime }}\}_{k\in\mathbb{{Z}}}$.  The
left-hand side of $\mathrm{\eqref{key-est}}$, raised to the\ power $\varrho $
is just 
\begin{equation*}
\sup \sum\limits_{j=-\infty }^{\infty }\int_{\mathbb{R}^{n}}t_{j}^{\varrho
}(x)\left( \mathcal{M}(f_{j})(x)\right) ^{\varrho }|g_{j}(x)|dx=\sup V,
\end{equation*}%
where the supremum is taken over all sequence of functions $\{g_{j}\}_{j\in 
	\mathbb{Z}}$ such that 
\begin{equation}
\{g_{j}\}_{j\in \mathbb{Z}}\in \dot{K}_{(\frac{p}{\varrho} )^{\prime }}^{-\alpha \varrho
	,(\frac{q}{\varrho} )^{\prime }}(\ell _{(\frac{\beta}{\varrho} )^{\prime }}),\quad \big\|%
\{g_{j}\}_{j\in \mathbb{Z}}\big\|_{\dot{K}_{(\frac{p}{\varrho} )^{\prime }}^{-\alpha \varrho
		,(\frac{q}{\varrho} )^{\prime }}(\ell _{(\frac{\beta}{\varrho} )^{\prime }})}\leq 1.
\label{last}
\end{equation}%
By Lemma \ref{FS-lemma} and H\"{o}lder's inequality, $V$ is bounded by 
\begin{align*}
& \ c\int_{\mathbb{R}^{n}}\sum\limits_{j=-\infty }^{\infty }\big|f_{j}(x)%
\big|^{\varrho }\mathcal{M}(t_{j}^{\varrho }g_{j})(x)dx \\
& \lesssim \Big\|\Big(\sum\limits_{j=-\infty }^{\infty }t_{j}^{\beta
}|f_{j}|^{\beta }\Big)^{1/\beta }\Big\|_{\dot{K}_{p}^{\alpha ,q}}^{\varrho } 
\Big\|\Big(\sum\limits_{j=-\infty }^{\infty }t_{j}^{-\varrho (\frac{\beta}{\varrho}
	)^{\prime }}(\mathcal{M}(t_{j}^{\varrho }g_{j}))^{(\frac{\beta}{\varrho} 
	)^{\prime }}\Big)^{1/(\frac{\beta}{\varrho} )^{\prime }}\Big\|_{\dot{K}_{(\frac{p}{\varrho}
		)^{\prime }}^{-\alpha \varrho ,(\frac{q}{\varrho} )^{\prime }}}.
\end{align*}%
By Step 1, the second term is bounded because of $(\frac{\beta}{\varrho} )^{\prime
}<(\frac{p}{\varrho} )^{\prime }$, 
\begin{equation*}
t_{k}^{-\varrho (\frac{p}{\varrho} )^{\prime }}\in A_{(\frac{p}{\varrho} )^{\prime }}(%
\mathbb{R}^{n}),\quad k\in \mathbb{Z},
\end{equation*}%
and from Lemma \ref{Ap-Property}/(iv) there exists a $1<\mu <(\frac{p}{\varrho}
)^{\prime }<\infty $ such that 
\begin{equation*}
t_{k}^{-\varrho (\frac{p}{\varrho} )^{\prime }}\in A_{\frac{(\frac{p}{\varrho} )^{\prime }}{\mu}
}(\mathbb{R}^{n}),\quad k\in \mathbb{Z}.
\end{equation*}%
Since $t_k^{-\varrho(\frac{p}{\varrho})'}\in A_{(\frac{p}{\varrho})'}(\mathbb{R}^n)$, $k\in\mathbb{Z}$, there exist constants $C>0$ and $\omega>0$, depending only on $(\frac{p}{\varrho})'$ and the $A_{(\frac{p}{\varrho})'}(\mathbb{R}^n)$ constant of $t_k^{-\varrho(\frac{p}{\varrho})'}$, such that, for every cube $Q$,
\[
M_{Q,1+\omega}\big(t_k^{-\varrho(\frac{p}{\varrho})'}\big)
\leq C M_Q\big(t_k^{-\varrho(\frac{p}{\varrho})'}\big).
\]

Let $\varepsilon$ be as in Theorem 2.6, corresponding to
$t_k^{p}\in A_{p}(\mathbb{R}^n),k\in \mathbb{Z}$.
In view of Theorem 7.2.2 in \cite{L. Graf14}, we may choose $\varepsilon$ sufficiently small. In particular, we choose $\varepsilon<\omega$.

Observing that
\[
\delta_{\{t_k^{-\varrho(\frac{p}{\varrho})'}\}_{k\in\mathbb Z}}
=1-\frac{1}{1+\omega}
\]
and
\[
\delta=1-\frac{1}{1+\varepsilon},
\]
we obtain
\[
\delta<\delta_{\{t_k^{-\varrho(\frac{p}{\varrho})'}\}_{k\in\mathbb Z}}.
\]

Consequently, by \eqref{varrho1}, we have
\[
-\frac{n\delta_{\{t_k^{-\varrho(\frac{p}{\varrho})'}\}_{k\in\mathbb Z}}}{(\frac{p}{\varrho})'}
<
-\frac{n\delta}{(\frac{p}{\varrho})'}
<
-\alpha\rho
<
n-\frac{n}{(\frac{p}{\varrho})'}.
\]

Thus, the desired estimate follows by Step 1 and $\mathrm{\eqref{last}}$.

\textit{Substep 2.2}. $0<\alpha <n-\frac{n}{p}$ . Since $t_{j}^{p}\in A_{p}(\mathbb{R}^{n}),j\in \mathbb{Z}$, there exist
constants $C_{1}>0$ and $\varepsilon _{1}>0,$\ depending only on $p$ and the 
$A_{p}(\mathbb{R}^{n})$-constant of $t_{j}^{p},j\in \mathbb{Z}$, such that
for every cube $Q$, 
\begin{equation}
M_{Q,1+\varepsilon _{1}}(t_{j}^{p})\leq C_{1}M_{Q}(t_{j}^{p}),\quad j\in 
\mathbb{Z} \label{ap-constant}
\end{equation}%
and $t_{j}^{p(1+\varepsilon _{1})}\in A_{p}(\mathbb{R}^{n}),j\in \mathbb{Z}$. Again, since $t_{j}^{p(1+\varepsilon _{1})}\in A_{p}(\mathbb{R}^{n}),j\in 
\mathbb{Z}$ there exist constants $C_{2}>0$ and $\varepsilon _{2}>0,$\
depending only on $p$ and the $A_{p}(\mathbb{R}^{n})$-constant of $%
t_{j}^{p},j\in \mathbb{Z}$, such that for every cube $Q$, 
\begin{equation*}
M_{Q,1+\varepsilon _{2}}(t_{j}^{p(1+\varepsilon _{1})})\leq
C_{2}M_{Q}(t_{j}^{p(1+\varepsilon _{1})}),\quad j\in \mathbb{Z}
\end{equation*}%
and $t_{j}^{p(1+\varepsilon _{1})(1+\varepsilon _{2})}\in A_{p}(\mathbb{R}%
^{n}),j\in \mathbb{Z}$. This 
leads to%
\begin{equation*}
(M_{Q}(t_{j}^{p(1+\varepsilon _{1})(1+\varepsilon _{2})}))^{1/(1+\varepsilon
	_{1})(1+\varepsilon _{2})}\leq CM_{Q,1+\varepsilon _{1}}(t_{j}^{p}),\quad j\in \mathbb{Z}.
\end{equation*}%
Repeating this argument, there exist constants $C_{m}>0$ and $\varepsilon
_{m}>0,m\in \mathbb{N}$\ depending only on $p$ and the $A_{p}(\mathbb{R}^{n})
$-constant of $t_{j}^{p},j\in \mathbb{Z}$, such that for every cube $Q$, 
\begin{equation*}
M_{Q,1+\varepsilon _{l}}(t_{j}^{p\Pi _{m=1}^{l-1}(1+\varepsilon _{m})})\leq
CM_{Q}(t_{j}^{p\Pi _{m=1}^{l-1}(1+\varepsilon _{m})}),\quad j\in \mathbb{Z}%
,l\in \mathbb{N}
\end{equation*}%
and $t_{j}^{p\Pi _{m=1}^{l}(1+\varepsilon _{m})}\in A_{p}(\mathbb{R}%
^{n}),j\in \mathbb{Z}$. This leads to%
\begin{equation*}
(M_{Q}(t_{j}^{p\Pi _{m=1}^{l}(1+\varepsilon _{m})}))^{1/\Pi
	_{m=1}^{l}(1+\varepsilon _{m})}\leq CM_{Q,1+\varepsilon _{1}}(t_{j}^{p}),\quad j\in \mathbb{Z}.
\end{equation*}%
We choose $0<\varepsilon <\min (\varepsilon _{1},...,\varepsilon _{l},\frac{\beta}{%
	p }-1,\theta -1).$ Then together with \eqref{ap-constant} we obtain
\begin{equation*}
(M_{Q}(t_{j}^{p(1+\varepsilon )^{l}}))^{1/(1+\varepsilon )^{l}}\leq
CM_{Q}(t_{j}^{p}),\quad j\in \mathbb{Z}
\end{equation*}%
and $t_{j}^{p(1+\varepsilon )^{l}}\in A_{p}(\mathbb{R}^{n}),j\in \mathbb{Z}$%
. Assume that 
\begin{equation*}
l-1\leq \frac{\log \beta -\log p}{\log (1+\varepsilon )}<l.
\end{equation*}%
Hence $p(1+\varepsilon )^{l-1}\leq \beta <p(1+\varepsilon )^{l}$. Now, $\mathrm{\eqref{key-est}}$ follows by  repeating
the same arguments of Substep 1.1.1 but with $\tau =p(1+\varepsilon )^{l}$
and $s=\frac{p}{1+\varepsilon }$.
The proof is complete.
\end{proof}

\begin{rem}
\label{r-estimates}$\mathrm{(i)}$ We would like to mention that the result
of this lemma is true if we assume that $t_{j}\in A_{\frac{p}{\theta} }(\mathbb{R}%
^{n})$,\ $j\in \mathbb{Z}$, $1<p<\infty $ with 
\begin{equation*}
A_{\frac{p}{\theta}}(t_{j})\leq c,\quad j\in \mathbb{Z},
\end{equation*}%
where $c>0$ independent of $j$.$\newline
\mathrm{(ii)}$ The proof of Lemma \ref{key-estimate1}\ for $t_{j}^{p}=\omega 
$, $j\in \mathbb{Z}$ is given in \cite{TDachun} but with  $-\frac{%
	n\delta }{p}<\alpha <n-\frac{n}{\theta}.$ $\newline
\mathrm{(iii)}$ In view of Lemma \ref{Ap-Property}/(iv) we can assume that $%
t_{j}^{p}\in A_{p}(\mathbb{R}^{n})$,\ $j\in \mathbb{Z}$, $1<p<\infty $ with 
\begin{equation*}
A_{p}(t_{j}^{p})\leq c,\quad j\in \mathbb{Z},
\end{equation*}%
where $c>0$ independent of $j$.
\end{rem}

\begin{rem}
	We retain the same notation as that used in Substep 2.1. Under some additional assumptions, we obtain $\omega>\varepsilon$. Indeed,  by the property (4) and (6) of Proposition 7.1.5 in  \cite{L. Graf14} we obtain
	\begin{equation*}
	 A_{(\frac{p}{\varrho})'}(t_k^{-\varrho(\frac{p}{\varrho})'})=\big(A_{\frac{p}{\varrho}}(t_k^{p})\big)^{^{\frac{1}{\frac{p}{\varrho}-1}}}\leq \big(A_{\frac{p}{\theta}}(t_k^{p})\big)^{^{\frac{1}{\frac{p}{\varrho}-1}}}.
	 \end{equation*}
	 In view of Theorem 7.2.2 in  \cite{L. Graf14} we have $\omega=\frac{1}{2}\frac{-\log h}{\log\frac{2^{n}}{s}}$, where $0<s<1$ and
	 \begin{equation*}
	  h=1-\frac{(1-s)^{(\frac{p}{\varrho})'}}{A_{(\frac{p}{\varrho})'}(t_k^{-\varrho(\frac{p}{\varrho})'})}.
	  \end{equation*}
	  Observing that $A_{\frac{p}{\theta}}(t_k^{p})\geq 1 $. Assume in addition that $p>1+\theta$ and  $\varrho<\min(\frac{p}{2},p-\theta) $.  Then
	  \begin{equation*}
	  A_{(\frac{p}{\varrho})'}(t_k^{-\varrho(\frac{p}{\varrho})'})\leq A_{\frac{p}{\theta}}(t_k^{p}).
	  \end{equation*}
	  This leads to
	  \begin{equation*}
	  h<1-\frac{(1-s)^{\frac{p}{\theta}}}{A_{\frac{p}{\theta}}(t_k^{p})}=u.
	  \end{equation*}
	  Therefore
	   \begin{equation*}
	   \omega>\frac{1}{2}\frac{-\log u}{\log\frac{2^{n}}{s}}=\varepsilon.
	   \end{equation*}
	   
\end{rem}
\begin{lem}
	\label{maximalHerz}  Let $1<\theta \leq p<\infty $ and $0< q\leq \infty $. Let $\{t_{j}\}_{j\in \mathbb{Z}}$\ be a $p$-admissible\ weight\ sequence\ such that $%
	t_{j}^{p}\in A_{\frac{p}{\theta }}(\mathbb{R}^{n})$, $j\in \mathbb{Z}$.
	Assume that $t_{j}^{p}$,\ $j\in \mathbb{Z}$ have the same Muckenhoupt
	constant.  Let $\delta $ be as in Lemma \ref{Ap-Property} and  $-\frac{%
		n\delta }{p}<\alpha <n-\frac{n}{p}$.
	\newline
	$\mathrm{\bullet }$ If $t_{j}^{p}$,\ $j\in \mathbb{Z}$ satisfies $\mathrm{\ %
		\eqref{Asum1}}$ with $\sigma _{1}=\theta\left( \frac{p}{\theta}\right) ^{\prime }$, then 
	\begin{equation}
	\big\|t_{k}\mathcal{M(}f_{j})\big\|_{\dot{K}_{p}^{\alpha ,q}}\leq c\text{ }%
	2^{\alpha _{1}(k-j)}\big\|t_{j}f_{j}\big\|_{\dot{K}_{p}^{\alpha ,q}}  \label{herz1}
	\end{equation}%
	holds for all sequences of functions $f_{j}$, 
	$j\in \mathbb{Z}$  such that $t_{j}f_{j}\in$ $\dot{K}_{p}^{\alpha ,q}$ and $j\geq k$, where $c>0$ is independent of $k$ and $j$. $\newline
	\mathrm{\bullet }$ If $t_{j}^{p}$,\ $j\in \mathbb{Z}$ satisfies $\mathrm{\ %
		\eqref{Asum2}}$ with $\sigma _{2}\geq p$, then 
	\begin{equation}
	\big\|t_{k}\mathcal{M(}f_{j})\big\|_{\dot{K}_{p}^{\alpha ,q}}\leq c\text{ }%
	2^{\alpha _{2}(k-j)}\big\|t_{j}f_{j}\big\|_{\dot{K}_{p}^{\alpha ,q}} \label{herz2}
	\end{equation}%
	holds for all sequences of functions $f_{j}$, 
	$j\in \mathbb{Z}$  such that $t_{j}f_{j}\in$ $\dot{K}_{p}^{\alpha ,q}$ and $k\geq j $, where $c>0$ is independent of $k$ and $j$.
\end{lem}

\begin{proof} We proceed into two steps.
	
\textit{Step 1}. We will prove \eqref{herz1}.	We split the operator $\mathcal{M}(f_{j}),j\in 
	\mathbb{Z}$ into 
	\begin{equation}
	\mathcal{M}(f_{j})\leq \mathcal{M}(f_{j}\chi _{B_{m-2}})+\mathcal{M}%
	(f_{j}\chi _{\widetilde{R}_{m}})+\mathcal{M}(f_{j}\chi _{{\mathbb{R}^{n}}%
		\setminus B_{m+2}}),  \label{splitMherz}
	\end{equation}%
	where $\widetilde{R}_{m}=\left\{ x\in {\mathbb{R}^{n}}:\,2^{m-2}\leq
	|x|<2^{m+2}\right\} $ with $j,m\in \mathbb{Z}$.
	
	\emph{Estimation of }$\mathcal{M}$\emph{\emph{$(f_{j}$}$\chi _{B_{m-2}})$}. Observing that
	
	\begin{equation*}
	\Big\| t_{k} \mathcal{M}(f_{j}\chi _{B_{m-2}})\chi _{m}\Big\|_{p},\quad m\in \mathbb{Z%
	}
	\end{equation*}
	can be estimated from above by 
	\begin{equation*}
	\sup \int_{R_{m}}t_{k}(x)\mathcal{M}(f_{j}\chi _{B_{m-2}})(x)|g(x)|dx=\sup T_{j,m,k},
	\end{equation*}%
	where the supremum is taken over all sequence of functions $g\chi_{m}\in L_{p^{\prime }}(\mathbb{R}^{n})$ with 
	\begin{equation}
	\big\|g\chi_{m}\big\|_{p^{\prime }}\leq 1.\label{dualherz}
	\end{equation}%
	 For $%
	x\in R_{m}$ and $y\in B_{m-2}$ we have $|x-y|\geq |x|-|y|>2^{m-2}$. As in Lemma \ref{key-estimate1} we obtain
	\begin{equation}
	2^{m\alpha }T_{j,m,k}\lesssim \sum_{i=-\infty }^{m-2}2^{(m-i)\alpha
		-mn}\int_{R_{m}}t_{k}(x)|g_{j}(x)|dx\int_{R_{i}}2^{i\alpha }|f_{j}(y)|\,dy. \label{Herz4}
	\end{equation}%
	and 
	\begin{align}
	\int_{R_{m}}t_{k}(x)|g(x)|dx& \leq |R_{m}|M_{B_{m},\tau }(t_{k})M_{B_{m},\tau ^{\prime }}(g\chi _{m})
	\notag
	\end{align}%
	with $\tau =p(1+\varepsilon )$ where $\varepsilon $ as in Theorem \ref%
	{reverse Holder inequality}, which is possible since $t_{j}^{p}\in
	A_{\frac{p}{\theta }}(\mathbb{R}^{n})$,  $j \in \mathbb{Z}$. Also, we have 
	\begin{equation*}
	M_{B_{m},\tau }(t_{k})\lesssim 
	M_{B_{m},p}(t_{k})
	\lesssim 2^{(i-m)(\frac{n}{p}-n)}M_{B_{i},p}(t_{k}),\quad i\leq m-2,m,k\in \mathbb{Z}.
	\end{equation*}%
	Therefore, 
	\begin{equation*}
	\int_{R_{m}}t_{k}(x)|g(x)|dx\lesssim 2^{(i-m)(\frac{n}{p}%
		-n)}|R_{m}|M_{B_{i},p}(t_{k})M_{B_{m},\tau ^{\prime }}(g\chi _{m}).
	\end{equation*}%
	By H\"{o}lder's inequality and $\mathrm{\ %
		\eqref{Asum1}}$ with $\sigma _{1}=\theta\left( \frac{p}{\theta}\right) ^{\prime }$, we obtain 
	\begin{equation*}
	M_{R_{i}}(f_{j})\leq M_{B_{i},\sigma_{1} }(t_{j}^{-1})M_{B_{i},\frac{p}{\theta}}(t_{j}f_{j}\chi _{i})\lesssim 2^{\alpha _{1}(k-j)}(M_{B_{i},p}(t_{k}))^{-1}M_{B_{i},\frac{p}{\theta}}(t_{j}f_{j}\chi _{i}) .
	\end{equation*}%
	Hence 
	\begin{equation*}
	\int_{R_{i}}|f_{j}(y)|\,\chi
	_{i}(y)dy\int_{R_{m}}t_{k}(x)|g(x)|dx,
	\end{equation*}%
 is bounded above by
  \begin{equation*}
 c  2^{\alpha _{1}(k-j)+(i-m)(\frac{n}{p}%
 	-n)}|R_{m}||R_{i}|M_{B_{i},\frac{p}{\theta}}(t_{j}f_{j}\chi _{i})M_{B_{m},\tau ^{\prime
 }}(g\chi _{m}).
 \end{equation*}
	The last expression is bounded by 
	\begin{align}
	&c2^{\alpha _{1}(k-j)+(i-m)(\frac{n}{p}%
		-n)}|R_{m}|\int_{R_{i}}M_{B_{i},\frac{p}{\theta}}(t_{j}f_{j}\chi
	_{i})M_{B_{m},\tau ^{\prime }}(g\chi _{m})dx\notag \\ 
	&\lesssim 2^{\alpha _{1}(k-j)}2^{(i-m)(\frac{n}{p}-n)}|R_{m}|\int_{%
		\mathbb{R}^{n}}\mathcal{M}_{\frac{p}{\theta}}(t_{j}f_{j})(x)\mathcal{M}_{\tau ^{\prime
	}}(g\chi _{m})(x)dx.\label{Herz6}
	\end{align}%
	Substituting \eqref{Herz6} into \eqref{Herz4}, we obtain
	\begin{equation*}
	2^{m\alpha }T_{j,m,k}\lesssim 2^{\alpha _{1}(k-j)}
	\sum_{i=-\infty }^{m-2}2^{(m-i)(\alpha -\frac{n}{p}+n)}\int_{\mathbb{R}%
		^{n}}\mathcal{M}_{\frac{p}{\theta}}(2^{\alpha
		i}t_{j}f_{j}\chi _{i})(x)\mathcal{M}_{\tau ^{\prime }}(g\chi _{m})(x)dx.
	\end{equation*}%
	By  H\"{o}lder's inequality $2^{m\alpha
	}T_{j,m,k}$ can be estimated from above by 
	\begin{align*}
	& c2^{\alpha _{1}(k-j)}\sum_{i=-\infty }^{m-2}2^{(m-i)(\alpha -\frac{n}{p}+n)}2^{\alpha i}\big\|%
	\mathcal{M}_{\frac{p}{\theta}}\big(%
	t_{j}f_{j}\chi _{i}\big)\big\|_{p}\big\|\mathcal{M}_{\tau
		^{\prime }}(g\chi _{m})%
	\big\|_{p^{\prime }} \\
	& \lesssim 2^{\alpha _{1}(k-j)}\sum_{i=-\infty }^{m-2}2^{(m-i)(\alpha -\frac{n}{p}+n)}2^{\alpha
		i}\big\|%
	\mathcal{M}_{\frac{p}{\theta}}\big(%
t_{j}f_{j}\chi _{i}\big)\big\|_{p} \\
	& \lesssim 2^{\alpha _{1}(k-j)}\sum_{i=-\infty }^{m-2}2^{(m-i)(\alpha -\frac{n}{p}+n)}2^{\alpha
		i}\big\|%
t_{j}f_{j}\chi _{i}\big\|_{p},
	\end{align*}%
	where in the second and third inequality we have used the boundedness of maximal function in $L^{p}(\mathbb R^{n})$ spaces, since $\tau^{\prime}<p$ and $ \frac{p}{\theta}<p$. Applying Lemma \ref%
	{lem:lq-inequality} we obtain
	\begin{equation}
	\big\|\big\{2^{m\alpha }\sup T_{j,m,k}\big\}%
	_{m\in \mathbb{Z}}\big\|_{\ell _{q}}\lesssim  2^{\alpha _{1}(k-j)} \big\|t_{j}f_{j}\big\|%
	_{\dot{K}_{p}^{\alpha ,q}},
	\end{equation}%
		where the supremum is taken over all sequence of functions $\{g\chi_{m}\}\in L_{p^{\prime }}(\mathbb{R}^{n})$ such that \eqref{dualherz} holds.
	
	\emph{Estimation of }$\mathcal{M}$\emph{$(f_{j}\chi _{\widetilde{R}_{m}})$}. In view of Remark \ref{r-estimatescopy(2)}  we have 
	\begin{align*}
	\Big(\sum\limits_{m=\mathbb{-\infty }}^{\infty }2^{m{\alpha q}}\big\|t_{k}%
	\mathcal{M}(f_{j}\chi _{\widetilde{R}_{m}})\big\|_{p }^{q}\Big)%
	^{1/q}\lesssim & 2^{\alpha _{1}(k-j)}\Big(\sum\limits_{m=\mathbb{-\infty }}^{\infty }2^{m{\alpha
			q}}\big\|t_{j}f_{j}\chi _{\widetilde{R}_{m}}\big\|_{p}^{q}\Big)^{1/q} \\
	\lesssim & 2^{\alpha _{1}(k-j)}\big\|t_{j}f_{j}\big\|_{\dot{K}_{p}^{\alpha ,q}},
	\end{align*}%
	with the usual modification when $q=\infty .$
	
	\emph{Estimation of }$\mathcal{M}$\emph{$(f_{j}\chi _{{\mathbb{R}^{n}}\setminus
			B_{m+2}})$}. Noting that $|x-y|>2^{i-3}$ for $x\in R_{m}$ and $y\in R_{i}$,
	with $i\geq m+3$. This implies that 
	\begin{align*}
	2^{m\alpha }\big\|t_{k}\mathcal{M}(f_{j}\chi _{{\mathbb{R}^{n}}\setminus B_{m+2}})\chi _{m}\big\|_{p }&
	\lesssim \big\|t_{k}\chi _{m}\big\|_{p }\sum_{i=m+3}^{\infty }2^{(m-i)\alpha -in}2^{i\alpha
	}\int_{R_{i}}|f_{j}(y)|\,dy \\
	& \lesssim \sum_{i=m+3}^{\infty }2^{(m-i)\alpha -in}2^{i\alpha }\big\|t_{j}%
	f_{j}\,\chi _{i}\big\|_{p }\big\|t_{k}\chi _{m}\big\|_{p} \big\|t^{-1}_{j}\chi _{i}\big\|_{p^{\prime}},
	\end{align*}%
where we have used the H\"{o}lder inequality.	Since $t_{j}^{p }\in A_{p}(\mathbb{R}^{n})$, $j\in \mathbb{Z}$,
	Lemma \ref{Ap-Property}/(vi) gives%
	\begin{equation*}
	\big\|t_{k}\chi _{B_{m}}\big\|_{p} \lesssim 2^{\frac{%
			(m-i)n\delta }{p }}	\big\|t_{k}\chi _{B_{i}}\big\|_{p },\quad i \geq m+3.
	\end{equation*}%
	In addition  by H\"{o}lder's inequality we  obtain
\begin{equation*}
2^{-i\frac{n}{p^{\prime}}}\big\|t^{-1}_{j}\chi _{i}\big\|_{p^{\prime} }\lesssim 2^{-i\frac{n}{\sigma_{1}}}\big\|t^{-1}_{j}\chi _{i}\big\|_{\sigma_{1}}.
\end{equation*}
Multiplying both sides by $ 2^{-i\frac{n}{p}} 	\big\|t_{k}\chi _{B_{i}}\big\|_{p }$  and applying $\mathrm{\eqref{Asum1}}$ with $\sigma _{1}=\theta\left( \frac{p}{\theta}\right) ^{\prime }$, we obtain
\begin{equation*}
	2^{ -in}\big\|t_{k}\chi _{m}\big\|_{p} \big\|t^{-1}_{j}\chi _{i}\big\|_{p^{\prime} }\lesssim 2^{\frac{(m-i)n\delta }{p }+\alpha _{1}(k-j)}.
	\end{equation*}
 Lemma \ref%
{lem:lq-inequality} gives
\begin{align*}
\Big(\sum\limits_{m=\mathbb{-\infty }}^{\infty }2^{m{\alpha q}}\big\|
t_{k}\mathcal{M}(f_{j}\chi _{{\mathbb{R}^{n}}\setminus B_{m+2}})\chi _{m}\big\|_{p}^{q}\Big)%
^{1/q}
\lesssim & 2^{\alpha _{1}(k-j)}\big\|t_{j}f_{j}\big\|_{\dot{K}_{p}^{\alpha ,q}},
\end{align*}
with the usual modification when $q=\infty .$.

\textit{Step 2}. We will prove \eqref{herz2}. From Lemma \ref{Ap-Property}/(iii)  and  $\mathrm{
	\eqref{Asum2}}$ with $\sigma _{2}\geq p$, we obtain%
\begin{equation*}
M_{B_{m},p}(t_{k})\leq C2^{(i-m)(\frac{n}{p}-n)}M_{B_{i},p}(t_{k})\leq C2^{(i-m)(\frac{n}{p}-n)+\alpha _{2}(k-j)}M_{B_{i},p}(t_{j}),\quad i\leq m-2.
\end{equation*}%
Now, using   $\mathrm{
	\eqref{Asum2}}$ with $\sigma _{2}\geq p$, we obtain
	\begin{equation*}
\big\|t_{k}\chi _{B_{m}}\big\|_{p} \lesssim 2^{\frac{%
		(m-i)n\delta }{p }}	\big\|t_{k}\chi _{B_{i}}\big\|_{p}\lesssim 2^{\frac{%
		(m-i)n\delta }{p }+\alpha _{2}(k-j)}	\big\|t_{j}\chi _{B_{i}}\big\|_{p} ,\quad i \geq m+3.
\end{equation*}
Hence, \eqref{herz2} follows by the same arguments as in Step 1.
\end{proof}

	\begin{rem}
		\ In Lemma \ref{maximalHerz}, we can assume that $%
		t_{j}^{p}\in A_{p}(\mathbb{R}^{n})$,\ $j\in \mathbb{Z}$, $1<p<\infty $ with 
		\begin{equation*}
		A_{p}(t_{j}^{p})\leq c,\quad j\in \mathbb{Z},
		\end{equation*}%
		which follows by Lemma \ref{Ap-Property}/(iii).
		
\end{rem}
Next, we consider the case $\beta=\infty$ in Lemma \ref{key-estimate1}.
\begin{lem}
	\label{key-estimate3}Let $1<\theta \leq p<\infty $\ and $1\leq q\leq \infty $.
	Let $\{t_{k}\}_{k\in \mathbb{Z}}$\ be a $p$-admissible\ weight\ sequence\ such that $%
	t_{k}^{p}\in A_{\frac{p}{\theta }}(\mathbb{R}^{n})$, $k\in \mathbb{Z}$.
	Assume that $t_{k}^{p}$,\ $k\in \mathbb{Z}$ have the same Muckenhoupt
	constant, $A_{\frac{p}{\theta }}(t_{k})=c,k\in \mathbb{Z}$. Let $\delta $ be as in Lemma \ref{Ap-Property}. Let $-\frac{%
		n\delta }{p}<\alpha <n-\frac{n}{p}$. Then
	\begin{equation}
	\Big\|\sup_{k\in \mathbb{Z}}\big( t_{k}\mathcal{M}(f_{k})%
	\big)\Big\|_{\dot{K}_{p}^{\alpha ,q}}\lesssim \Big\|\sup_{k\in \mathbb{Z}}\big( t_{k}|f_{k}|)\Big\|_{\dot{K}_{p}^{\alpha ,q}} \label{key-est1}
	\end{equation}%
	holds for all sequences of functions $\{t_{k}f_{k}\}_{k\in \mathbb{Z}}\in
	\dot{K}_{p}^{\alpha ,q}(\ell _{\infty })$.
\end{lem}

\begin{proof}
	 We adopt the approach used in \cite{DRu25}. For  a sequence of functions $\{g_{k}\}_{k\in \mathbb{Z}}$  and $\alpha \in \mathbb{R}$ we define 
	\begin{equation}
	T(\{g_{k}\}_{k\in \mathbb{Z}})=\{\mathcal{M}(|x|^{\alpha}t_{k}\mathcal{M}(|x|^{-\alpha}t_{k}^{-1}g_{k}))\}_{k\in \mathbb{Z}}.
	\label{operator1}
	\end{equation}%
	$T$ is a sublinear operator.
	
	\textit{Step 1.} Let $\{g_{k}\}_{k\in \mathbb{Z}}\subset L_{p}(\ell
	_{\infty })$, $t_{k}^{p}\in A_{\frac{p}{\theta }}(\mathbb{R}^{n})$, $k\in \mathbb{Z}$ and $-\frac{%
		n\delta }{p}<\alpha <n-\frac{n}{p}$. We prove that $|x|^{\alpha}t_{k}\mathcal{M}(|x|^{-\alpha}t_{k}^{-1}g_{k})$, 
	$k\in \mathbb{Z}$\ are locally integrable\ over $\mathbb{R}^{n}$. Let $%
	K\subset \mathbb{R}^{n}$ be a compact set. Using H\"{o}lder's inequality and
	Lemma \ref{maximalHerz}, we obtain 
	\begin{align*}
	\int_{K}|x|^{\alpha}t_{k}(x)\mathcal{M}(|x|^{-\alpha}t_{k}^{-1}g_{k})(x)dx\leq & \big\|\chi
	_{K}\big\|_{p^{\prime }}\big\||x|^{\alpha}t_{k}\mathcal{M}%
	(|x|^{-\alpha}t_{k}^{-1}g_{k})\big\|_{p} \\
	\lesssim & \big\|g_{k}\big\|_{p}\\
	\lesssim & \big\|\sup_{j\in \mathbb{Z}}(|g_{j}|)%
	\big\|_{p},
	\end{align*}%
	where the implicit constant is independent of $k$. This yields the desired
	result. In addition it holds $|x|^{\alpha}t_{k}\mathcal{M}(|x|^{-\alpha}t_{k}^{-1}g_{k})\leq \mathcal{%
		M}(|x|^{\alpha}t_{k}\mathcal{M}(|x|^{-\alpha}t_{k}^{-1}g_{k}))$, $k\in \mathbb{Z}$, almost everywhere.

	\textit{Substep 1.1. }Proceeding as in Substep 1.2 of the proof of Lemma \ref{key-estimate1}, we obtain that there exists  $-\frac{%
		n\delta }{p_{1}}<\alpha <n-\frac{n}{p_{1}}$, $t_{k}^{p}\in A_{p_{1}}(\mathbb{R}^{n})$
	and $t_{k}^{p_{1}}\in A_{p_{1}}(\mathbb{R}^{n})$. Arguing as in the proof of Lemma 3 in  \cite{DRu25} and applying Lemma \ref{maximalHerz}, we conclude that   $T$ maps $L_{p_{1}}(\ell _{\infty })$ to $L_{p_{1},\infty }(\ell
	_{\infty })$.
	
	\textit{Substep 1.2. } Using again the fact that $t_{k}^{p}\in A_{p}(\mathbb{R}%
	^{n}) $ for all $k\in \mathbb{Z}$, and proceeding as in Substep 1.2 of the proof of Lemma \ref{key-estimate1}, we obtain that there exists  an $\varepsilon >0$ independent of  $k$, such that $%
	t_{k}^{p(1+\varepsilon )}\in A_{p}(\mathbb{R}^{n})\subset
	A_{p(1+\varepsilon )}(\mathbb{R}^{n}),k\in \mathbb{Z}$ and  $-\frac{%
		n\delta }{p(1+\varepsilon )}<\alpha <n-\frac{n}{p(1+\varepsilon )}$. Using again the same arguments of Lemma 3 in \cite{DRu25} with the help of Lemma \ref{maximalHerz} we obtain $T$ maps $L_{p(1+\varepsilon )}(\ell _{\infty })$ to $%
	L_{p(1+\varepsilon ),\infty }(\ell _{\infty })$.
	
	\textit{Substep 1.3.\ Conclusion. }By Substeps 1.1-1.2\ and the
	vector-valued version of the Marcinkiewicz interpolation; see  \cite[Lemma 2.5]{GLD} we obtain the boundedness of $T$ from $%
	L_{p}(\ell _{\infty })$ to $L_{p}(\ell _{\infty })$ and 
	\begin{align*}
	\big\|\sup_{k\in \mathbb{Z}}(t_{k}\mathcal{M}(|x|^{-\alpha}t_{k}^{-1}g_{k}))|x|^{\alpha}\big\|_{p}\lesssim & \big\|\sup_{k\in \mathbb{Z}}(\mathcal{M}%
	(|x|^{\alpha}t_{k}\mathcal{M}(|x|^{-\alpha}t_{k}^{-1}g_{k})))\big\|_{p} \\
	\lesssim & \big\|\sup_{k\in \mathbb{Z}}(|g_{k}|)%
	\big\|_{p},
	\end{align*}%
	whenever $\{g_{k}\}_{k\in \mathbb{Z}}\subset L_{p}(\ell _{\infty })$.\\
	 In the last two estimates setting %
	$g_{k}=|x|^{\alpha}t_{k}f_{k},k\in \mathbb{Z}$\ we then obtain
	\begin{align*}
	\big\||x|^{\alpha}\sup_{k\in \mathbb{Z}}(t_{k}\mathcal{M}(f_{k}))\big\|_{p}
	\lesssim & \big\||x|^{\alpha}\sup_{k\in \mathbb{Z}}(|t_{k}f_{k}|)%
	\big\|_{p},
	\end{align*}%
	whenever $\{t_{k}f_{k}\}_{k\in \mathbb{Z}}\subset L_{p}(\ell _{\infty },|x|^{\alpha})$.
	
		\textit{Step 2}. We will prove \eqref{key-est1}. Let $h\in \mathbb{R}$ be such that  $0<h<\min\Big(n-\frac{n}{p}-\alpha,\alpha+\dfrac{n\delta}{p}\Big)$. We put $\alpha _{1}=\alpha -h$ and $\alpha _{2}=\alpha+h$. Then
	
	\begin{equation*}
	-\frac{n\delta }{p}<\alpha _{i}<n-\frac{n}{p},\quad i\in \{1,2\}\label{alpha4}
	\end{equation*}
	and%
	\begin{equation*}
	\alpha =(1-\frac{1}{2} )\alpha _{1}+\frac{1}{2} \alpha _{2}\label{alpha5},\quad \alpha _{1}\neq\alpha _{2}.
	\end{equation*}
From Substep 1.3, we obtain
		\begin{align*}
	\big\|\sup_{k\in \mathbb{Z}}(t_{k}\mathcal{M}(f_{k}))\big\|_{\dot{K}_{p}^{\alpha_{i},p}}
	\lesssim & \big\|\sup_{k\in \mathbb{Z}}(|t_{k}f_{k}|)%
	\big\|_{\dot{K}_{p}^{\alpha_{i},p}}
	\end{align*}
	for any $\{t_{k}f_{k}\}_{k\in \mathbb{Z}}\in
	\dot{K}_{p}^{\alpha_{i} ,p}(\ell _{\infty }),i\in\{1.2\}$. As in Lemma \ref{key-estimate1}, Substep 1.2, by Theorems \ref{Lech} and  \ref{InterpolationHerz}, we obtain
	\begin{align*}
	\big\|\sup_{k\in \mathbb{Z}}(t_{k}\mathcal{M}(f_{k}))\big\|_{\dot{K}_{p}^{\alpha,q}}
	\lesssim & \big\|\sup_{k\in \mathbb{Z}}(|t_{k}f_{k}|)%
	\big\||_{\dot{K}_{p}^{\alpha,q}}
	\end{align*}
	for any $\{t_{k}f_{k}\}_{k\in \mathbb{Z}}\in
	\dot{K}_{p}^{\alpha ,q}(\ell _{\infty })$. This completes the proof of Lemma \ref{key-estimate3}.
\end{proof}
We shall need the following lemma concerning a discrete convolution inequality. The proof follows the same lines as that of  \cite{D20.1}.

\begin{lem}
\label{lq-inequality1}\textit{Let }$\alpha \in \mathbb{R},0<a<1,1\leq
p,q<\infty ,1< r\leq \infty $ \textit{and }$0<\beta <\infty $\textit{.
Let }$\left\{ f_{k}\right\}_{k\in \mathbb{Z}} $ \textit{\ and }$\left\{ g_{k}\right\}_{k\in \mathbb{Z}} $ 
\textit{be two sequences of positive\ real\ functions\ and denote} 
\begin{equation*}
\delta _{k}=\sum_{j=-\infty }^{k+m}a^{k-j}\big\|g_{k}f_{j}\big\|_{1}^{1/\beta },\quad k,m\in \mathbb{Z}
\end{equation*}%
and\textit{\ } 
\begin{equation*}
\eta _{k}=\sum_{j=k+m}^{\infty }a^{j-k}\big\|g_{k}f_{j}
\big\|_{1}^{1/\beta },\quad k,m\in \mathbb{Z}.
\end{equation*}%
Then there exists a constant $c>0\ $\textit{depending only on }$a$\textit{\
and }$q$ such that 
\begin{equation}
\sum\limits_{k=-\infty }^{\infty }\delta _{k}^{\beta
}+\sum\limits_{k=-\infty }^{\infty }\eta _{k}^{\beta }\leq c\Big\|\Big(%
\sum\limits_{k=-\infty }^{\infty }f_{k}^{r}\Big)^{1/r}\Big\|_{\dot{K}_{p}^{\alpha
,q}}\Big\|\Big(\sum\limits_{k=-\infty }^{\infty }g_{k}^{r^{\prime }}\Big)%
^{1/r^{\prime }}\Big\|_{\dot{K}_{p^{\prime }}^{-\alpha ,q^{\prime }}}.
\label{estimate1}
\end{equation}
\end{lem}

Using the same type of arguments as in Lemma \ref{lq-inequality1} combined
with \cite{D20.1} it is easy to prove the following lemma.

\begin{lem}
\label{lq-inequality1copy(1)}\textit{Let }$\alpha \in \mathbb{R}$, \textit{Let }$a>0,1\leq p,q< \infty  $ and  $1\leq r\leq \infty $ \textit{and }$0<\beta <\infty $\textit{. Let }$\left\{
f_{k}\right\}_{k\in \mathbb{Z}} $\textit{\ and }$\left\{ g_{k}\right\}_{k\in \mathbb{Z}} $ \textit{be two
sequences of positive\ real\ functions\ and denote} 
\begin{equation*}
\delta _{k}=\sum_{j=k}^{k+m}a^{k-j}\big\|g_{k}f_{j} %
\big\|_{1}^{1/\beta},\quad k\in \mathbb{Z},m\in \mathbb{N}
\end{equation*}
and\textit{\ } 
\begin{equation*}
\eta _{k}=\sum_{j=k+l}^{k}a^{j-k}\big\|g_{k}f_{j}%
\big\|_{1} ^{1/\beta},\quad k\in \mathbb{Z},l\leq 0.
\end{equation*}
Then there exists a constant $c>0\ $\textit{depending only on }$a,m,l$ 
\textit{\ and }$q$ such that \eqref{estimate1} holds.
\end{lem}

The following lemma plays an important role in the study of our function spaces. 

\begin{lem}
\label{key-estimate1.1herz}Let $m\in \mathbb{Z}$, $K\geq 0,1<\theta \leq
p<\infty ,1<q,\beta<\infty \ $ and $\upsilon=(\upsilon _{1},\upsilon _{2})\in \mathbb{R}%
^{2}$. Let $\{t_{k}\}_{k\in \mathbb{Z}}\in \dot{X}_{\upsilon ,\sigma ,p}$ be a $p$-admissible
weight sequence with $\sigma =(\sigma _{1}=\theta \left( \frac{p}{\theta} \right)
^{\prime },\sigma _{2}\geq p)$. Let $\delta $ be as in Lemma \ref{Ap-Property}. Let $-\frac{%
	n\delta }{p}<\alpha <n-\frac{n}{p}$. Then for all sequence of functions $%
\{t_{k}f_{k}\}_{k\in \mathbb{Z}}\in \dot{K}_{p}^{\alpha
	,q}(\ell _{\beta})$, 
\begin{align}
& \Big\|\Big(\sum\limits_{j=-\infty }^{\infty }t_{j}^{\beta}\Big(\sum_{l=-\infty
}^{j+m}2^{(l-j)K}\mathcal{M}(f_{l})\Big)^{\beta}\Big)^{1/\beta}
\Big\|_{\dot{K}_{p}^{\alpha,q}}  \notag \\
& \lesssim \Big\|\Big(\sum\limits_{j=-\infty }^{\infty }t_{j}^{\beta}\left\vert
f_{j}\right\vert ^{\beta}\Big)^{1/\beta}\Big\|_{\dot{K}_{p}^{\alpha
		,q}}
\label{key-est1.1.1}
\end{align}%
if $K>\alpha _{2}$ and 
\begin{align}
& \Big\|\Big(\sum\limits_{j=-\infty }^{\infty }t_{j}^{\beta}\Big(%
\sum_{l=j+m}^{\infty }2^{(l-j)K}\mathcal{M}(f_{l})\Big)^{\beta}\Big)^{1/\beta}\Big\|_{\dot{K}_{p}^{\alpha
		,q}}  \notag \\
& \lesssim \Big\|\Big(\sum\limits_{j=-\infty }^{\infty }t_{j}^{\beta}\left\vert
f_{j}\right\vert ^{\beta}\Big)^{1/\beta}\Big\|_{\dot{K}_{p}^{\alpha
		,q}}
\label{key-est1.1.2.1}
\end{align}%
if $K<\alpha _{1}$.
\end{lem}
\begin{proof} The proof is divided into two steps.
	
	\textit{Step 1. Proof of \eqref{key-est1.1.1}}. We split the operator $\mathcal{M}(f_{l}),l\in 
	\mathbb{Z}$ into 
	\begin{equation*}
	\mathcal{M}(f_{l})\leq \mathcal{M}(f_{l}\chi _{B_{k-2}})+\mathcal{M}%
	(f_{l}\chi _{\widetilde{R}_{k}})+\mathcal{M}(f_{l}\chi _{{\mathbb{R}^{n}}%
		\setminus B_{k+2}}),  \label{splitMherz1}
	\end{equation*}%
	where $\widetilde{R}_{k}=\left\{ x\in {\mathbb{R}^{n}}:\,2^{k-2}\leq
	|x|<2^{k+2}\right\} $ with $l,k\in \mathbb{Z}$. 
	
	\textit{Subtep 1.1.} By  Lemma 2.18 in \cite{D20.1} we obtain
	
	\begin{align*}
	& 2^{\alpha k}\Big\|\Big(\sum\limits_{j=-\infty }^{\infty }t_{j}^{\beta}\Big(\sum_{l=-\infty
	}^{j+m}2^{(l-j)K}\mathcal{M}(f_{l}\chi_{\widetilde{R}_{k}})\Big)^{\beta}\Big)^{1/\beta}\chi_{k}
	\Big\|_{p}  \notag \\
	& \lesssim 2^{\alpha k}\Big\|\Big(\sum\limits_{j=-\infty }^{\infty }t_{j}^{\beta}\left\vert
	f_{j}\right\vert ^{\beta}\Big)^{1/\beta}\chi_{\widetilde{R}_{k}}\Big\|_{p}
	\end{align*}
	for any $k\in\mathbb{Z}$. Clearly that the last term in $\ell_{q}$-norm  is bounded by
		\begin{align*}
	c\Big\|\Big(\sum\limits_{j=-\infty }^{\infty }t_{j}^{\beta}\left\vert
	f_{j}\right\vert ^{\beta}\Big)^{1/\beta}\Big\|_{\dot{K}_{p}^{\alpha,q}}.
		\end{align*}
				\textit{Subtep 1.2.} We see that 
		\begin{align*}
		& 2^{\alpha k}\Big\|\Big(\sum\limits_{j=-\infty }^{\infty }t_{j}^{\beta}\Big(\sum_{l=-\infty
		}^{j+m}2^{(l-j)K}\mathcal{M}(f_{l}\chi _{B_{k-2}})\Big)^{\beta}\Big)^{1/\beta}\chi_{k}
		\Big\|_{p}  
		\end{align*}
		is just
		\begin{equation*}
		\sup \sum\limits_{j=-\infty }^{\infty }\int_{R_{k}} 2^{\alpha k}t_{j}(x)\Big(\sum_{l=-\infty
		}^{j+m}2^{(l-j)K}\mathcal{M}(f_{l}\chi _{B_{k-2}})(x)\Big)|g_{j}(x)|dx=\sup \sum\limits_{j=-\infty }^{\infty }\sum_{l=-\infty
	}^{j+m}V_{j,k,l},
		\end{equation*}%
		where the supremum is taken over all sequence of functions $\{g_{j}\chi_{k}\}_{j\in 
			\mathbb{Z}}\in L_{p^{\prime }}(\ell _{\beta ^{\prime }})$ with 
		\begin{equation}
		\big\|\{g_{j}\chi_{k}\}_{j\in \mathbb{Z}}\big\|_{L_{p^{\prime }}(\ell _{\beta
				^{\prime }})}\leq 1,\label{dualherz1}
		\end{equation}%
		where $p^{\prime }$ and $\beta ^{\prime }$ are the conjugate exponent of $p$
		and $\beta $, respectively and
		\begin{equation*}
		V_{j,k,l}=2^{(l-j)K+\alpha k}\int_{R_{k}}t_{j}(x)\mathcal{M}(f_{l}\chi _{B_{k-2}})(x)|g_{j}(x)|dx.
		\end{equation*}
			In view the proof of Lemma 	\ref{maximalHerz},  we obtain 
		\begin{equation*}
		V_{j,k,l}\lesssim \Lambda_{j,l}
		\sum_{i=-\infty }^{k-2}2^{(k-i)(\alpha -\frac{n}{p}+n)}\int_{\mathbb{R}%
			^{n}}\mathcal{M}_{\frac{p}{\theta}}(2^{\alpha
			i}t_{l}f_{l}\chi _{i})(x)\mathcal{M}_{\tau ^{\prime }}(g_{j}\chi _{k})(x)dx,
		\end{equation*}%
		where  $\tau=p(1+\epsilon)$ and 
		\begin{equation*}
		\Lambda _{j,l}=\left\{ 
		\begin{array}{ccc}
		2^{(\alpha _{2}-K)(j-l)}, & \text{if} & l\leq j, \\ 
		2^{(\alpha _{1}-K)(j-l)}, & \text{if} & j<l\leq j+m.%
		\end{array}
		\right.
		\end{equation*}
		Since $K>\alpha _{2}$ by Lemmas  	\ref{lq-inequality1} and \ref{lq-inequality1copy(1)} we obtain that
		\begin{equation*}
	\int_{\mathbb{R}^{n}}\sum\limits_{j=-\infty }^{\infty }\Big(\sum_{l=-\infty
		}^{j+m}\Lambda _{j,l}\mathcal{M}_{\frac{p}{\theta}}(t_{l}f_{l}\chi _{i})(x)\Big)\mathcal{M}_{\tau ^{\prime }}(g_{j}\chi _{k})(x)dx.
		\end{equation*}%
		can be estimated from above by
		\begin{align*}
		& c\Big\|\Big(\sum\limits_{j=-\infty }^{\infty }\big(\mathcal{M}_{\frac{p}{\theta}
		}(t_{j}\mathcal{M}(f_{j}\chi_{i}))\big)^{\beta}\Big)^{1/\beta}\Big\|_{p } 
		\Big\|\Big(\sum\limits_{j=-\infty }^{\infty }\big(\mathcal{M}_{\tau^{\prime}
			}(g_{j}\chi_{k})\big)^{\beta^{\prime }}\Big)^{1/\beta^{\prime }}\Big\|_{p^{\prime} }\\
		& \lesssim \Big\|\Big(\sum\limits_{j=-\infty }^{\infty }t_{j}^{\beta}|f_{j}|^{\beta} %
		\Big)^{1/\beta}\chi_{i}\Big\|_{p },
		\end{align*}
		where we used the vector-valued maximal inequality of Fefferman and Stein $%
		\mathrm{\eqref{Fe-St71}}$, Lemma \ref{key-estimate2} and \eqref{dualherz1}. Consequently, $\sup \sum\limits_{j=-\infty }^{\infty }\sum_{l=-\infty
		}^{j+m}V_{j,k,l}$ is bounded by
	
	\begin{align*}
	c \sum_{i=-\infty }^{k-2}2^{(k-i)(\alpha -\frac{n}{p}+n)+\alpha i}\Big\|\Big(\sum\limits_{j=-\infty }^{\infty }t_{j}^{\beta}|f_{j}|^{\beta} %
	\Big)^{1/\beta}\chi_{i}\Big\|_{p },
	\end{align*}
	By Lemma \ref{lem:lq-inequality} the last expression in $\ell_{q}$-norm is bounded by 
	\begin{align*}
	c \Big\|\Big(\sum\limits_{j=-\infty }^{\infty }t_{j}^{\beta}\left\vert
	f_{j}\right\vert ^{\beta}\Big)^{1/\beta}\Big\|_{\dot{K}_{p}^{\alpha,q}}.
	\end{align*}
	\textit{Subtep 1.3.}  We have
		\begin{align*}
		& 2^{\alpha k}\Big\|\Big(\sum\limits_{j=-\infty }^{\infty }t_{j}^{\beta}\Big(\sum_{l=-\infty
		}^{j+m}2^{(l-j)K}\mathcal{M}(f_{l}\chi _{{\mathbb{R}^{n}}%
		\setminus B_{k+2}})\Big)^{\beta}\Big)^{1/\beta}\chi_{k}
		\Big\|_{p}  
		\end{align*}
		is just
		\begin{equation*}
		\sup \sum\limits_{j=-\infty }^{\infty }\int_{R_{k}} 2^{\alpha k}t_{j}(x)\Big(\sum_{l=-\infty
		}^{j+m}2^{(l-j)K}\mathcal{M}(f_{l}\chi _{{\mathbb{R}^{n}}%
		\setminus B_{k+2}})(x)\Big)|g_{j}(x)|dx=\sup \sum\limits_{j=-\infty }^{\infty }\sum_{l=-\infty
		}^{j+m}H_{j,k,l},
		\end{equation*}%
		where the supremum is taken over all sequence of functions $\{g_{j}\chi_{k}\}_{j\in 
			\mathbb{Z}}\in L_{p^{\prime }}(\ell _{\beta ^{\prime }})$ with 
		\begin{equation}
		\big\|\{g_{j}\chi_{k}\}_{j\in \mathbb{Z}}\big\|_{L_{p^{\prime }}(\ell _{\beta
				^{\prime }})}\leq 1,\label{dualherz2}
		\end{equation}%
		where $p^{\prime }$ and $\beta ^{\prime }$ are the conjugate exponent of $p$
		and $\beta $, respectively and
		\begin{equation*}
		H_{j,k,l}=2^{(l-j)K+\alpha k}\int_{R_{k}}t_{j}(x)\mathcal{M}(f_{l}\chi _{{\mathbb{R}^{n}}%
			\setminus B_{k+2}})(x)|g_{j}(x)|dx.
		\end{equation*}
	In view the proof of Lemma 	\ref{maximalHerz},  we obtain 
		\begin{equation*}
		H_{j,k,l}\lesssim \Lambda _{j,l} 2^{(l-j)K+\alpha k}\sum_{i=k+3}^{\infty }2^{(k-i)\alpha+i\alpha+\frac{(k-i)n\delta }{%
				\tau }}\int_{R_{i}} M_{B_{i},\frac{p}{\theta }}(t_{j}f_{l}\chi _{i})M_{B_{i},\tau ^{\prime
		}}(g_{j}\chi _{k})dx.
		\end{equation*}
	The desired estimate follows by the same arguments as in Subtep 1.2. The proof is complete.
	
		\textit{Step 2. Proof \eqref{key-est1.1.2.1}}. The proof is essentially the same as that of Step 1.
	
\end{proof}

Combining the arguments used in the proof of Lemma \ref{key-estimate3} with those of Lemma 4 in \cite{DRu25}, we obtain the following result.
\begin{lem}
	\label{key-estimate1.1}Let $m\in \mathbb{Z}$, $K\geq 0,1<\theta \leq
	p<\infty ,1<q\leq \infty \ $and $\upsilon =(\upsilon _{1},\upsilon _{2})\in 
	\mathbb{R}^{2}$. Let $\{t_{k}\}_{k\in \mathbb{Z}}\in \dot{X}_{\upsilon ,\sigma ,p}$ be a $p$
	-admissible weight sequence with $\sigma =(\sigma _{1}=\theta \left( \frac{p%
	}{\theta }\right) ^{\prime },\sigma _{2}\geq p)$. Let $\delta $ be as in Lemma \ref{Ap-Property}. Let $-\frac{%
	n\delta }{p}<\alpha <n-\frac{n}{p}$. Then for all sequence of
	functions  $%
	\{t_{k}f_{k}\}_{k\in \mathbb{Z}}\in \dot{K}_{p}^{\alpha
		,q}(\ell _{\beta})$, 
	\begin{align}
	& \Big\|\sup_{k\in \mathbb{Z} }\Big(t_{k}\sum_{j=-\infty
	}^{k+m}2^{(j-k)K}\mathcal{M}(f_{j})\Big)\Big\|_{\dot{K}_{p}^{\alpha,q}}  \notag \\
	& \lesssim \big\|\sup_{k\in \mathbb{Z} }\big(t_{k}\left\vert
	f_{k}\right\vert \big)\big\|_{\dot{K}_{p}^{\alpha,q}}
	\label{key-est1.1.1herz}
	\end{align}%
	if $K>\alpha _{2}$ and 
	\begin{align}
	& \Big\|\sup_{k\in \mathbb{Z} }\Big(t_{k}
	\sum_{j=k+m}^{\infty }2^{(j-k)K}\mathcal{M}(f_{j})\Big)%
	\Big\|_{\dot{K}_{p}^{\alpha,q}}  \notag \\
	& \lesssim  \big\|\sup_{k\in \mathbb{Z} }\big(t_{k}\left\vert
	f_{k}\right\vert \big)\big\|_{\dot{K}_{p}^{\alpha,q}}
	\label{key-est1.1.2.}
	\end{align}%
	if $K<\alpha _{1}$, with the usual modifications if $q=\infty $..
\end{lem}
\begin{proof}
	We adopt the approach used in \cite{DRu25}. For  a sequence of functions $\{g_{k}\}_{k\in \mathbb{Z}}$  and $\alpha \in \mathbb{R}$ we define 
	\begin{equation}
	T(\{g_{k}\}_{k\in \mathbb{Z}})=\Big\{\mathcal{M}\big(|x|^{\alpha}t_{k}\sum_{j=-\infty }^{k+m}2^{(j-k)K}%
	\mathcal{M}(|x|^{-\alpha}t_{j}^{-1}g_{j})\big)\Big\}_{k\in \mathbb{Z}}.  \label{operator2}
	\end{equation}%
	$T$ is a sublinear operator. The remainder of the proof follows by adapting the arguments in \cite{DRu25} and those used in Lemma \ref{key-estimate3}.
	\end{proof}
\begin{rem}
\label{r-estimates copy(1)}Let $i\in \mathbb{Z},1<\theta \leq p<\infty
,1<q,\beta<\infty \ $and $\upsilon =(\upsilon _{1},\upsilon _{2})\in \mathbb{R}^{2}$.
Let $\{t_{k}\}\in \dot{X}_{\upsilon ,\sigma ,p}$ be a $p$-admissible weight
sequence with $\sigma =(\sigma _{1}=\theta \left( \frac{p}{\theta} \right) ^{\prime
},\sigma _{2}\geq p)$. From Lemma \ref{key-estimate1.1} we easily obtain 
\begin{equation*}
\Big\|\Big(\sum\limits_{k=-\infty }^{\infty }t_{k}^{\beta}\big(\mathcal{M}%
(f_{k+i})\big)^{\beta}\Big)^{1/\beta}\Big\|_{\dot{K}_{p}^{\alpha,q}}\lesssim \Big\|\Big(%
\sum\limits_{k=-\infty }^{\infty }t_{k}^{\beta}\left\vert f_{k}\right\vert ^{\beta}%
\Big)^{1/\beta}\Big\|_{\dot{K}_{p}^{\alpha,q}}
\end{equation*}%
holds for all sequence of functions $\{t_{k}f_{k}\}_{k\in \mathbb{Z}} \in {\dot{K}_{p}^{\alpha,q}}(\ell _{\beta})$,
where the implicit constant depends on $i$. Indeed, we have 
\begin{equation*}
\mathcal{M}(f_{k+i})\leq \sum\limits_{j=-\infty }^{k+i}2^{(j-k-i)M}\mathcal{M%
}(f_{j}),\quad M>\alpha _{2},k\in \mathbb{Z}.
\end{equation*}%
Lemma \ref{key-estimate1.1} yields the desired result.
\end{rem}

\section{The spaces ${\dot{K}^{\alpha,q}_{p}\dot{A}_{\beta}(\mathbb{R}^{n},\{t_{k}\}_{k\in \mathbb{Z}})}$}

In this section, we introduce the Fourier-analytic definition of Herz-type Besov and Triebel-Lizorkin spaces with variable smoothness and investigate their fundamental properties, in analogy with the classical Besov and Triebel-Lizorkin spaces.

Select a pair of Schwartz functions  $\varphi $ and $\psi $ satisfying
\begin{equation}
\supp \mathcal{F}\varphi ,\supp\mathcal{F}\psi \subset \big\{\xi :\frac{1}{2}%
\leq |\xi |\leq 2\big\},  \label{Ass1}
\end{equation}%
\begin{equation}
|\mathcal{F}\varphi (\xi )|,|\mathcal{F}\psi (\xi )|\geq c\quad \text{if}%
\quad \frac{3}{5}\leq |\xi |\leq \frac{5}{3}  \label{Ass2}
\end{equation}%
and 
\begin{equation}
\sum_{k=-\infty }^{\infty }\overline{\mathcal{F}\varphi (2^{-k}\xi )}%
\mathcal{F}\psi (2^{-k}\xi )=1\quad \text{if}\quad \xi \neq 0,  \label{Ass3}
\end{equation}%
where $c>0$. Throughout the paper, for all $k$ $\in \mathbb{Z}$ and $x\in 
\mathbb{R}^{n}$, we put $\varphi _{k}(x):=2^{kn}\varphi (2^{k}x)$. Let $%
\varphi \in \mathcal{S}(\mathbb{R}^{n})$ be a function satisfying $\mathrm{\ %
	\eqref{Ass1}}$-$\mathrm{\eqref{Ass2}}$. Recall that there exists a function $%
\psi \in \mathcal{S}(\mathbb{R}^{n})$ satisfying $\mathrm{\eqref{Ass1}}$-$%
\mathrm{\eqref{Ass3}}$, see \cite[Lemma (6.9)]{FrJaWe01}. We set%
\begin{equation*}
\mathcal{S}_{\infty }(\mathbb{R}^{n})=\Big\{\varphi \in \mathcal{S}(\mathbb{%
	R}^{n}):\int_{\mathbb{R}^{n}}x^{\beta }\varphi (x)dx=0\text{ for all
	multi-indices }\beta \in \mathbb{N}_{0}^{n}\Big\}.
\end{equation*}%
Let $\mathcal{S}_{\infty }^{\prime }(\mathbb{R}^{n})$ be the topological
dual of $\mathcal{S}_{\infty }(\mathbb{R}^{n})$, namely, the set of all
continuous linear functionals on $\mathcal{S}_{\infty }(\mathbb{R}^{n})$.

Now, we define the spaces under consideration.

\begin{defn}
	\label{B-F-def}Let $0<p,q,\beta\leq \infty $. Let $\{t_{k}\}_{k\in \mathbb{Z}}$
	be a $p$-admissible weight sequence, and $\varphi \in \mathcal{S}(\mathbb{R}%
	^{n})$\ satisfy $\mathrm{\eqref{Ass1}}$ and $\mathrm{\eqref{Ass2}}$.\newline
	$\mathrm{(i)}$ The Herz-type Besov space $\dot{K}^{\alpha,q}_{p}\dot{B}_{\beta}(\mathbb{R}^{n},\{t_{k}\}_{k\in \mathbb{Z}})$\ is
	the collection of all $f\in \mathcal{S}_{\infty }^{\prime }(\mathbb{R}^{n})$
	\ such that 
	\begin{equation*}
	\big\|f\big\|_{\dot{K}^{\alpha,q}_{p}\dot{B}_{\beta}(\mathbb{R}^{n},\{t_{k}\}_{k\in \mathbb{Z}})}=\Big(%
	\sum\limits_{k=-\infty }^{\infty }\big\|t_{k}(\varphi _{k}\ast f)\big\|_{\dot{K}^{\alpha,q}_{p}}^{\beta}\Big)^{1/\beta}<\infty
	\end{equation*}%
	with the usual modification if $\beta=\infty $.\newline
	$\mathrm{(ii)}$ Let $0<p,q<\infty $. The Herz-type Triebel-Lizorkin space $\dot{K}^{\alpha,q}_{p}\dot{F}_{\beta}(\mathbb{R}^{n},\{t_{k}\}_{k\in \mathbb{Z}})$\ is the collection of all $f\in \mathcal{S}%
	_{\infty }^{\prime }(\mathbb{R}^{n})$\ such that 
	\begin{equation*}
	\big\|f\big\|_{\dot{K}^{\alpha,q}_{p}\dot{F}_{\beta}(\mathbb{R}^{n},\{t_{k}\}_{k\in \mathbb{Z}})}=\Big\|\Big(%
	\sum\limits_{k=-\infty }^{\infty }t_{k}^{\beta}|\varphi _{k}\ast f|^{\beta}\Big)^{
		1/\beta}\Big\|_{\dot{K}^{\alpha,q}_{p}}<\infty
	\end{equation*}%
	with the usual modification if $\beta=\infty $.
\end{defn}
\begin{rem}
	Using the system $\{\varphi _{k}\}_{k\in \mathbb{Z}}$ we can define the
	quasi-norms%
	\begin{equation*}
	\big\|f\big\|_{\dot{B}_{p,q}(\mathbb{R}^{n},\{t_{k}\}_{k\in\mathbb{Z}})}=\Big(%
	\sum\limits_{k=-\infty }^{\infty }\big\|t_{k}(\varphi _{k}\ast f)\big\|_{p}^{\beta}\Big)^{1/\beta}<\infty
\end{equation*}%
with the usual modification if $\beta=\infty $
	and%
	\begin{equation*}
	\big\|f\big\|_{\dot{F}_{p,q}(\mathbb{R}^{n},\{t_{k}\}_{k\in\mathbb{Z}})}=\Big\|\Big(%
	\sum\limits_{k=-\infty }^{\infty }t_{k}^{\beta}|\varphi _{k}\ast f|^{\beta}\Big)^{\frac{1}{\beta}}\Big\|_{p}
	\end{equation*}%
	where  $0<p,q,\beta\leq \infty $ with $0<p,q<\infty $
	in the $\dot{F}$-case. The   weighted Besov space $\dot{B}_{p,\beta}(\mathbb{R}^{n},\{t_{k}\}_{k\in\mathbb{Z}})$\
	consist of all distributions $f\in \mathcal{S}_{\infty }^{\prime }(\mathbb{R}%
	^{n})$ for which $\big\|f\big\|_{\dot{B}_{p,\beta}(\mathbb{R}^{n},\{t_{k}\}_{k\in\mathbb{Z}})}<\infty .$
	The weighted Triebel-Lizorkin space $\dot{F}_{p,\beta}(\mathbb{R}^{n},\{t_{k}\}_{k\in\mathbb{Z}})$\ consist of
	all distributions $f\in \mathcal{S}_{\infty }^{\prime }(\mathbb{R}^{n})$ for
	which $\big\|f\big\|_{\dot{F}_{p,\beta}(\mathbb{R}^{n},\{t_{k}\}_{k\in\mathbb{Z}})}<\infty .$
	
Several fundamental properties of these function spaces, such as their \(\varphi\)-transform characterization in the sense of Frazier and Jawerth, duality, complex interpolation, and smooth atomic and molecular decompositions, were established in \cite{D20}, \cite{D20.1}, \cite{D20.2}, \cite{DRu25}, \cite{D20.3} and \cite{DFA26}.
	
	One recognizes immediately that if $p=q$ and $\alpha=0$, then 
	\begin{equation*}
	\dot{K}^{0,p}_{p}\dot{B}_{\beta}(\mathbb{R}^{n},\{t_{k}\}_{k\in \mathbb{Z}})=\dot{B}_{p,\beta}(\mathbb{R}^{n},\{t_{k}\}_{k\in \mathbb{Z}}),\quad \dot{K}^{0,p}_{p}\dot{F}_{\beta}(\mathbb{R}^{n},\{t_{k}\}_{k\in \mathbb{Z}})=\dot{F}_{p,\beta}(\mathbb{R}^{n},\{t_{k}\}_{k\in \mathbb{Z}}).
	\end{equation*}%
	Moreover, 
		\begin{equation*}
	\dot{K}^{\alpha,q}_{p}\dot{B}_{\beta}(\mathbb{R}^{n},\{2^{ks}\}_{k\in \mathbb{Z}})=\dot{K}^{\alpha,q}_{p}\dot{B}^{s}_{\beta},\quad \dot{K}^{\alpha,q}_{p}\dot{F}_{\beta}(\mathbb{R}^{n},\{2^{ks}\}_{k\in \mathbb{Z}})=\dot{K}^{\alpha,q}_{p}\dot{F}^{s}_{\beta},
	\end{equation*}
		where $\dot{K}^{\alpha,q}_{p}\dot{B}^{s}_{\beta}$ and $\dot{K}^{\alpha,q}_{p}\dot{F}^{s}_{\beta}$, $s\in \mathbb{R}$ are the  Herz-type Besov and Triebel-Lizorkin spaces, respectively. In addition
		\begin{equation*}
		\dot{K}^{0,p}_{p}\dot{B}_{\beta}(\mathbb{R}^{n},\{2^{ks}\}_{k\in \mathbb{Z}})=\dot{B}^{s}_{p,\beta},\quad \dot{K}^{0,p}_{p}\dot{F}_{\beta}(\mathbb{R}^{n},\{2^{ks}\}_{k\in \mathbb{Z}})=\dot{F}^{s}_{p,\beta},
		\end{equation*}
		where $\dot{B}^{s}_{p,\beta}$ and $\dot{F}^{s}_{p,\beta}$ are the   Besov and Triebel-Lizorkin spaces, respectively. Further details on
		the classical theory of  Besov and Triebel-Lizorkin spaces, can
		be found in \cite{FJ86}, \cite{FJ90}, \cite{FrJaWe01}, \cite{SchTr87}, \cite%
		{T1}, \cite{T2} and \cite{T3}.
		Further results, concerning, for instance lifting properties,
		Fourier multiplier and local means characterizations of Herz type Besov and Triebel-Lizorkin spaces can be found in \cite%
		{Drihem1.13}, \cite{Drihem2.13}, \cite{Drmana}, \cite{XuYang05} and \cite{Xu05}%
		. Some applications of such function spaces in semilinear parabolic
		equations are given in \cite{Dr22.Banach}.
	\end{rem}

Let $\varphi $, $\psi \in \mathcal{S}(\mathbb{R}^{n})$ satisfying $\mathrm{\ %
\eqref{Ass1}}$\ through\ $\mathrm{\eqref{Ass3}}$. Recall that the $\varphi $%
-transform $S_{\varphi }$ is defined by setting 
\begin{equation*}
(S_{\varphi}f)_{k,m}=\langle f,\varphi _{k,m}\rangle,
\end{equation*}
where $\varphi _{k,m}(x)=2^{k \frac{n}{2}}\varphi (2^{k}x-m)$, $m\in \mathbb{%
Z}^{n}$ and $k\in \mathbb{Z}$. The inverse $\varphi $-transform $T_{\psi }$
is defined by 
\begin{equation*}
T_{\psi }\lambda =\sum_{k=-\infty }^{\infty }\sum_{m\in \mathbb{Z}
^{n}}\lambda _{k,m}\psi _{k,m},
\end{equation*}%
where $\lambda =\{\lambda _{k,m}\}_{k\in \mathbb{Z},m\in \mathbb{Z}
^{n}}\subset \mathbb{C}$, see \cite{FJ90}.

Now we introduce the corresponding sequence spaces of $\dot{K}^{\alpha,q}_{p}\dot{A}_{\beta}(\mathbb{R}^{n},\{t_{k}\}_{k\in \mathbb{Z}})$.

\begin{defn}
\label{sequence-space}Let $0<p,q,\beta\leq\infty $. Let $%
\{t_{k}\}_{k\in \mathbb{Z}} $ be a $p$-admissible weight sequence. \newline
$\mathrm{(i)}$ For all complex valued
sequences $\lambda =\{\lambda _{k,m}\}_{k\in \mathbb{Z},m\in \mathbb{Z}
	^{n}}\subset \mathbb{C}$ we define 
\begin{equation*}
\dot{K}^{\alpha,q}_{p}\dot{b}_{\beta}(\mathbb{R}^{n},\{t_{k}\}_{k\in \mathbb{Z}})=\Big\{\lambda :\big\|\lambda \big\|_{\dot{K}^{\alpha,q}_{p}\dot{b}_{\beta}(\mathbb{R}^{n},\{t_{k}\}_{k\in \mathbb{Z}})}<\infty \Big\},
\end{equation*}
where 
\begin{equation*}
\big\|\lambda\big\|_{\dot{K}_{p}^{\alpha,q}\dot{b}_{\beta }(\mathbb{R}^{n},\{t_{k}\}_{k\in \mathbb{Z}})}=\Big(%
\sum_{k=-\infty }^{\infty }2^{\frac{kn\beta}{2}}\Big\| \sum\limits_{m\in \mathbb{Z}%
	^{n}}\lambda _{k,m}t_{k}\chi _{k,m}\Big\|_{\dot{K}_{p}^{\alpha,q}}^{\beta}\Big)^{1/\beta},
\end{equation*}
with the usual modifications if $\beta=\infty $.\newline

$\mathrm{(ii)}$ Let $0<p,q<\infty $ and $0<\beta\leq \infty $. For all complex valued
sequences $\lambda =\{\lambda _{k,m}\}_{k\in \mathbb{Z},m\in \mathbb{Z}
^{n}}\subset \mathbb{C}$ we define 
\begin{equation*}
\dot{K}^{\alpha,q}_{p}\dot{f}_{\beta}(\mathbb{R}^{n},\{t_{k}\}_{k\in \mathbb{Z}})=\Big\{\lambda :\big\|\lambda \big\|_{\dot{K}^{\alpha,q}_{p}\dot{F}_{\beta}(\mathbb{R}^{n},\{t_{k}\}_{k\in \mathbb{Z}})}<\infty \Big\},
\end{equation*}
where 
\begin{equation*}
\big\|\lambda \big\|_{\dot{K}^{\alpha,q}_{p}\dot{F}_{\beta}(\mathbb{R}^{n},\{t_{k}\}_{k\in \mathbb{Z}})}=\Big\|\Big( %
\sum_{k=-\infty }^{\infty }\sum\limits_{m\in \mathbb{Z}
^{n}}2^{\frac{kn\beta}{2}}t_{k}^{\beta}|\lambda _{k,m}|^{\beta}\chi _{k,m}\Big)^{1/\beta}\Big\|_{\dot{K}^{\alpha,q}_{p}},
\end{equation*}
with the usual modifications if $\beta=\infty $.
\end{defn}
Allowing the smoothness $t_{k}$, $k\in \mathbb{Z}$ to vary from point to
point will raise extra difficulties\ to study these function spaces. But by
the following lemma the problem can be reduced to the case of fixed
smoothness, which is a slight variant of \cite{D20.2}. For the convenience
of the reader, we give some details.
\begin{prop}
	\label{Equi-norm1herz}Let $0<\theta \leq p<\infty $, $0<q,\beta\leq\infty $ and $0<\kappa
	\leq 1$. Assume that $\{t_{k}\}_{k\in \mathbb{Z}}$ satisfying $\mathrm{\eqref{Asum1}}$ with $%
	\sigma _{1}=\theta \left( \frac{p}{\theta }\right) ^{\prime }$ and $j=k$.  Let $\delta $ be as in Lemma \ref{Ap-Property} and $\alpha>-\frac{
		n\delta }{p}$.\newline
	$\mathrm{(i)}$ Then%
	\begin{equation*}
	\big\|\lambda\big\|_{\dot{K}_{p}^{\alpha,q}\dot{b}_{\beta,\kappa }(\mathbb{R}^{n},\{t_{k}\}_{k\in \mathbb{Z}})}^{\ast }:=\Big(%
	\sum_{k=-\infty }^{\infty }2^{kn\beta(\frac{1}{2}+\frac{1}{\kappa p})}\Big\| \sum\limits_{m\in \mathbb{Z}%
		^{n}}|\lambda _{k,m}|t_{k,m,\kappa p}\chi _{k,m}\Big\|_{\dot{K}_{p}^{\alpha,q}}^{\beta}\Big)^{\frac{1}{\beta}},
	\end{equation*}%
	\textrm{\ }is an equivalent quasi-norm in $\dot{K}_{p}^{\alpha,q}\dot{b}_{\beta }(\mathbb{R}^{n},\{t_{k}\}_{k\in \mathbb{Z}})$.\newline
	
	$\mathrm{(ii)}$ Let  $0<q<\infty $ . Assume that $\{t_{k}\}_{k\in \mathbb{Z}}$ satisfying $\mathrm{\eqref{Asum1}}$ with $%
	\sigma _{1}=\theta \left( \frac{p}{\theta }\right) ^{\prime }$ and $j=k$.
	Then%
	\begin{equation*}
	\big\|%
	\lambda |\big\|_{\dot{K}_{p}^{\alpha,q}\dot{f}_{\beta,\kappa }(\mathbb{R}^{n},\{t_{k}\}_{k\in \mathbb{Z}})}^{\ast }:=\Big\|\Big(\sum_{k=-\infty }^{\infty }\sum\limits_{m\in \mathbb{Z}%
		^{n}}2^{kn\beta(\frac{1}{2}+\frac{1}{\kappa p})}t_{k,m,\kappa p}^{\beta}|\lambda
	_{k,m}|^{\beta}\chi _{k,m}\Big)^{\frac{1}{\beta}}\Big\|_{\dot{K}_{p}^{\alpha,q}},
	\end{equation*}%
	is an equivalent quasi-norm in $\dot{K}_{p}^{\alpha,q}\dot{f}_{\beta }(\mathbb{R}^{n},\{t_{k}\}_{k\in \mathbb{Z}})$,
	where%
	\begin{equation*}
	t_{k,m,\kappa p}:=\big\|t_{k}|L_{\kappa p}(Q_{k,m})\big\|,\quad k\in \mathbb{%
		Z},m\in \mathbb{Z}^{n}.
	\end{equation*}
\end{prop}

\begin{proof}
	It suffices to prove (ii), since (i) can be established by a similar argument. We proceed in two steps.
	
	\textit{Step 1.} Let us prove that%
	\begin{equation}
	\big\|\lambda \big\|_{\dot{K}_{p}^{\alpha,q}\dot{f}_{\beta}(\mathbb{R}^{n},\{t_{k}\}_{k\in \mathbb{Z}})}\lesssim \big\|%
	\lambda |\big\|_{\dot{K}_{p}^{\alpha,q}\dot{f}_{\beta,\kappa }(\mathbb{R}^{n},\{t_{k}\}_{k\in \mathbb{Z}})}^{\ast }.
	\label{new-norm}
	\end{equation}%
	Let $\eta>0$ be such that
	\begin{equation}
	0<\eta <\min \Big(\frac{pn\delta_{\{t^{-1}\}_{k\in\mathbb{{Z}}}}}{n\delta_{\{t_{k}^{-1}\}_{k\in\mathbb{{Z}}}}
		+p\alpha },q,\beta,\theta \Big),\label{varrho}
	\end{equation}
	where $ \delta_{\{t_{k}^{-1}\}_{k\in\mathbb{{Z}}}} $ is obtained by applying  Lemma \ref{Ap-Property}/(vi) to the sequence $\{t_{k}^{-1}\}_{k\in\mathbb{{Z}}}$. We have $\big\|\lambda \big\|^{\eta}_{\dot{K}_{p}^{\alpha,q}\dot{f}_{\beta,\kappa}(\mathbb{R}^{n},\{t_{k}\}_{k\in \mathbb{Z}})}$ is just%
	\begin{equation*}
	\sup \sum_{k=-\infty }^{\infty }\sum\limits_{m\in \mathbb{Z}^{n}}2^{\frac{%
			n\eta }{2}k}|\lambda _{k,m}|^{\eta }\int_{Q_{k,m}}t_{k}^{\eta
	}(x)|g_{k}(x)|dx=\sup \sum_{k=-\infty }^{\infty }S_{k},
	\end{equation*}%
	where the supremum is taken over all sequences of functions $\{g_{k}\}_{k\in \mathbb{Z}}\in \dot{K}_{(\frac{p}{\eta})^{\prime}}^{-\alpha\eta,(\frac{q}{\eta})^{\prime}}
	(\ell _{(\frac{\beta}{\eta })^{\prime }})$ such that
	\begin{equation}
	\big\|\{g_{k}\}_{k\in\mathbb{Z}}\big\|_{\dot{K}_{(\frac{p}{\eta})^{\prime}}^{-\alpha\eta,(\frac{q}{\eta})^{\prime}}
		(\ell _{(\frac{\beta}{\eta })^{\prime }})}\leq 1.  \label{dual}
	\end{equation}%
	By H\"{o}lder's inequality,%
	\begin{equation*}
	1=\Big(\frac{1}{|Q_{k,m}|}\int_{Q_{k,m}}t_{k}^{-\eta h}(x)t_{k}^{\eta h}(x)dx%
	\Big)^{\frac{1}{h}}\leq M_{Q_{k,m},\varrho }\mathcal{(}t_{k}^{-\eta })\big(%
	M_{Q_{k,m},\kappa p}\mathcal{(}t_{k})\big)^{\eta }
	\end{equation*}%
	for any $h>0$, with $\frac{1}{h}=\frac{1}{\varrho }+\frac{\eta }{\kappa p}$.
	Therefore,%
	\begin{equation*}
	\int_{Q_{k,m}}t_{k}^{\eta }(x)|g_{k}(x)|dx\lesssim |Q_{k,m}|\big(%
	M_{Q_{k,m},\kappa p}\mathcal{(}t_{k})\big)^{\eta }M_{Q_{k,m},\varrho }%
	\mathcal{(}t_{k}^{-\eta }\mathcal{M}(t_{k}^{\eta }g_{k})).
	\end{equation*}%
	We set%
	\begin{equation*}
	f_{k}(x)=\sum\limits_{m\in \mathbb{Z}^{n}}2^{\frac{kn}{2}}|\lambda
	_{k,m}||Q_{k,m}|^{-\frac{1}{\kappa p}}t_{k,m,\kappa p}\chi _{k,m}(x),\quad
	x\in \mathbb{R}^{n},k\in \mathbb{Z}.
	\end{equation*}%
	Take $0<\varrho <\min \big((\frac{p}{\eta })^{\prime },(\frac{\beta}{\eta })^{\prime
	}\big)$, we find by H\"{o}lder's inequality that the sum $\sum_{k=-\infty
	}^{\infty }S_{k}$ can be estimated from above by 
	\begin{align*}
	& c\int_{\mathbb{R}^{n}}\sum_{k=-\infty }^{\infty }f_{k}^{\eta }(x)\mathcal{M%
	}_{\varrho }(t_{k}^{-\eta }\mathcal{M}(t_{k}^{\eta }g_{k}))(x)dx \\
	\lesssim & \Big\|\Big(\sum_{k=-\infty }^{\infty }f_{k}^{\beta}\Big)^{\eta/ 
		\beta}\Big\|_{\dot{K}_{\frac{p}{\eta}}^{\alpha\eta,\frac{q}{\eta}}
	}\Big\|\Big(\sum_{k=-\infty
	}^{\infty }\left( \mathcal{M}_{\varrho }(t_{k}^{-\eta }\mathcal{M}%
	(t_{k}^{\eta }g_{k}))\right) ^{(\frac{\beta}{\eta })^{\prime }}\Big)^{1/
			\frac{\beta}{\eta })^{\prime }}
	\Big\|_{\dot{K}_{(\frac{p}{\eta})^{\prime}}^{-\alpha\eta,(\frac{q}{\eta})^{\prime}}}.
	\end{align*}%
	Observe that%
	\begin{equation*}
	t_{k}^{p}\in A_{\frac{p}{\theta }}\subset A_{\frac{p}{\eta }},\quad k\in 
	\mathbb{Z},
	\end{equation*}%
	which yields that%
	\begin{equation*}
	t_{k}^{-\eta (\frac{p}{\eta })^{\prime }}\in A_{(\frac{p}{\eta })^{\prime
	}},\quad k\in \mathbb{Z},
	\end{equation*}%
	by Lemma \ref{Ap-Property}/(i)-(ii). By the vector-valued maximal inequality
	of Fefferman and Stein $\mathrm{\eqref{Fe-St71}}$, Lemma \ref{key-estimate1}
	and $\mathrm{\eqref{dual}}$, we obtain that%
	\begin{equation*}
	\sum_{k=-\infty }^{\infty }S_{k}\lesssim \Big\|\Big(\sum_{k=-\infty
	}^{\infty }f_{k}^{\beta}\Big)^{\frac{1}{\beta}}\Big\|_{\dot{K}_{p}^{\alpha,q}}^{\eta },
	\end{equation*}%
	which yields $\mathrm{\eqref{new-norm}.}$
	
	\textit{Step 2.} We prove the opposite inequality of $\mathrm{%
		\eqref{new-norm}}$. We can write 
	\begin{equation*}
	|\lambda _{k,m}|^{\nu \theta }=\frac{1}{|Q_{k,m}|}\int_{Q_{k,m}}|\lambda
	_{k,m}|^{\nu \theta }\chi _{k,m}(y)dy,\quad k\in \mathbb{Z},m\in \mathbb{Z}%
	^{n}.
	\end{equation*}%
	Again\ by H\"{o}lder's inequality, since $\frac{\theta }{\sigma _{1}}+\frac{%
		\theta }{p}=1$, 
	\begin{align*}
	|\lambda _{k,m}|\leq & M_{Q_{k,m},\nu p}(\lambda _{k,m}t_{k})M_{Q_{k,m},\nu
		\sigma _{1}}(t_{k}^{-1}) \\
	\lesssim & \mathcal{M}_{\nu p}\big(\sum\limits_{m\in \mathbb{Z}^{n}}\lambda
	_{k,m}t_{k}\chi _{k,m}\big)(x)M_{Q_{k,m},\sigma _{1}}(t_{k}^{-1})
	\end{align*}%
	for any $x\in Q_{k,m}$ with $0<\nu <\min \big(1,\frac{n-\alpha p}{n},\frac{q}{p},\frac{\beta}{p}\big)$. Observing that%
	\begin{equation*}
	M_{Q_{k,m},\kappa p}(t_{k})\leq M_{Q_{k,m},p}(t_{k}),\quad k\in \mathbb{Z}%
	,m\in \mathbb{Z}^{n}.
	\end{equation*}%
	Hence by $\mathrm{\eqref{Asum1}}$\textrm{, }we find that%
	\begin{equation*}
	|\lambda _{k,m}|\lesssim 2^{-k\frac{n}{\kappa p}}t_{k,m,\kappa p}^{-1}%
	\mathcal{M}_{\nu p}\big(\sum\limits_{m\in \mathbb{Z}^{n}}t_{k}\lambda
	_{k,m}\chi _{k,m}\big)(x)
	\end{equation*}%
	for any $x\in Q_{k,m}$. Therefore,%
	\begin{align*}
	\big\|\lambda \big\|_{\dot{K}_{p}^{\alpha,q}\dot{f}_{\beta,\kappa }(\mathbb{R}^{n},\{t_{k}\}_{k\in \mathbb{Z}})}^{\ast }
	\lesssim & \Big\|\Big(\sum_{k=-\infty }^{\infty }\big(\mathcal{M}_{\nu
		p}(\sum\limits_{m\in \mathbb{Z}^{n}}2^{\frac{kn}{2}}t_{k}\lambda _{k,m}\chi
	_{k,m})\big)^{\beta}\Big)^{1/\beta}\Big\|_{\dot{K}_{p}^{\alpha,q} }\\
	\lesssim & \big\|\lambda \big\|_{\dot{K}_{p}^{\alpha,q}\dot{f}_{\beta}(\mathbb{R}^{n},\{t_{k}\}_{k\in \mathbb{Z}})},
	\end{align*}%
	where we used the vector-valued maximal inequality of Fefferman and Stein $%
	\mathrm{\eqref{Fe-St71}}$. The proof is complete.
\end{proof}

The following lemma is a slight variant of \cite{D20} and \cite{D20.1}.

\begin{thm}
\label{Inv-phi-trans}Let $\upsilon=(\upsilon _{1},\upsilon _{2})\in \mathbb{R}
^{2},0<\theta \leq p<\infty $ and $0<q,\beta\leq \infty $. Let $\{t_{k}\}_{k\in \mathbb{Z}}\in \dot{%
X	}_{\upsilon ,\sigma ,p}$ be a $p$-admissible weight sequence with $\sigma
=(\sigma _{1}=\theta \left( \frac{p}{\theta} \right) ^{\prime },\sigma _{2}\geq p)$.  
Let\ $\psi \in \mathcal{S}(\mathbb{R}^{n})$ satisfying $\mathrm{\eqref{Ass1}}
$ and $\mathrm{\eqref{Ass2}}$. Then for all $\lambda \in \dot{K}_{p}^{\alpha,q}\dot{a}_{\beta}(\mathbb{R}^{n},\{t_{k}\}_{k\in \mathbb{Z}})$ 
\begin{equation*}
T_{\psi }\lambda =\sum_{k=-\infty }^{\infty }\sum_{m\in \mathbb{Z}
^{n}}\lambda _{k,m}\psi _{k,m},
\end{equation*}
converges in $\mathcal{S}_{\infty }^{\prime }(\mathbb{R}^{n})$; moreover, $%
T_{\psi }:\dot{K}_{p}^{\alpha,q}\dot{a}_{\beta}(\mathbb{R}^{n},\{t_{k}\}_{k\in \mathbb{Z}})\rightarrow \mathcal{S}
_{\infty }^{\prime }(\mathbb{R}^{n})$ is continuous.
\end{thm}

For a sequence $\lambda =\{\lambda _{k,m}\}_{k\in \mathbb{Z},m\in \mathbb{Z}
^{n}}\subset \mathbb{C},0<r\leq \infty $ and a fixed $d>0$, set%
\begin{equation*}
\lambda _{k,m,r,d}^{\ast }=\Big(\sum_{h\in \mathbb{Z}^{n}}\frac{|\lambda
_{k,h}|^{r}}{(1+2^{k}|2^{-k}h-2^{-k}m|)^{d}}\Big)^{1/r}
\end{equation*}%
and $\lambda _{r,d}^{\ast }:=\{\lambda _{k,m,r,d}^{\ast }\}_{k\in \mathbb{Z}
,m\in \mathbb{Z}^{n}}\subset \mathbb{C}$. By adapting the arguments of \cite{D20.1}, and using  Lemma \ref{key-estimate1}, we obtain

\begin{lem}
\label{lamda-equi}Let $\upsilon =(\upsilon  _{1},\upsilon  _{2})\in \mathbb{R}
^{2},0<\theta \leq p<\infty ,0<q,\beta\leq\infty ,\gamma \in \mathbb{Z}$ and $d>n$.
Let $\{t_{k}\}_{k\in \mathbb{Z}}$ be a $p$-admissible weight sequence satisfying $\mathrm{\ %
\eqref{Asum1}}$ with $\sigma _{1}=\theta \left( p/\theta \right) ^{\prime }$
and $\upsilon _{1}\in \mathbb{R}$. Let $\delta $ be as in Lemma \ref{Ap-Property} and $\alpha>-\frac{
	n\delta }{p}$. Then 
\begin{equation}
\big\|\lambda _{r,d}^{\ast } \big\|_{\dot{K}_{p}^{\alpha,q}\dot{a}_{\beta}(\mathbb{R}%
	^{n},\{t_{k-\gamma }\}_{k\in \mathbb{Z}})}\approx \big\|\lambda  \big\|_{\dot{K}_{p}^{\alpha,q}\dot{a}_{\beta}(\mathbb{R}^{n},\{t_{k}\}_{k\in \mathbb{Z}})}.  \label{First}
\end{equation}
In addition if $\{t_{k}\}_{k\in \mathbb{Z}}$ satisfies $\mathrm{\eqref{Asum2}}$ with $\sigma
_{2}\geq p$ and $\alpha _{2}\in \mathbb{R}$, then 
\begin{equation}
\big\|\lambda _{r,d}^{\ast } \big\|_{\dot{K}_{p}^{\alpha,q}\dot{a}_{\beta}(\mathbb{R}%
	^{n},\{t_{k-\gamma }\}_{k\in \mathbb{Z}})}\leq \big\|\lambda  \big\|_{\dot{K}_{p}^{\alpha,q}\dot{a}_{\beta}(\mathbb{R}^{n},\{t_{k}\}_{k\in \mathbb{Z}})},  \label{Second}
\end{equation}
where
$r=\min(p,\beta)$ if $\dot{a}=\dot{f}$ and $r=p$ if $\dot{a}=\dot{b}$.
\end{lem}

We now state the following $\varphi$-transform characterization in the sense of Frazier and Jawerth, which plays a crucial role in the remainder of the paper. The proof follows closely the arguments in \cite{D20.1}, together with Lemma \ref{lamda-equi}.
\begin{thm}
\label{phi-tran}Let $\upsilon =(\upsilon _{1},\upsilon _{2})\in \mathbb{R}
^{2},0<\theta \leq p<\infty $ and$\ 0<q,\beta\leq\infty $. Let $\{t_{k}\}_{k\in \mathbb{Z}}\in \dot{X}
_{\upsilon ,\sigma ,p}$ be a $p$-admissible weight sequence with $\sigma
=(\sigma _{1}=\theta \left( \frac{p}{\theta} \right) ^{\prime },\sigma _{2}\geq p)$
. Let $\delta $ be as in Lemma \ref{Ap-Property} and $\alpha>-\frac{
	n\delta }{p}$. Let $\varphi $, $\psi \in \mathcal{S}(\mathbb{R}^{n})$ satisfying $%
\mathrm{\eqref{Ass1}}$\ through\ $\mathrm{\eqref{Ass3}}$. The operators 
\begin{equation*}
S_{\varphi }:\dot{K}_{p}^{\alpha,q}\dot{A}_{\beta}(\mathbb{R}^{n},\{t_{k}\}_{k\in \mathbb{Z}})\rightarrow \dot{K}_{p}^{\alpha,q}\dot{a}_{\beta}(\mathbb{R}^{n},\{t_{k}\}_{k\in \mathbb{Z}})
\end{equation*}
and 
\begin{equation*}
T_{\psi }:\dot{K}_{p}^{\alpha,q}\dot{a}_{\beta}(\mathbb{R}^{n},\{t_{k}\}_{k\in \mathbb{Z}})\rightarrow \dot{K}_{p}^{\alpha,q}\dot{A}_{\beta}(\mathbb{R}^{n},\{t_{k}\}_{k\in \mathbb{Z}})
\end{equation*}
are bounded. Furthermore, $T_{\psi }\circ S_{\varphi }$ is the identity on $%
\dot{K}_{p}^{\alpha,q}\dot{A}_{\beta}(\mathbb{R}^{n},\{t_{k}\}_{k\in \mathbb{Z}})$.
\end{thm}

\begin{rem}
This theorem can then be exploited to obtain a variety of results for the $%
\dot{K}_{p}^{\alpha,q}\dot{A}_{\beta}(\mathbb{R}^{n},\{t_{k}\}_{k\in \mathbb{Z}})$ spaces, where arguments can be
equivalently transferred to the sequence space, which is often more
convenient to handle. More precisely, under the same hypothesis of the last
theorem, 
\begin{equation*}
\big\|\{\langle f,\varphi _{k,m}\rangle \}_{k\in \mathbb{Z},m\in \mathbb{Z}
^{n}}\big\|_{\dot{K}_{p}^{\alpha,q}\dot{a}_{\beta}(\mathbb{R}^{n},\{t_{k}\}_{k\in \mathbb{Z}})}\approx \big\|f\big\|_{\dot{K}_{p}^{\alpha,q}\dot{a}_{\beta}(\mathbb{R}^{n},\{t_{k}\}_{k\in \mathbb{Z}})}.
\end{equation*}
\end{rem}

From Theorem \ref{phi-tran}, we obtain the next important property of the
function spaces\ $\dot{K}_{p}^{\alpha,q}\dot{A}_{\beta}(\mathbb{R}^{n},\{t_{k}\}_{k\in \mathbb{Z}})$.

\begin{cor}
\label{Indpendent}Let $\upsilon =(\upsilon _{1},\upsilon _{2})\in \mathbb{R}
^{2},0<\theta \leq p<\infty $ and $0<q,\beta\leq\infty $. Let $\{t_{k}\}_{k\in \mathbb{Z}}\in \dot{X}
_{\upsilon ,\sigma ,p}$ be a $p$-admissible weight sequence with $\sigma
=(\sigma _{1}=\theta \left( \frac{p}{\theta} \right) ^{\prime },\sigma _{2}\geq p)$. Let $\delta $ be as in Lemma \ref{Ap-Property} and $\alpha>-\frac{
	n\delta }{p}$.
The definition of the spaces $\dot{K}_{p}^{\alpha,q}\dot{A}_{\beta}(\mathbb{R}^{n},\{t_{k}\}_{k\in \mathbb{Z}})$ is
independent of the choices of $\varphi \in \mathcal{S}(\mathbb{R}^{n})$
satisfying $\mathrm{\eqref{Ass1}}$\ and\ $\mathrm{\eqref{Ass2}}$.
\end{cor}
By adapting the arguments of Lemma 3.11 in \cite{D20.1}, we obtain the following statement.
\begin{lem}
Let $0<\theta \leq p<\infty $ and $0<q,\beta\leq\infty $. Let $\{t_{k}\}_{k\in \mathbb{Z}}$ be a $p$
-admissible weight sequence.  $\dot{K}_{p}^{\alpha,q}\dot{a}_{\beta}(\mathbb{R}^{n},\{t_{k}\}_{k\in \mathbb{Z}})$ are
quasi-Banach spaces. They are Banach spaces if $1\leq p<\infty $ and $1\leq
q,\beta\leq\infty $.
\end{lem}

Applying this lemma and Theorem \ref{phi-tran} we obtain the following
useful properties of these function spaces, where the  proof follows closely the arguments in \cite{D20.1}.

\begin{thm}
Let $\upsilon =(\upsilon _{1},\upsilon _{2})\in \mathbb{R}^{2},0<\theta \leq
p<\infty $ and $0<q,\beta\leq\infty $. Let $\{t_{k}\}_{k\in \mathbb{Z}}\in \dot{X}_{\upsilon ,\sigma ,p}$
be a $p$-admissible weight sequence with $\sigma =(\sigma _{1}=\theta \left(
\frac{p}{\theta}\right) ^{\prime },\sigma _{2}\geq p)$. $\dot{K}_{p}^{\alpha,q}\dot{A}_{\beta}(\mathbb{R}^{n},\{t_{k}\}_{k\in \mathbb{Z}})$ are quasi-Banach spaces. They are Banach spaces if $1\leq
p<\infty $ and $1\leq q,\beta\leq\infty $.
\end{thm}

\begin{thm}
\label{embeddings-S-inf}Let $0<\theta \leq p<\infty $ and $0<q,\beta\leq\infty $.
Let\ $\{t_{k}\}_{k\in \mathbb{Z}}\in \dot{X}_{\upsilon ,\sigma ,p}$ be a $p$-admissible weight
sequence with $\sigma =(\sigma _{1}=\theta \left( \frac{p}{\theta} \right) ^{\prime
},\sigma _{2}\geq p)$ and $\upsilon =(\upsilon _{1},\upsilon _{2})\in \mathbb{R}
^{2}$. Let $\delta $ be as in Lemma \ref{Ap-Property} and $\alpha>-\frac{
	n\delta }{p}$. \newline
$\mathrm{(i)}$ We have the embedding 
\begin{equation}
\mathcal{S}_{\infty }(\mathbb{R}^{n})\hookrightarrow \dot{K}_{p}^{\alpha,q}\dot{A}_{\beta}(\mathbb{R}^{n},\{t_{k}\}_{k\in \mathbb{Z}}).  \label{embedding}
\end{equation}
Let $0<q,\beta<\infty $. Then $\mathcal{S}_{\infty }(\mathbb{R}^{n})$ is dense in $\dot{K}_{p}^{\alpha,q}\dot{A}_{\beta}(\mathbb{R}^{n},\{t_{k}\}_{k\in \mathbb{Z}})$.\newline
$\mathrm{(ii)}$ We have the embedding 
\begin{equation*}
\dot{K}_{p}^{\alpha,q}\dot{A}_{\beta}(\mathbb{R}^{n},\{t_{k}\}_{k\in \mathbb{Z}})\hookrightarrow \mathcal{S}_{\infty
}^{\prime }(\mathbb{R}^{n}).
\end{equation*}
\end{thm}

\begin{proof}
The proof is a variant of the arguments used in \cite{D20} and \cite{D20.1}, together with Lemma \ref{key-estimate1}.
\end{proof}
\section{Embeddings}

For our spaces introduced above we want to show Sobolev embedding theorems.
We recall that a quasi-Banach space $A_{1}$ is continuously embedded in
another quasi-Banach space $A_{2}$, $A_{1}\hookrightarrow A_{2}$, if $%
A_{1}\subset A_{2}$ and there is a $c>0$ such that 
\begin{equation*}
\big\|f\big\|_{A_{2}}\leq c\big\|f\big\|_{A_{1}}
\end{equation*}
for all $f\in A_{1}$. We begin with the following elementary embeddings.

The following theorem establishes the basic embeddings of the spaces  $\dot{K}_{p}^{\alpha,q}\dot{A}_{\beta}(\mathbb{R}^{n},\{t_{k}\}_{k\in \mathbb{Z}})$.  Its proof is very similar to those for  Herz-type Besov and Triebel-Lizorkin spaces; see \cite{Drihem1.13}

\begin{thm}
	\label{embeddings1.1}\textit{Let $\upsilon =(\upsilon _{1},\upsilon _{2})\in \mathbb{R}
		^{2},0<\theta \leq p<\infty $ and $0<q,\beta\leq\infty $. Let $\{t_{k}\}_{k\in \mathbb{Z}}\in \dot{X}
		_{\upsilon ,\sigma ,p}$ be a $p$-admissible weight sequence with $\sigma
		=(\sigma _{1}=\theta \left( \frac{p}{\theta} \right) ^{\prime },\sigma _{2}\geq p)$. Let $\delta $ be as in Lemma \ref{Ap-Property} and $\alpha>-\frac{
			n\delta }{p}$.}\newline
	$\mathrm{(i)}$\textit{\ If }$0<\beta _{1}\leq\beta _{2}\leq \infty $, then%
	\begin{equation*}
	\dot{K}_{p}^{\alpha,q}\dot{A}_{\beta_{1}}(\mathbb{R}^{n},\{t_{k}\}_{k\in \mathbb{Z}})\hookrightarrow \dot{K}_{p}^{\alpha,q}\dot{A}_{\beta_{2}}(\mathbb{R}^{n},\{t_{k}\}_{k\in \mathbb{Z}}).  \label{embed2}
	\end{equation*}%
	$\mathrm{(ii)}$\textit{\ If }$0<q_{1}\leq q_{2}\leq \infty $, then%
	\begin{equation*}
	\dot{K}_{p}^{\alpha,q_{1}}\dot{A}_{\beta}(\mathbb{R}^{n},\{t_{k}\}_{k\in \mathbb{Z}})\hookrightarrow \dot{K}_{p}^{\alpha,q_{2}}\dot{A}_{\beta}(\mathbb{R}^{n},\{t_{k}\}_{k\in \mathbb{Z}}).  \label{embed3}
	\end{equation*}%
	$\mathrm{(iii)}$\textit{\ If }$0<p_{1}\leq p_{2}\leq \infty $, then%
	\begin{equation*}
	\dot{K}_{p_{2}}^{\alpha,q}\dot{A}_{\beta}(\mathbb{R}^{n},\{t_{k}\}_{k\in \mathbb{Z}})\hookrightarrow \dot{K}_{p_{1}}^{r,q}\dot{A}_{\beta}(\mathbb{R}^{n},\{t_{k}\}_{k\in \mathbb{Z}})
	\end{equation*}%
	where $r=\alpha -n\big(\frac{1}{p_{1}}-\frac{1}{p_{2}}\big).$	\newline
	$\mathrm{(iv)}$ We have
	\begin{equation*}
	\dot{K}_{p}^{\alpha,q}\dot{B}_{\min \left( \beta ,p,q\right)}(\mathbb{R}^{n},\{t_{k}\}_{k\in \mathbb{Z}})\hookrightarrow 
	\dot{K}_{p}^{\alpha,q}\dot{F}_{\beta}(\mathbb{R}^{n},\{t_{k}\}_{k\in \mathbb{Z}})\hookrightarrow \dot{K}_{p}^{\alpha,q}\dot{A}_{\min \left( \beta ,p,q\right)}(\mathbb{R}^{n},\{t_{k}\}_{k\in \mathbb{Z}}).  \label{aux7}
	\end{equation*}
\end{thm}

The following lemma is the $\dot{K}_{p}^{\alpha ,q}$-version of the
Plancherel-Polya-Nikolskij inequality; see \cite{D20}.

\begin{lem}
	\label{Bernstein-Herz-ine1}\textit{Let }$R>0,\alpha _{1},\alpha _{2}\in 
	\mathbb{R}\mathit{\ }$\textit{and} $0<p_{1},p_{2},r\leq \infty $. \textit{Suppose
		that }%
	\begin{equation*}
	\alpha _{1}+\frac{n}{p_{1}}>0,\quad 0<p_{2}\leq p_{1}\leq \infty \quad \text{and\quad }%
	\alpha _{2}\geq \alpha _{1}.
	\end{equation*}%
	\textit{Then there exists a positive constant }$c>0$\textit{\ independent of 
	}$R$\textit{\ such that for all }$f\in \dot{K}_{p_{2}}^{\alpha _{2},\mu }\cap 
	\mathcal{S}^{\prime }(\mathbb{R}^{n})$\textit{\ with }$\mathrm{\ supp}$%
	\textit{\ }$\mathcal{F}f\subset \{\xi \in \mathbb{R}^{n}:|\xi |\leq R\}$%
	\textit{, we have} 
	\begin{equation*}
	\big\|f\big\|_{\dot{K}_{p_{1}}^{\alpha _{1},r}}\leq c\text{ }R^{\frac{n}{p_{2}}-%
		\frac{n}{p_{1}}+\alpha _{2}-\alpha _{1}}\big\|f\big\|_{\dot{K}_{p_{2}}^{\alpha
			_{2},\mu }},
	\end{equation*}%
	where 
	\begin{equation*}
	\mu =\left\{ 
	\begin{array}{ccc}
	r, & \text{if} & \alpha _{2}=\alpha _{1}, \\ 
	\infty, & \text{if} & \alpha _{2}>\alpha _{1}.%
	\end{array}%
	\right.
	\end{equation*}
\end{lem}
In the previous lemma we have not treated the case $p_{1}< p_{2}$. The next lemma
gives a positive answer; see \cite{D20}.

\begin{lem}
	\label{Bernstein-Herz-ine2}\textit{Let }$R>0,\alpha _{1},\alpha _{2}\in 
	\mathbb{R}\mathit{\ }$\textit{and} $0<p,p_{1},p_{2},r\leq \infty $. \textit{Suppose
		that }$\alpha _{1}+\frac{n}{p_{1}}>0,0<p_{1}\leq p_{2}\leq \infty $. Assume that 
	\begin{equation*}
	\alpha _{2}>\alpha _{1}+\frac{n}{p_{1}}-\frac{n}{p_{2}}\quad \text{or}\quad \alpha
	_{2}=\alpha _{1}+\frac{n}{p_{1}}-\frac{n}{p_{2}}\quad \text{and}\quad r=p.
	\end{equation*}%
	\textit{Then there exists a positive constant }$c$\textit{\ independent of }$%
	R$\textit{\ such that for all }$f\in \dot{K}_{p_{2}}^{\alpha _{2},p}\cap 
	\mathcal{S}^{\prime }(\mathbb{R}^{n})$\textit{\ with }$\mathrm{supp}$\textit{%
		\ }$\mathcal{F}f\subset \{\xi \in \mathbb{R}^{n}:|\xi |\leq R\}$\textit{, we
		have} 
	\begin{equation*}
	\big\|f\big\|_{\dot{K}_{p_{1}}^{\alpha _{1},r}}\leq c\text{ }R^{\frac{n}{p_{2}}-%
		\frac{n}{p_{1}}+\alpha _{2}-\alpha _{1}}\big\|f\big\|_{\dot{K}_{p_{2}}^{\alpha
			_{2},p}}. 
	\end{equation*}
\end{lem}
\begin{rem} \label{Bernstein-Herz-ine3}
In Lemmas \ref{Bernstein-Herz-ine1} and \ref{Bernstein-Herz-ine2} we can replace $\big\| f\big\|_{\dot{K}_{p_{1}}^{\alpha _{1},r}}$ by $\big\|(\eta _{R,N}\ast |f|^{d})^{\frac{1}{d}}\big\|_{\dot{K}_{p_{1}}^{\alpha _{1},r}}$ where $\eta _{R,N}=R^{n}(1+R\left\vert \cdot \right\vert )^{-N}$, $N $ large enough and $0<d<\min \big(p_{2},\frac{n}{\frac{n}{p_{2}}+\alpha }\big).$
\end{rem}

The first  main result of this subsection is the following Sobolev-type embeddings. In the classical setting this was done in  \cite{T1}. We set%
\begin{equation*}
w_{k,Q}(p_{1})=\Big(\int_{Q}w_{k}^{p_{1}}(x)dx\Big)^{1/p_{1}}\quad \text{
	and\quad }t_{k,Q}(p_{0})=\Big(\int_{Q}t_{k}^{p_{0}}(x)dx\Big)^{1/p_{0}},
\end{equation*}%
where $Q\in \mathcal{Q}$ with $\ell (Q)=2^{-k}, k \in \mathbb{Z}.$

We begin by the  Sobolev embeddings for the spaces $\dot{K}_{p}^{\alpha,q}\dot{b}_{\beta}(\mathbb{R}^{n},\{t_{k}\}_{k\in \mathbb{Z}})$.
\begin{thm}
	\label{Sobolev-embedding-sequence}Let $0<\theta \leq p_{1}<\infty ,0<\theta \leq p_{2}<\infty $
	and $0<r,\beta\leq\infty $. Let\ $\{t_{k}\}_{k\in \mathbb{Z}}$ be a $p_{2}$-admissible weight
	sequence satisfying $\mathrm{\eqref{Asum1}}$\ with $p=p_{2}$, $\sigma
	_{1}=\theta \left( \frac{p_{2}}{\theta }\right) ^{\prime }$ and $j=k$. Let\ $%
	\{w_{k}\}_{k\in \mathbb{Z}}$ be a $p_{1}$-admissible weight sequence satisfying $\mathrm{\ %
		\eqref{Asum1}}$\ with $p=p_{1}$, $\sigma _{1}=\theta \left( \frac{p_{1}}{\theta}
	\right) ^{\prime }$ and $j=k$. Let $\delta_{\{t_{k}\}_{k\in \mathbb{Z}}} $ be as in Lemma \ref{Ap-Property} and $\alpha_{2}>-\frac{
		n\delta_{\{t_{k}\}_{k\in \mathbb{Z}}} }{p_{2}}$. Let $\delta_{\{w_{k}\}_{k\in \mathbb{Z}}} $ be as in Lemma \ref{Ap-Property} and $\alpha_{1}>-\frac{
		n\delta_{\{w_{k}\}_{k\in \mathbb{Z}}} }{p_{1}}$. Assume that 
	\begin{equation*}
	2^{-k\alpha_{1}}w_{k,Q}(p_{1})\lesssim 2^{-k\alpha_{2}} t_{k,Q}(p_{2})
	\end{equation*}
	
	for all $Q\in \mathcal{Q}$ with $\ell (Q)=2^{-k}, k \in \mathbb{Z}$. \newline
	$\mathrm{(i)}$ Let $\alpha _{1},\alpha _{2}\in 
	\mathbb{R}$ and $0<\theta \leq p_{2}<p_{1}<\infty$. \textit{Suppose
		that }%
	\begin{equation*}
	\alpha _{1}+\frac{n}{p_{1}}>0,\quad 0<p_{2}\leq p_{1}\leq \infty \quad \text{and\quad }%
	\alpha _{2}\geq \alpha _{1}.
	\end{equation*}%
	Then we have
	\begin{equation*}
	\dot{K}_{p_{2}}^{\alpha_{2},\mu}\dot{b}_{\beta}(\mathbb{R}^{n},\{t_{k}\}_{k\in \mathbb{Z}})\hookrightarrow \dot{K}_{p_{1}}^{\alpha_{1},r}\dot{b}_{\beta}(\mathbb{R}^{n},\{w_{k}\}_{k\in \mathbb{Z}}).
	\end{equation*}
	where 
	\begin{equation*}
	\mu =\left\{ 
	\begin{array}{ccc}
	r, & \text{if} & \alpha _{2}=\alpha _{1}, \\ 
	\infty, & \text{if} & \alpha _{2}>\alpha _{1}.%
	\end{array}%
	\right.
	\end{equation*}
	\newline
	$\mathrm{(ii)}$ Let $\alpha _{1},\alpha _{2}\in 
	\mathbb{R}$, $0<\theta \leq p_{1}<p_{2}<\infty $
	and $0<p_{3}\leq\infty $. Assume that 
	\begin{equation*}
	\alpha _{2}>\alpha _{1}+\frac{n}{p_{1}}-\frac{n}{p_{2}}\quad \text{or}\quad \alpha
	_{2}=\alpha _{1}+\frac{n}{p_{1}}-\frac{n}{p_{2}}\quad \text{and}\quad r=p_{3}.
	\end{equation*}%
	Then we have
	\begin{equation*}
	\dot{K}_{p_{2}}^{\alpha_{2},p_{3}}\dot{b}_{\beta}(\mathbb{R}^{n},\{t_{k}\}_{k\in \mathbb{Z}})\hookrightarrow \dot{K}_{p_{1}}^{\alpha_{1},r}\dot{b}_{\beta}(\mathbb{R}^{n},\{w_{k}\}_{k\in \mathbb{Z}}).
	\end{equation*}

\end{thm}
\begin{proof}
	By similarity, we prove only (i). Let  $\lambda =\{\lambda _{k,m}\}_{k\in \mathbb{Z},m\in \mathbb{Z}^{n}}\in \dot{K}_{p_{2}}^{\alpha_{2},\mu}\dot{b}_{\beta}(\mathbb{R}^{n},\{t_{k}\}_{k\in \mathbb{Z}})$. We have
	\begin{align*}
	|\lambda _{k,m}|^{u}=& 2^{nk}\int_{\mathbb{R}^{n}}|\lambda _{k,m}|^{u}\chi
	_{k,m}(y)dy \\
	\lesssim & \eta _{k,N}\ast(|\lambda _{k,m}|^{u}\chi _{k,m})
	\end{align*}
	for any $x\in Q_{k,m} $ and any $N,u>0$, where $\eta _{k,N}=2^{nk}(1+2^{k}\left\vert \cdot \right\vert )^{-N}$.  This leads to
	
	\begin{align*}
	\sum_{m\in \mathbb{Z}^{n}}|\lambda _{k,m}|^{u}(w_{k,m,p_{1}})^{u}\chi _{k,m}
	\lesssim & \eta _{k,N}\ast(\sum_{m\in \mathbb{Z}^{n}}|\lambda _{k,m}|^{u}(w_{k,m,p_{1}})^{u}\chi _{k,m})
	\end{align*}
Choose  $0<u<\min \big(p_{2},\frac{n}{\frac{n}{p_{2}}+\alpha }\big)$ and $N$ large enough.	Applying Lemma 	\ref{Bernstein-Herz-ine1} combined with Remark 	\ref{Bernstein-Herz-ine3}, we obtain
		\begin{align*}
	\Big\| \sum\limits_{m\in \mathbb{Z}%
		^{n}}|\lambda _{k,m}|w_{k,m, p_{1}}\chi _{k,m}\Big\|_{\dot{K}_{p_{1}}^{\alpha_{1},r}}
		\end{align*}
		can be estimated by 
		\begin{align*}
	c2^{k(\frac{n}{p_{2}}-\frac{n}{p_{1}})}	\Big\| \sum\limits_{m\in \mathbb{Z}%
			^{n}}|\lambda _{k,m}|t_{k,m, p_{1}}\chi _{k,m}\Big\|_{\dot{K}_{p_{2}}^{\alpha_{2},\mu}}
		\end{align*}
		for any $k\in \mathbb{Z}$ where the positive constant $c>0$ is independent of $k$. This leads to desired estimate.
\end{proof}	
From Theorems {\ref{phi-tran} and \ref{Sobolev-embedding-sequence}}, we
infer the following Sobolev-type embedding for the spaces $\dot{K}_{p}^{\alpha,q}\dot{B}_{\beta}(\mathbb{R}^{n},\{t_{k}\}_{k\in \mathbb{Z}})$.
\begin{thm}
	\label{Sobolev-embedding-sequenceherz2}
	Let $0<\theta \leq p_{1}<\infty ,0<\theta \leq p_{2}<\infty $
	and $0<r,\beta\leq\infty $. Let\ $\{t_{k}\}_{k\in \mathbb{Z}}\in \dot{X}_{\upsilon _{2},\sigma ,p_{2}} $be a $p_{2}$-admissible weight
	sequence  with $\sigma =(\sigma _{1}=\theta \left(
	\frac{p_{2}}{\theta}\right) ^{\prime },\sigma _{2}\geq p_{2})$ and $\upsilon
	_{2}=(\upsilon _{1,2},\upsilon _{2,2})\in \mathbb{R}^{2}$.
	Let\ $\{w_{k}\}_{k\in \mathbb{Z}}$ $\in 
	\dot{X}_{\upsilon _{1},\sigma ,p_{1}}$ be a $p_{1}$-admissible weight sequence
	with $\sigma =(\sigma _{1}=\theta \left( \frac{p_{1}}{\theta} \right) ^{\prime
	},\sigma _{2}\geq p_{1})$ and $\upsilon _{1}=(\upsilon _{1,1},\upsilon _{2,1})\in 
	\mathbb{R}^{2}$.
	Let $\delta_{\{t_{k}\}_{k\in \mathbb{Z}}} $ be as in Lemma \ref{Ap-Property} and $\alpha_{2}>-\frac{
		n\delta_{\{t_{k}\}_{k\in \mathbb{Z}}} }{p_{2}}$. Let $\delta_{\{w_{k}\}_{k\in \mathbb{Z}}} $ be as in Lemma \ref{Ap-Property} and $\alpha_{1}>-\frac{
		n\delta_{\{w_{k}\}_{k\in \mathbb{Z}}} }{p_{1}}$.  Assume that 
	\begin{equation*}
	2^{-k\alpha_{1}}w_{k,Q}(p_{1})\lesssim 2^{-k\alpha_{2}} t_{k,Q}(p_{2})
	\end{equation*}
	
	for all $Q\in \mathcal{Q}$ with $\ell (Q)=2^{-k}, k \in \mathbb{Z}$. \newline
	$\mathrm{(i)}$ Let $0<\theta \leq p_{2}<p_{1}<\infty $. \textit{Suppose
		that }%
	$\alpha _{2}\geq \alpha _{1}$.

	Then we have
	\begin{equation*}
	\dot{K}_{p_{2}}^{\alpha_{2},\mu}\dot{B}_{\beta}(\mathbb{R}^{n},\{t_{k}\}_{k\in \mathbb{Z}})\hookrightarrow \dot{K}_{p_{1}}^{\alpha_{1},r}\dot{B}_{\beta}(\mathbb{R}^{n},\{w_{k}\}_{k\in \mathbb{Z}}).
	\end{equation*}
	where 
	\begin{equation*}
	\mu =\left\{ 
	\begin{array}{ccc}
	r, & \text{if} & \alpha _{2}=\alpha _{1}, \\ 
	\infty, & \text{if} & \alpha _{2}>\alpha _{1}.%
	\end{array}%
	\right.
	\end{equation*}
	\newline
	$\mathrm{(ii)}$ Let  $0<\theta \leq p_{1}<p_{2}<\infty $
	and $0<p_{3}\leq\infty $. Assume that 
	\begin{equation*}
	\alpha _{2}>\alpha _{1}+\frac{n}{p_{1}}-\frac{n}{p_{2}}\quad \text{or}\quad \alpha
	_{2}=\alpha _{1}+\frac{n}{p_{1}}-\frac{n}{p_{2}}\quad \text{and}\quad r=p_{3}.
	\end{equation*}%
	Then we have
	\begin{equation*}
	\dot{K}_{p_{2}}^{\alpha_{2},p_{3}}\dot{B}_{\beta}(\mathbb{R}^{n},\{t_{k}\}_{k\in \mathbb{Z}})\hookrightarrow \dot{K}_{p_{1}}^{\alpha_{1},r}\dot{B}_{\beta}(\mathbb{R}^{n},\{w_{k}\}_{k\in \mathbb{Z}}).
	\end{equation*}

\end{thm}
From Theorem \ref{Sobolev-embedding-sequenceherz2}\ and the fact that $\dot{K}_{p}^{0,p}\dot{B}_{\beta}(\mathbb{R}^{n},\{t_{k}\}_{k\in \mathbb{Z}})=\dot{B}_{p,\beta }(\mathbb{R}^{n},\{t_{k}\}_{k\in \mathbb{Z}})$\ we immediately arrive at the following
result.
\begin{thm}
	\label{Sobolev-embedding-sequenceherz3}Let $\theta, p_{2},p_{1},\beta  $ $\{t_{k}\}_{k\in \mathbb{Z}}$ and $ \{w_{k}\}_{k\in \mathbb{Z}},\alpha_{2} $ be as in Theorem  \ref{Sobolev-embedding-sequence}.  Assume that 
	\begin{equation*}
	w_{k,Q}(p_{1})\lesssim 2^{-k\alpha_{2}} t_{k,Q}(p_{2})
	\end{equation*}
	
	for all $Q\in \mathcal{Q}$ with $\ell (Q)=2^{-k}, k \in \mathbb{Z}$. \newline
	$\mathrm{(i)}$ Let  $0<\theta \leq p_{2}<p_{1}<\infty $. \textit{Suppose
		that $\alpha _{2}\geq 0$}.Then we have
	\begin{equation*}
	\dot{K}_{p_{2}}^{\alpha_{2},\mu}\dot{B}_{\beta}(\mathbb{R}^{n},\{t_{k}\}_{k\in \mathbb{Z}})\hookrightarrow \dot{B}_{p_{1},\beta}(\mathbb{R}^{n},\{w_{k}\}_{k\in \mathbb{Z}}).
	\end{equation*}
	where 
	\begin{equation*}
	\mu =\left\{ 
	\begin{array}{ccc}
	p_{1}, & \text{if} & \alpha _{2}=0, \\ 
	\infty, & \text{if} & \alpha _{2}>0.%
	\end{array}%
	\right.
	\end{equation*}
	\newline
	$\mathrm{(ii)}$ Let $0<\theta \leq p_{1}<p_{2}<\infty $
	and $0<p_{3}\leq\infty $. Assume that 
	\begin{equation*}
	\alpha _{2}>\frac{n}{p_{1}}-\frac{n}{p_{2}}\quad \text{or}\quad \alpha
	_{2}=\frac{n}{p_{1}}-\frac{n}{p_{2}}\quad \text{and}\quad p_{1}=p_{3}.
	\end{equation*}%
	Then we have
	\begin{equation*}
	\dot{K}_{p_{2}}^{\alpha_{2},p_{3}}\dot{B}_{\beta}(\mathbb{R}^{n},\{t_{k}\}_{k\in \mathbb{Z}})\hookrightarrow \dot{B}_{p_{1},\beta}(\mathbb{R}^{n},\{w_{k}\}_{k\in \mathbb{Z}}).
	\end{equation*}
	
\end{thm}
\begin{cor}
	Let $\theta, p_{2},p_{1},\beta, $ $\{t_{k}\}_{k\in \mathbb{Z}}, \{w_{k}\}_{k\in \mathbb{Z}}$ be as in Theorem  \ref{Sobolev-embedding-sequence}.  Assume that 
	\begin{equation*}
	w_{k,Q}(p_{1})\lesssim t_{k,Q}(p_{2})
	\end{equation*}
	
	for all $Q\in \mathcal{Q}$ with $\ell (Q)=2^{-k}, k \in \mathbb{Z}$. Assume that $0<\theta \leq p_{2}\leq p_{1}< \infty$. Then%
	
	\begin{equation*}
	\dot{B}_{p_{2},\beta}(\mathbb{R}^{n},\{t_{k}\}_{k\in \mathbb{Z}})\hookrightarrow\dot{K}_{p_{2}}^{0,p_{1}}\dot{B}_{\beta}(\mathbb{R}^{n},\{t_{k}\}_{k\in \mathbb{Z}}) \hookrightarrow\dot{B}_{p_{1},\beta}(\mathbb{R}^{n},\{w_{k}\}_{k\in \mathbb{Z}}).
\end{equation*}
	
\end{cor}

\begin{proof}
	To prove these embeddings, it sufficient to choose in Theorem \ref%
	{Sobolev-embedding-sequenceherz3}, $\mu =p_{1}$ and $\alpha_{2} =0$.\ However the desired embeddings
	are an immediate consequence of the fact that 
	\begin{align*}
	\dot{B}_{p_{2},\beta}(\mathbb{R}^{n},\{t_{k}\}_{k\in \mathbb{Z}})=&\dot{K}_{p_{2}}^{0,p_{2}}\dot{B}_{\beta}(\mathbb{R}^{n},\{t_{k}\}_{k\in \mathbb{Z}}) \\
	 \hookrightarrow 
	  &\dot{K}_{p_{2}}^{0,p_{1}}\dot{B}_{\beta}(\mathbb{R}^{n},\{t_{k}\}_{k\in \mathbb{Z}}) \\
	\hookrightarrow
	&\dot{B}_{p_{1},\beta}(\mathbb{R}^{n},\{w_{k}\}_{k\in \mathbb{Z}}).
	\end{align*}%
	The proof is complete.
\end{proof}

Now, we present the  Sobolev embeddings for the spaces $\dot{K}_{p}^{\alpha,q}\dot{f}_{\beta}(\mathbb{R}^{n},\{t_{k}\}_{k\in \mathbb{Z}})$.

\begin{thm}
	\label{embeddings3copy(1)}
	Let $0< p_{2},p_{1}<\infty $, $0<r\leq p_{3}\leq\infty $
	and $0<\beta\leq\infty $. Let\ $\{t_{k}\}_{k\in \mathbb{Z}}$ be a $p_{2}$-admissible weight
	sequence satisfying $\mathrm{\eqref{Asum1}}$\ with $p=p_{2}$, $\sigma
	_{1}=\theta \left( \frac{p_{2}}{\theta} \right) ^{\prime }$ and $j=k$. Let\ $%
	\{w_{k}\}_{k\in \mathbb{Z}}$ be a $p_{1}$-admissible weight sequence satisfying $\mathrm{\ %
		\eqref{Asum1}}$\ with $p=p_{1}$, $\sigma _{1}=\theta \left( \frac{p_{1}}{\theta}
	\right) ^{\prime }$ and $j=k$. Let $\delta_{\{t_{k}\}_{k\in \mathbb{Z}}} $ be as in Lemma \ref{Ap-Property} and $\alpha_{2}>-\frac{
		n\delta_{\{t_{k}\}_{k\in \mathbb{Z}}} }{p_{2}}$. Let $\delta_{\{w_{k}\}_{k\in \mathbb{Z}}} $ be as in Lemma \ref{Ap-Property} and $\alpha_{1}>-\frac{
		n\delta_{\{w_{k}\}_{k\in \mathbb{Z}}} }{p_{1}}$. Assume that 
	\begin{equation}
	2^{-k\alpha_{1}}w_{k,Q}(p_{1})\lesssim 2^{-k\alpha_{2}} t_{k,Q}(p_{2})\label{newexp1.1herz}
	\end{equation}
	
	for all $Q\in \mathcal{Q}$ with $\ell (Q)=2^{-k}, k \in \mathbb{Z}$.
	\textit{Let }$0<\theta \leq p_{2}<p_{1}<\infty $ and $\alpha _{2}\geq \alpha _{1}$ or $0<\theta \leq p_{1}\leq
	p_{2}<\infty $ and 
	\begin{equation*}
	\alpha _{2}+\frac{n}{p_{2}}\geq \alpha _{1}+\frac{n}{p_{1}}.  \label{newexp1.2}
	\end{equation*}%
	Then%
	\begin{equation*}
	\dot{K}_{p_{2}}^{\alpha_{2},r}\dot{f}_{\mu}(\mathbb{R}^{n},\{t_{k}\}_{k\in \mathbb{Z}})\hookrightarrow \dot{K}_{p_{1}}^{\alpha_{1},p_{3}}\dot{f}_{\beta}(\mathbb{R}^{n},\{w_{k}\}_{k\in \mathbb{Z}}),
\end{equation*}
 where%
	\begin{equation*}
	\mu =\left\{ 
	\begin{array}{ccc}
	\beta , & \text{if} & 0<\theta \leq p_{1}\leq p_{2}<\infty \text{ and }\alpha _{2}+\frac{n}{p_{2}}%
	=\alpha _{1}+\frac{n}{p_{1}}, \\ 
	\infty , &  & \text{otherwise.}%
	\end{array}%
	\right.
	\end{equation*}
\end{thm}
\begin{proof} Following the arguments in \cite{Drihem2.13}, we provide a detailed proof for the reader's convenience.
	
\textit{Step 1.} We consider the case $0<\theta \leq p_{2}<p_{1}<\infty $ and $\alpha _{2}\geq
\alpha _{1}$. In view of the embedding $\ell ^{r}\hookrightarrow \ell ^{p_{3}}$,
it is sufficient to prove that 
\begin{equation*}
\dot{K}_{p_{2}}^{\alpha_{2},r}\dot{f}_{\infty}(\mathbb{R}^{n},\{t_{k}\}_{k\in \mathbb{Z}})\hookrightarrow \dot{K}_{p_{1}}^{\alpha_{1},r}\dot{f}_{\beta}(\mathbb{R}^{n},\{w_{k}\}_{k\in \mathbb{Z}}),
\end{equation*}
For simplicity, we only consider the case $\beta =1$. Let $\lambda \in \dot{K}_{p_{2}}^{\alpha_{2},r}\dot{f}_{\infty}(\mathbb{R}^{n},\{t_{k}\}_{k\in \mathbb{Z}})$. We have%
\begin{align}
\big\|\lambda \big\|_{\dot{K}_{p_{1}}^{\alpha_{1},r}\dot{f}_{1}(\mathbb{R}^{n},\{w_{k}\}_{k\in \mathbb{Z}})}\lesssim & %
\Big(\sum\limits_{k=-\infty }^{c_{n}+1}2^{k\alpha _{1}r}\Big\|%
\sum_{v=0}^{\infty }2^{vn(\frac{1}{2}+\frac{1}{ p_{1}})}\sum_{m\in \mathbb{Z}%
	^{n}}w_{v,m, p_{1}}|\lambda _{v,m}|\chi _{v,m}\chi _{k}\Big\|_{p_{1}}^{r}\Big)^{1/r}  \notag \\
& +\Big(\sum\limits_{k=c_{n}+2}^{\infty }2^{k\alpha _{1}r}\Big\|%
\sum_{v=0}^{\infty }2^{vn(\frac{1}{2}+\frac{1}{ p_{1}})}\sum_{m\in \mathbb{Z}%
	^{n}}w_{v,m, p_{1}}|\lambda _{v,m}|\chi _{v,m}\chi _{k}\Big\|_{p_{1}}^{r}\Big)^{1/r}.
\label{est 1}
\end{align}%
Here $c_{n}=1+\lfloor \log _{2}(2\sqrt{n}+1)\rfloor $. The first term can be
estimated \ from above by%
\begin{align*}
& c\Big(\sum\limits_{k=-\infty }^{c_{n}+1}2^{k\alpha _{1}r}\Big\|%
\sum_{v=0}^{c_{n}-k+1}2^{vn(\frac{1}{2}+\frac{1}{ p_{1}})}\sum_{m\in \mathbb{Z}%
	^{n}}w_{v,m, p_{1}}|\lambda _{v,m}|\chi _{v,m}\chi _{k}\Big\|_{p_{1}}^{r}\Big)^{1/r} \\
& +c\Big(\sum\limits_{k=-\infty }^{c_{n}+1}2^{k\alpha _{1}r}\Big\|%
\sum_{v=c_{n}-k+2}^{\infty }2^{vn(\frac{1}{2}+\frac{1}{ p_{1}})}\sum_{m\in \mathbb{Z}%
	^{n}}w_{v,m, p_{1}}|\lambda _{v,m}|\chi _{v,m}\chi _{k}\Big\|_{p_{1}}^{r}\Big)^{1/r} \\
=& I+II.
\end{align*}

\textit{Estimation of }$I$\textit{.} Let $x\in R_{k}\cap Q_{v,m}$ and $y\in
Q_{v,m}$\ such that $0\leq v\leq c_{n}-k+1$. We have 
\begin{equation*}
|x-y|\leq 2\sqrt{n}2^{-v}<2^{c_{n}-v}
\end{equation*}%
and from this it follows that 
\begin{equation*}
|y|<2^{c_{n}-v}+2^{k}\leq 2^{c_{n}-v+2},
\end{equation*}%
which implies that $y$ is located in the ball $B(0,2^{c_{n}-v+2})$. Then%
\begin{align*}
(w_{v,m, p_{1}})^{h}|\lambda _{v,m}|^{h} =&2^{nv}\int_{\mathbb{R}^{n}}(w_{v,m, p_{1}})^{h}|\lambda _{v,m}|^{h}\chi
_{v,m}(y)dy \\
\leq &2^{nv}\int_{B(0,2^{c_{n}-v+2})}(w_{v,m, p_{1}})^{h}|\lambda _{v,m}|^{h}\chi _{v,m}(y)dy,
\end{align*}%
if $x\in R_{k}\cap Q_{v,m}$ and $h>0$. Therefore for any $x\in R_{k}$, we
obtain that 
\begin{align*}
\sum_{m\in \mathbb{Z}^{n}}(w_{v,m, p_{1}})^{h}|\lambda _{v,m}|^{h}\chi _{v,m}(x)\leq &
2^{nv}\int_{B(0,2^{c_{n}-v+2})}\sum_{m\in \mathbb{Z}^{n}}(w_{v,m, p_{1}})^{h}|\lambda
_{v,m}|^{h}\chi _{v,m}(y)dy \\
=& 2^{nv}\Big\|\sum_{m\in \mathbb{Z}^{n}}w_{v,m, p_{1}}|\lambda _{v,m}|\chi _{v,m}\chi
_{B(0,2^{c_{n}-v+2})}\Big\|_{h}^{h}.
\end{align*}%
This yields%
\begin{align*}
& 2^{\alpha _{1}k}\Big\|\sum_{v=0}^{c_{n}-k+1}2^{vn(\frac{1}{2}+\frac{1}{ p_{1}})}\sum_{m\in \mathbb{Z}^{n}}w_{v,m, p_{1}}|\lambda _{v,m}|\chi _{v,m}\chi _{k}\Big\|_{p_{1}} \\
\lesssim & 2^{(\alpha _{1}+\frac{n}{p_{1}})k}\sum_{v=0}^{c_{n}-k+1}2^{vn(\frac{1}{2}+\frac{1}{ p_{1}})+\frac{nv}{h}}\Big\|\sum_{m\in \mathbb{Z}^{n}}w_{v,m, p_{1}}|\lambda
_{v,m}|\chi _{v,m}\chi _{B(0,2^{c_{n}-v+2})}\Big\|_{h},
\end{align*}%
where the implicit\ constant is independent of $k$. We may choose $h>0$ such
that $\frac{1}{h}>\max (\frac{1}{p_{2}},\frac{1}{p_{2}}+\frac{\alpha _{2}}{n})$.
Using \eqref{newexp1.1herz} and Lemma \ref{lem:lq-inequality} we estimate $I^{r}$
by%
\begin{align*}
& c\sum_{v=0}^{\infty }2^{v(\frac{n}{h}-\alpha _{1}+\frac{n%
	}{2})r}\Big\|\sum_{m\in \mathbb{Z}^{n}}w_{v,m, p_{1}}|\lambda _{v,m}|\chi _{v,m}\chi
_{B(0,2^{c_{n}-v+2})}\Big\|_{h}^{r} \\
\leq & c\sum_{v=0}^{\infty }2^{v(\frac{n}{h}-\alpha _{2}+%
	\frac{n}{2})r}\Big(\sum_{i\leq -v}\Big\|\sum_{m\in \mathbb{Z}^{n}}t_{v,m, p_{2}}|\lambda
_{v,m}|\chi _{v,m}\chi _{i+c_{n}+2}\Big\|_{h}^{\varkappa }\Big)^{r/\varkappa
} \\
\leq & c\sum_{v=0}^{\infty }2^{v\frac{nr}{d}}\Big(\sum_{i\leq -v}2^{i(\frac{n%
	}{d}+\alpha _{2})\varkappa }\Big\|\sup_{j\in \mathbb{N}_{0}}2^{(\frac{n}{p_{2}}+\frac{n%
	}{2})j}\sum_{m\in \mathbb{Z}^{n}}t_{v,m, p_{2}}|\lambda _{j,m}|\chi _{j,m}\chi _{i+c_{n}+2}%
\Big\|_{p_{2}}^{\varkappa }\Big)^{r/\varkappa },
\end{align*}%
by H\"{o}lder's inequality, with $\varkappa =\min (1,h)$ and $\frac{n}{d}=%
\frac{n}{h}-\frac{n}{p_{2}}-\alpha _{2}$. Again, by applying Lemma \ref%
{lem:lq-inequality} we get%
\begin{equation*}
I^{r}\leq c\sum_{i=0}^{\infty }2^{-\alpha _{2}ir}\Big\|\sup_{j\in \mathbb{N}%
	_{0}}2^{(\frac{n%
	}{2}+\frac{n}{p_{2}})j}\sum_{m\in \mathbb{Z}^{n}}t_{v,m, p_{2}}|\lambda _{j,m}|\chi
_{j,m}\chi _{-i+c_{n}+2}\Big\|_{p_{2}}^{r}\leq c\big\|\lambda \big\|_{\dot{K}_{p_{2}}^{\alpha_{2},r}\dot{f}_{\infty}(\mathbb{R}^{n},\{t_{k}\}_{k\in \mathbb{Z}})}^{r}.
\end{equation*}

\textit{Estimation of }$II$\textit{.} We see that it suffices to show that
for any $k\leq c_{n}+1$%
\begin{align*}
& 2^{k\alpha _{1}}\Big\|\sum_{v=c_{n}-k+2}^{\infty }2^{vn(\frac{1}{2}+\frac{1}{ p_{1}})}\sum_{m\in \mathbb{Z}^{n}}w_{v,m, p_{1}}|\lambda _{v,m}|\chi _{v,m}\chi _{k}\Big\|_{p_{1}} \\
\leq & C\text{ }2^{k\alpha _{2}}\Big\|\sup_{v\geq c_{n}-k+2}2^{vn(\frac{1}{2}+\frac{1}{ p_{2}})}\sum_{m\in \mathbb{Z}^{n}}t_{v,m, p_{2}}|\lambda _{v,m}|\chi _{v,m}\chi _{%
	\widetilde{R}_{k}}\Big\|_{p_{2}} \\
=& \delta ,
\end{align*}%
for some suitable constant $C$, where\ $\widetilde{R}_{k}=\{x\in \mathbb{R}%
^{n}:2^{k-2}\leq\left\vert x\right\vert \leq2^{k+2}\}$. This claim can be
reformulated as showing that%
\begin{equation}
\int_{R_{k}}2^{k\alpha _{1}p_{1}}\delta ^{-p_{1}}\Big(\sum_{v=c_{n}-k+2}^{\infty
}2^{vn(\frac{1}{2}+\frac{1}{ p_{1}})}\sum_{m\in \mathbb{Z}^{n}}w_{v,m, p_{1}}|\lambda _{v,m}|\chi
_{v,m}(x)\Big)^{p_{1}}dx\leq 1.  \label{sobolevfor}
\end{equation}%
The left-hand side of \eqref{sobolevfor}\textit{\ }can be rewritten as 
\begin{align}
& \int_{R_{k}}\delta ^{-p_{1}}\Big(\sum_{v=c_{n}-k+2}^{\infty }2^{v(\frac{n}{p_{1}}-%
	\frac{n}{p_{2}})+(\alpha _{1}-\alpha _{2})(v+k)+v(\frac{n}{p_{2}}+\alpha _{2}-\alpha _{1}+\frac{n}{2})+k\alpha
	_{2}}\sum_{m\in \mathbb{Z}^{n}}w_{v,m, p_{1}}|\lambda _{v,m}|\chi _{v,m}(x)\Big)^{p_{1}}dx 
\notag \\
\leq & \int_{R_{k}}\delta ^{-p_{1}}\Big(\sum_{v=c_{n}-k+2}^{\infty }2^{v(\frac{n%
	}{p_{1}}-\frac{n}{p_{2}})+v(\frac{n}{p_{2}}+\frac{n}{2})+k\alpha _{2}}\sum_{m\in \mathbb{Z}%
	^{n}}t_{v,m, p_{2}}|\lambda _{v,m}|\chi _{v,m}(x)\Big)^{p_{1}}dx  \label{sumover} \\
=& T_{k},  \notag
\end{align}%
since $\alpha _{2}\geq \alpha _{1}$. Let us prove that $T_{k}\leq 1$ for any 
$k\leq c_{n}+1$. Our estimate use partially some decomposition techniques
already used in \cite[Theorem 3.1]{V}. We distinguish two cases.

\textit{Case 1.} $\sup_{v\geq c_{n}-k+2}2^{v(\frac{n}{p_{2}}+\frac{n}{2})+k\alpha _{2}}\sum_{m\in \mathbb{Z}%
	^{n}}t_{v,m, p_{2}}|\lambda _{v,m}|\chi _{v,m}(x)\leq \delta $.
In this case we obtain%
\begin{align*}
T_{k}& \leq C^{p_{1}}\int_{R_{k}}\Big(\delta ^{-1}\sup_{v\geq c_{n}-k+2}2^{v(\frac{n}{p_{2}}+\frac{n}{2})+k\alpha _{2}}\sum_{m\in \mathbb{Z}%
	^{n}}t_{v,m, p_{2}}|\lambda _{v,m}|\chi _{v,m}(x)\Big)^{p_{2}}dx \\
& \leq 1.
\end{align*}

\textit{Case 2.} $\sup_{v\geq c_{n}-k+2}2^{v(\frac{n}{p_{2}}+\frac{n}{2})+k\alpha _{2}}\sum_{m\in \mathbb{Z}%
	^{n}}t_{v,m, p_{2}}|\lambda _{v,m}|\chi _{v,m}(x)>\delta $. We
can distinguish two cases as follows:

$\bullet $ $\delta ^{-1}\sup_{v\geq c_{n}-k+2}2^{v(\frac{n}{p_{2}}+\frac{n}{2})+k\alpha _{2}}\sum_{m\in \mathbb{Z}%
	^{n}}t_{v,m, p_{2}}|\lambda _{v,m}|\chi _{v,m}(x)=\infty $, then there is nothing to prove.

$\bullet $ $\delta <\sup_{v\geq c_{n}-k+2}2^{v(\frac{n}{p_{2}}+\frac{n}{2})+k\alpha _{2}}\sum_{m\in \mathbb{Z}%
	^{n}}t_{v,m, p_{2}}|\lambda _{v,m}|\chi _{v,m}(x)<\infty $. Let $%
N\in \mathbb{N}$ be such that 
\begin{equation*}
2^{\frac{nN}{p_{2}}}\leq \delta ^{-1}\sup_{v\geq c_{n}-k+2}2^{v(\frac{n}{p_{2}}+\frac{n}{2})+k\alpha _{2}}\sum_{m\in \mathbb{Z}%
	^{n}}t_{v,m, p_{2}}|\lambda _{v,m}|\chi _{v,m}(x)<2^{%
	\frac{(N+1)n}{p_{2}}}.
\end{equation*}

\textit{Subcase 2.1.} $k\geq c_{n}-N+2$. We split the sum in \eqref{sumover}%
\ into two parts,%
\begin{equation*}
\sum_{v=c_{n}-k+2}^{\infty }2^{v(\frac{n%
	}{p_{1}}-\frac{n}{p_{2}})+v(\frac{n}{p_{2}}+\frac{n}{2})+k\alpha _{2}}\sum_{m\in \mathbb{Z}%
	^{n}}t_{v,m, p_{2}}|\lambda _{v,m}|\chi _{v,m}=J_{1,k}+J_{2,k},
\end{equation*}%
where 
\begin{equation*}
J_{1,k}=\sum_{v=c_{n}-k+2}^{N}2^{v(\frac{n%
	}{p_{1}}-\frac{n}{p_{2}})+v(\frac{n}{p_{2}}+\frac{n}{2})+k\alpha _{2}}\sum_{m\in \mathbb{Z}%
	^{n}}t_{v,m, p_{2}}|\lambda _{v,m}|\chi _{v,m}
\end{equation*}%
and%
\begin{equation*}
J_{2,k}=\sum_{v=N+1}^{\infty }2^{v(\frac{n%
	}{p_{1}}-\frac{n}{p_{2}})+v(\frac{n}{p_{2}}+\frac{n}{2})+k\alpha _{2}}\sum_{m\in \mathbb{Z}%
	^{n}}t_{v,m, p_{2}}|\lambda _{v,m}|\chi _{v,m}.
\end{equation*}%
Let $x\in R_{k}\cap Q_{v,m}$ and $y\in Q_{v,m}$ with $v\geq c_{n}-k+2$.\ We
have $|x-y|\leq 2\sqrt{n}2^{-v}<2^{c_{n}-v}$ and from this it follows that $%
2^{k-2}<|y|<2^{c_{n}-v}+2^{k}<2^{k+2}$, which implies that $y$ is located in 
$\widetilde{R}_{k}$. Then%
\begin{align*}
(t_{v,m, p_{2}})^{p_{2}}|\lambda _{v,m}|^{p_{2}}=& 2^{nv}\int_{\mathbb{R}^{n}}(t_{v,m, p_{2}})^{p_{2}}|\lambda _{v,m}|^{p_{2}}\chi
_{v,m}(y)dy \\
\leq & 2^{nv}\int_{\widetilde{R}_{k}}(t_{v,m, p_{2}})^{p_{2}}|\lambda _{v,m}|^{p_{2}}\chi _{v,m}(y)dy,
\end{align*}%
if $x\in R_{k}\cap Q_{v,m}$. But this immediately implies that%
\begin{align*}
\sum_{m\in \mathbb{Z}^{n}}(t_{v,m, p_{2}})^{p_{2}}|\lambda _{v,m}|^{p_{2}}\chi _{v,m}(x)\leq &
2^{nv}\int_{\widetilde{R}_{k}}\sum_{m\in \mathbb{Z}^{n}}(t_{v,m, p_{2}})^{p_{2}}|\lambda
_{v,m}|^{p_{2}}\chi _{v,m}(y)dy \\
=& 2^{nv}\big\|\sum_{m\in \mathbb{Z}^{n}}t_{v,m, p_{2}}|\lambda _{v,m}|\chi _{v,m}\chi _{%
	\widetilde{R}_{k}}\big\|_{p_{2}}^{p_{2}} \\
\leq & 2^{-\frac{np_{2}}{2}v-\alpha _{2}p_{2}k}\delta ^{p_{2}}.
\end{align*}%
Therefore for any $x\in R_{k}$, we have%
\begin{equation*}
\delta ^{-1}J_{1,k}(x)\leq \sum_{v=0}^{N}2^{\frac{n}{p_{1}}v}\leq C\text{ }2^{%
	\frac{n}{p_{1}}N}
\end{equation*}%
and 
\begin{align*}
\delta ^{-1}J_{2,k}(x)=& \delta ^{-1}\sum_{v=N+1}^{\infty }2^{v(\frac{n}{p_{1}}-%
	\frac{n}{p_{2}})+v(\frac{n}{p_{2}}+\frac{n}{2})+k\alpha _{2}}\sum_{m\in \mathbb{Z}%
	^{n}}t_{v,m, p_{2}}|\lambda _{v,m}|\chi _{v,m}(x) \\
\leq & c\text{ }\delta ^{-1}\sup_{v\geq c_{n}-k+2}2^{v(\frac{n}{p_{2}}+\frac{n}{2}%
	)+k\alpha _{2}}\sum_{m\in \mathbb{Z}^{n}}t_{v,m, p_{2}}|\lambda _{v,m}|\chi
_{v,m}(x)\sum_{v=N+1}^{\infty }2^{v(\frac{n}{p_{1}}-\frac{n}{p_{2}})} \\
\leq & C\text{ }2^{\frac{n}{p_{1}}N}.
\end{align*}%
Hence%
\begin{align*}
T_{k}& \leq C^{p_{1}}\int_{R_{k}}2^{nN}dx \\
& \leq C^{p_{1}}\int_{R_{k}}\Big(\delta ^{-1}\sup_{v\geq c_{n}-k+2}2^{v(\frac{n}{p_{2}}+\frac{n}{2})+k\alpha _{2}}\sum_{m\in \mathbb{Z}%
	^{n}}t_{v,m, p_{2}}|\lambda _{v,m}|\chi _{v,m}(x)\Big)^{p_{2}}dx \\
& \leq 1.
\end{align*}

\textit{Subcase 2.2.} $k<c_{n}-N+2$. We use the same of arguments as in
Subcase 2.1, in view of the fact that 
\begin{equation*}
\sum_{v=c_{n}-k+2}^{\infty }2^{v(\frac{n%
	}{p_{1}}-\frac{n}{p_{2}})+v(\frac{n}{p_{2}}+\frac{n}{2})+k\alpha _{2}}\sum_{m\in \mathbb{Z}%
	^{n}}t_{v,m, p_{2}}|\lambda _{v,m}|\chi _{v,m}\leq
J_{2,k}.
\end{equation*}

\textit{Estimate of \eqref{est 1}}. The arguments here are quite similar to
those used in the estimation of $II$. This complete the proof of the first
case.

\textit{Step 2.} We consider the case $0<\theta \leq p_{1}\leq p_{2}<\infty $ and $\alpha _{2}+%
\frac{n}{p_{2}}>\alpha _{1}+\frac{n}{p_{1}}$. We only need to estimate the part $%
T_{k}$. H\"{o}lder's inequality implies that%
\begin{align*}
T_{k}\lesssim & \Big\|\delta ^{-1}\sum_{v=c_{n}-k+2}^{\infty }2^{(\frac{n}{p_{1}}-%
	\frac{n}{p_{2}}+\alpha _{1}-\alpha _{2})(v+k)+v(\frac{n}{p_{2}}+\alpha _{2}-\alpha _{1}+\frac{n}{2})+k\alpha
	_{2}}\sum_{m\in \mathbb{Z}^{n}}w_{v,m, p_{1}}|\lambda _{v,m}|\chi _{v,m}\Big\|%
_{p_{2}}^{p_{1}} \\
\lesssim & \Big\|\delta ^{-1}\sup_{v\geq c_{n}-k+2}\sum_{m\in \mathbb{Z}%
	^{n}}2^{v(\frac{n}{p_{2}}+\frac{n}{2})+k\alpha _{2}}t_{v,m, p_{2}}|\lambda _{v,m}|\chi _{v,m}\chi
_{k}\Big\|_{p_{2}}^{p_{1}} \\
\lesssim 1& ,
\end{align*}%
where the last inequality follows by the fact that $\alpha _{2}+\frac{n}{q}%
>\alpha _{1}+\frac{n}{s}$. The remaining case can be easily solved. The
proof is complete.
\end{proof}

From Theorems {\ref{phi-tran} and \ref{embeddings3copy(1)}}, we
infer the following Sobolev-type embedding for the spaces  $\dot{K}_{p}^{\alpha,q}\dot{F}_{\beta}(\mathbb{R}^{n},\{t_{k}\}_{k\in \mathbb{Z}})$.
\begin{thm}
	\label{Sobolev-embedding}
	Let $0< p_{1},p_{2}<\infty $,
 $0<r\leq p_{3}\leq\infty $
and $0<\beta\leq\infty $. Let\ $\{t_{k}\}_{k\in \mathbb{Z}}\in \dot{X}_{\upsilon _{2},\sigma ,p_{2}} $be a $p_{2}$-admissible weight
sequence  with $\sigma =(\sigma _{1}=\theta \left(
\frac{p_{2}}{\theta} \right) ^{\prime },\sigma _{2}\geq p_{2})$ and $\upsilon
_{2}=(\upsilon _{1,2},\upsilon _{2,2})\in \mathbb{R}^{2}$.
Let\ $\{w_{k}\}_{k\in \mathbb{Z}}$ $\in 
\dot{X}_{\upsilon _{1},\sigma ,p_{1}}$ be a $p_{1}$-admissible weight sequence
with $\sigma =(\sigma _{1}=\theta \left( \frac{p_{1}}{\theta} \right) ^{\prime
},\sigma _{2}\geq p_{1})$ and $\upsilon _{1}=(\upsilon _{1,1},\upsilon _{2,1})\in 
\mathbb{R}^{2}$.
Let $\delta_{\{t_{k}\}_{k\in \mathbb{Z}}} $ be as in Lemma \ref{Ap-Property} and $\alpha_{2}>-\frac{
	n\delta_{\{t_{k}\}_{k\in \mathbb{Z}}} }{p_{2}}$. Let $\delta_{\{w_{k}\}_{k\in \mathbb{Z}}} $ be as in Lemma \ref{Ap-Property} and $\alpha_{1}>-\frac{
	n\delta_{\{w_{k}\}_{k\in \mathbb{Z}}} }{p_{1}}$.  Assume that 
\begin{equation}
2^{-k\alpha_{1}}w_{k,Q}(p_{1})\lesssim 2^{-k\alpha_{2}} t_{k,Q}(p_{2})\label{newexp1.1}
\end{equation}

for all $Q\in \mathcal{Q}$ with $\ell (Q)=2^{-k}, k \in \mathbb{Z}$.
\textit{Let }$0<\theta\leq p_{2}<p_{1}<\infty $ and $\alpha _{2}\geq \alpha _{1}$ or $0<\theta \leq  p_{1}\leq
p_{2}<\infty $ and 
\begin{equation}
\alpha _{2}+\frac{n}{p_{2}}\geq \alpha _{1}+\frac{n}{p_{1}}.  \label{newexp1.2}
\end{equation}%
Then%
	\begin{equation*}
	\dot{K}_{p_{2}}^{\alpha_{2},r}\dot{F}_{\mu}(\mathbb{R}^{n},\{t_{k}\}_{k\in \mathbb{Z}})\hookrightarrow \dot{K}_{p_{1}}^{\alpha_{1},p_{3}}\dot{F}_{\beta}(\mathbb{R}^{n},\{w_{k}\}_{k\in \mathbb{Z}}),
	\end{equation*}
	where%
	\begin{equation*}
	\mu =\left\{ 
	\begin{array}{ccc}
	\beta , & \text{if} & 0<p_{1}\leq p_{2}<\infty \text{ and }\alpha _{2}+\frac{n}{p_{2}}%
	=\alpha _{1}+\frac{n}{p_{1}}, \\ 
	\infty , &  & \text{otherwise.}%
	\end{array}%
	\right.
	\end{equation*}
\end{thm}

From Theorem\ \ref{Sobolev-embedding}\ and the fact that $\dot{K}_{p_{1}}^{\alpha_{1},p_{1}}\dot{F}_{\beta}(\mathbb{R}^{n},\{w_{k}\}_{k\in \mathbb{Z}})=\dot{F}_{p_{1},\beta}(\mathbb{R}^{n},\{w_{k}\}_{k\in \mathbb{Z}})$ we immediately arrive at
the following results.

\begin{thm}
	\label{Sobolev-embedding-sequenceherz4}Let $\theta, p_{2},p_{1},\beta  $ $\{t_{k}\}_{k\in \mathbb{Z}}$ and $ \{w_{k}\}_{k\in \mathbb{Z}},\alpha_{2} $ be as in Theorem  \ref{Sobolev-embedding}.  Assume that 
	\begin{equation*}
	w_{k,Q}(p_{1})\lesssim 2^{-k\alpha_{2}} t_{k,Q}(p_{2})
	\end{equation*}
	
	for all $Q\in \mathcal{Q}$ with $\ell (Q)=2^{-k}, k \in \mathbb{Z}$. \newline
	$\mathrm{(i)}$ Let  $0<\theta \leq p_{2}<p_{1}<\infty $ and 
	$\alpha _{2}\geq 0$.
	Then we have
	\begin{equation*}
	\dot{K}_{p_{2}}^{\alpha_{2},p_{1}}\dot{F}_{\infty}(\mathbb{R}^{n},\{t_{k}\}_{k\in \mathbb{Z}})\hookrightarrow \dot{F}_{p_{1},\beta}(\mathbb{R}^{n},\{w_{k}\}_{k\in \mathbb{Z}}).
	\end{equation*}
	$\mathrm{(ii)}$ Let $0<\theta \leq p_{1}<p_{2}<\infty $
	 Assume that 
	\begin{equation*}
	\alpha _{2}\geq\frac{n}{p_{1}}-\frac{n}{p_{2}}.
	\end{equation*}%
	Then we have
	\begin{equation*}
	\dot{K}_{p_{2}}^{\alpha_{2},p_{1}}\dot{F}_{\mu}(\mathbb{R}^{n},\{t_{k}\}_{k\in \mathbb{Z}})\hookrightarrow \dot{F}_{p_{1},\beta}(\mathbb{R}^{n},\{w_{k}\}_{k\in \mathbb{Z}}).
	\end{equation*}
	where%
	\begin{equation*}
	\mu =\left\{ 
	\begin{array}{ccc}
	\beta , & \text{if} &\alpha _{2}+\frac{n}{p_{2}}%
	=\frac{n}{p_{1}}, \\ 
	\infty , &  & \text{otherwise.}%
	\end{array}%
	\right.
	\end{equation*}
\end{thm}

\begin{cor}
	Let $\theta, p_{2},p_{1},\beta, $ $\{t_{k}\}_{k\in \mathbb{Z}}, \{w_{k}\}_{k\in \mathbb{Z}}$ be as in Theorem  \ref{Sobolev-embedding-sequenceherz4}.  Assume that 
	\begin{equation*}
	w_{k,Q}(p_{1})\lesssim t_{k,Q}(p_{2})
	\end{equation*}
	
	for all $Q\in \mathcal{Q}$ with $\ell (Q)=2^{-k}, k \in \mathbb{Z}$. Assume that $0<\theta \leq p_{2}< p_{1}< \infty$. Then%
	
	\begin{equation*}
	\dot{F}_{p_{2},\infty}(\mathbb{R}^{n},\{t_{k}\}_{k\in \mathbb{Z}})\hookrightarrow\dot{K}_{p_{2}}^{0,p_{1}}\dot{F}_{\beta}(\mathbb{R}^{n},\{t_{k}\}_{k\in \mathbb{Z}}) \hookrightarrow\dot{F}_{p_{1},\beta}(\mathbb{R}^{n},\{w_{k}\}_{k\in \mathbb{Z}}).
	\end{equation*}
	
\end{cor}

\begin{proof}
	The desired embeddings follow immediately by setting  $\alpha_{2} =0$ in Theorem \ref
	{Sobolev-embedding-sequenceherz4}, together with the fact that
	\begin{align*}
	\dot{F}_{p_{2},\beta}(\mathbb{R}^{n},\{t_{k}\}_{k\in \mathbb{Z}})=&\dot{K}_{p_{2}}^{0,p_{2}}\dot{F}_{\beta}(\mathbb{R}^{n},\{t_{k}\}_{k\in \mathbb{Z}}) \\
	\hookrightarrow 
	&\dot{K}_{p_{2}}^{0,p_{1}}\dot{F}_{\beta}(\mathbb{R}^{n},\{t_{k}\}_{k\in \mathbb{Z}}) \\
	\hookrightarrow
	&\dot{F}_{p_{1},\beta}(\mathbb{R}^{n},\{w_{k}\}_{k\in \mathbb{Z}}).
	\end{align*}%
	The proof is complete.
\end{proof}
\section{Atomic and molecular decompositions}

We will use the notation of \cite{FJ90}. We shall say that an operator $A$
is associated with the matrix $\{a_{Q_{k,m}P_{v,h}}\}_{k,v\in \mathbb{Z}
,m,h\in \mathbb{Z}^{n}}$, if for all sequences $\lambda =\{\lambda
_{k,m}\}_{k\in \mathbb{Z},m\in \mathbb{Z}^{n}}\subset \mathbb{C}$,%
\begin{equation*}
A\lambda =\{(A\lambda )_{k,m}\}_{k\in \mathbb{Z},m\in \mathbb{Z}^{n}}=\Big\{%
\sum_{v=-\infty }^{\infty }\sum_{h\in \mathbb{Z}^{n}}a_{Q_{k,m}P_{v,h}}%
\lambda _{v,h}\Big\}_{k\in \mathbb{Z},m\in \mathbb{Z}^{n}}.
\end{equation*}%
Throughout the rest of the paper, we set
\begin{equation*}
J=\frac{n}{\min (1,p,q,\beta)}
\end{equation*}%
in the context of Herz-type Triebel-Lizorkin spaces, and
\begin{equation*}
J=\frac{n}{\min (1,p,q)}
\end{equation*}%
in the context of Herz-type Besov spaces. 
We say that $A$, with associated matrix $\{a_{Q_{k,m}P_{v,h}}\}_{k,v\in 
\mathbb{Z},m,h\in \mathbb{Z}^{n}}$, is almost diagonal on $\dot{f}_{p,q}(%
\mathbb{R}^{n},\{t_{k}\})$ if there exists $\varepsilon >0$ such that%
\begin{equation*}
\sup_{k,v\in \mathbb{Z},m,h\in \mathbb{Z}^{n}}\frac{|a_{Q_{k,m}P_{v,h}}|}{
\omega _{Q_{k,m}P_{v,h}}(\varepsilon )}<\infty ,
\end{equation*}%
where

\begin{align}
&\omega _{Q_{k,m}P_{v,h}}(\varepsilon )  \notag \\
&=\Big(1+\frac{|x_{Q_{k,m}}-x_{P_{v,h}}|}{\max \big(2^{-k},2^{-v}\big)}\Big)%
^{-J-\varepsilon }\left\{ 
\begin{array}{ccc}
2^{(v-k)(\alpha _{2}+(n+\varepsilon )/2)}, & \text{if} & v\leq k, \\ 
2^{(v-k)(\alpha _{1}-(n+\varepsilon )/2-J+n)}, & \text{if} & v>k.%
\end{array}%
\right.  \label{omega-assumption}
\end{align}

Using Lemma {\ref{key-estimate1.1} }the following theorem is a
generalization of {\cite[Theorem 3.3]{FJ90}.}

\begin{thm}
\label{almost-diag-est}Let $\upsilon _{1},\upsilon _{2}\in \mathbb{R}$, $%
0<\theta \leq p<\infty $ and $0<q,\beta\leq\infty $. Let $\{t_{k}\}_{k\in \mathbb{Z}}\in \dot{X}
_{\upsilon ,\sigma ,p}$ be a $p$-admissible weight sequence with $\sigma
_{1}=\theta \left( \frac{p}{\theta} \right) ^{\prime }$ and\ $\sigma _{2}\geq p$. Let $\delta_{\{t_{k}\}_{k\in \mathbb{Z}}} $ be as in Lemma \ref{Ap-Property} and $\alpha>-\frac{
	n\delta_{\{t_{k}\}_{k\in \mathbb{Z}}} }{p}$. 
Any almost diagonal operator $A$ on $\dot{K}_{p}^{\alpha,q}\dot{a}_{\beta}(\mathbb{R}^{n},\{t_{k}\}_{k\in \mathbb{Z}})$
is bounded.
\end{thm}

\begin{proof} By similarity, we consider only the spaces $\dot{K}_{p}^{\alpha,q}\dot{f}_{\beta}(\mathbb{R}^{n},\{t_{k}\}_{k\in \mathbb{Z}})$.
We write $A\equiv A_{0}+A_{1}$ with\ 
\begin{equation*}
(A_{0}\lambda )_{k,m}=\sum_{v=-\infty }^{k}\sum_{h\in \mathbb{Z}
^{n}}a_{Q_{k,m}P_{v,h}}\lambda _{v,h},\quad k\in \mathbb{Z},m\in \mathbb{Z}
^{n}
\end{equation*}
and 
\begin{equation*}
(A_{1}\lambda )_{k,m}=\sum_{v=k+1}^{\infty }\sum_{h\in \mathbb{Z}
^{n}}a_{Q_{k,m}P_{v,h}}\lambda _{v,h},\quad k\in \mathbb{Z},m\in \mathbb{Z}
^{n}.
\end{equation*}

\textit{Estimate of }$A_{0}$. From $\mathrm{\eqref{omega-assumption}}$,\ we
obtain 
\begin{align*}
\big|(A_{0}\lambda )_{k,m}\big|& \leq \sum_{v=-\infty }^{k}\sum_{h\in 
\mathbb{Z}^{n}}2^{(v-k)(\alpha _{2}+(n+\varepsilon )/2)}\frac{|\lambda
_{v,h}|}{\big(1+2^{v}|x_{k,m}-x_{v,h}|\big)^{J+\varepsilon }} \\
& =\sum_{v=-\infty }^{k}2^{(v-k)(\alpha _{2}+(n+\varepsilon )/2)}S_{k,v,m}.
\end{align*}%
In \cite{D20.1} we have proved 
\begin{equation*}
S_{k,v,m}\lesssim \mathcal{M}_{\tau }\big(\sum\limits_{h\in \mathbb{Z}%
^{n}}\lambda _{v,h}\chi _{v,h}\big)(x)
\end{equation*}%
for any $x\in Q_{k,m}$ and any $0\leq v\leq k$. Lemmas \ref{key-estimate1.1herz}  and \ref
{key-estimate1.1} we obtain 
\begin{equation*}
{{\big\|}A_{0}\lambda {\big\|_{\dot{K}_{p}^{\alpha,q}\dot{f}_{\beta}(\mathbb{R}^{n},\{t_{k}\}_{k\in \mathbb{Z}})}}}
\end{equation*}%
is bounded by 
\begin{equation*}
c{{\big\|}\lambda {\big\|_{\dot{K}_{p}^{\alpha,q}\dot{f}_{\beta}(\mathbb{R}^{n},\{t_{k}\}_{k\in \mathbb{Z}})}}.}
\end{equation*}

\textit{Estimate of }$A_{1}$. Again from $\mathrm{\eqref{omega-assumption}}$
,\ we see that 
\begin{align*}
\big|(A_{1}\lambda )_{v,h}\big|& \leq \sum_{v=k+1}^{\infty }\sum_{h\in 
\mathbb{Z}^{n}}2^{(v-k)(\alpha _{1}-\varepsilon /2-J+n/2)}\frac{|\lambda
_{v,h}|}{\big(1+2^{k}|x_{k,m}-x_{v,h}|\big)^{J+\varepsilon }} \\
& =\sum_{v=k+1}^{\infty }2^{(v-k)(\alpha _{1}-\varepsilon
/2-J+n/2)}T_{k,v,m}.
\end{align*}
In \cite{D20.1} we have proved  
\begin{equation*}
T_{k,v,m}\leq c2^{(v-k)n/\tau }\mathcal{M}_{\tau }\big(\sum\limits_{h\in 
\mathbb{Z}^{n}}\lambda _{v,h}\chi _{v,h}\big)(x),\quad v>k,x\in Q_{k,m},
\end{equation*}
where $n/(J+\frac{\varepsilon }{2})<\tau <\min (1,p,q,\beta)$ and the positive
constant $c$ is independent of $v$, $k$ and $m$. Again applying Lemmas \ref{key-estimate1.1herz}  and \ref
{key-estimate1.1} we obtain 
\begin{equation*}
{{\big\|}A_{1}\lambda {\big\|_{\dot{K}_{p}^{\alpha,q}\dot{f}_{\beta}(\mathbb{R}^{n},\{t_{k}\}_{k\in \mathbb{Z}})}}}
\end{equation*}
is bounded by 
\begin{equation*}
c{{\big\|}\lambda{\big\|_{\dot{K}_{p}^{\alpha,q}\dot{f}_{\beta}(\mathbb{R}^{n},\{t_{k}\}_{k\in \mathbb{Z}})}}.}
\end{equation*}
Hence the theorem is proved.
\end{proof}

\begin{defn}
\label{Atom-Def}Let\ $\upsilon _{1},\upsilon _{2},\alpha\in \mathbb{R},0<p<\infty $ and 
$0<q,\beta\leq\infty $. Let $\{t_{k}\}_{k\in \mathbb{Z}}$ be a $p$-admissible weight sequence. Let $%
N=\max \{J-n-\upsilon _{1},-1\}$ and $\upsilon _{2}^{\ast }=\upsilon _{2}-\lfloor
\upsilon _{2}\rfloor $.\newline
$\mathrm{(i)}$\ Let $k\in \mathbb{Z}$ and $m\in \mathbb{Z}^{n}$. A function $%
\varrho _{Q_{k,m}}$ is called an homogeneous smooth synthesis molecule for $%
\dot{K}_{p}^{\alpha,q}\dot{A}_{\beta}(\mathbb{R}^{n},\{t_{k}\}_{k\in \mathbb{Z}})$\ supported near $Q_{k,m}$ if there
exist a real number $\delta \in (\upsilon _{2}^{\ast },1]$ and a real number $%
M\in (J,\infty )$ such that 
\begin{equation}
\int_{\mathbb{R}^{n}}x^{\gamma }\varrho _{Q_{k,m}}(x)dx=0\text{\quad if\quad }
0\leq |\gamma |\leq N,  \label{mom-cond}
\end{equation}
\begin{equation}
|\varrho _{Q_{k,m}}(x)|\leq 2^{\frac{kn}{2}}(1+2^{k}|x-x_{Q_{k,m}}|)^{-\max
(M,M-\upsilon _{1})},  \label{cond1}
\end{equation}
\begin{equation}
|\partial ^{\gamma }\varrho _{Q_{k,m}}(x)|\leq 2^{k(|\gamma |+\frac{1}{2}
)}(1+2^{k}|x-x_{Q_{k,m}}|)^{-M}\quad \text{if}\quad |\gamma |\leq \lfloor
\upsilon _{2}\rfloor  \label{cond2}
\end{equation}
and 
\begin{align}
&|\partial ^{\gamma }\varrho _{Q_{k,m}}(x)-\partial ^{\gamma }\varrho
_{Q_{k,m}}(y)|  \label{cond3} \\
&\leq 2^{k(|\gamma |+\frac{1}{2}+\delta )}|x-y|^{\delta }\sup_{|z|\leq
|x-y|}(1+2^{k}|x-z-x_{Q_{k,m}}|)^{-M}\text{\quad if\quad }|\gamma |=\lfloor
\upsilon_{2}\rfloor .  \notag
\end{align}
A collection $\{\varrho _{Q_{k,m}}\}_{k\in \mathbb{Z},m\in \mathbb{Z}^{n}}$
is called a family of homogeneous smooth synthesis molecules for $\dot{K}_{p}^{\alpha,q}\dot{A}_{\beta}(\mathbb{R}^{n},\{t_{k}\}_{k\in \mathbb{Z}})$, if each $\varrho _{Q_{k,m}}$, $k\in 
\mathbb{Z},m\in \mathbb{Z}^{n}$, is an homogeneous smooth synthesis molecule
for $\dot{K}_{p}^{\alpha,q}\dot{A}_{\beta}(\mathbb{R}^{n},\{t_{k}\}_{k\in \mathbb{Z}})$ supported near $Q_{k,m}$. 
\newline
$\mathrm{(ii)}$\ Let $k\in \mathbb{Z}$ and $m\in \mathbb{Z}^{n}$. A function 
$b_{Q_{k,m}}$ is called an homogeneous smooth analysis molecule for $\dot{K}_{p}^{\alpha,q}\dot{A}_{\beta}(\mathbb{R}^{n},\{t_{k}\}_{k\in \mathbb{Z}})$ supported near $Q_{k,m}$ if there exist a $%
\kappa \in ((J-\upsilon_{2})^{\ast },1]$ and an $M\in (J,\infty )$ such that 
\begin{equation}
\int_{\mathbb{R}^{n}}x^{\gamma}b_{Q_{k,m}}(x)dx=0\text{\quad if\quad }0\leq
|\gamma |\leq \left\lfloor \upsilon _{2}\right\rfloor ,  \label{mom-cond2}
\end{equation}
\begin{equation}
|b_{Q_{k,m}}(x)|\leq 2^{\frac{kn}{2}}(1+2^{k}|x-x_{Q_{k,m}}|)^{-\max
(M,M+n+\upsilon _{2}-J)},  \label{cond1.1}
\end{equation}
\begin{equation}
|\partial ^{\gamma }b_{Q_{k,m}}(x)|\leq 2^{k(|\gamma|+\frac{n}{2}
)}(1+2^{k}|x-x_{Q_{k,m}}|)^{-M}\quad \text{if}\quad |\gamma |\leq N
\label{cond1.2}
\end{equation}
and 
\begin{align}
&|\partial ^{\gamma }b_{Q_{k,m}}(x)-\partial ^{\gamma }b_{Q_{k,m}}(y)|
\label{cond1.3} \\
&\leq 2^{k(|\gamma |+\frac{n}{2}+\kappa )}|x-y|^{\kappa }\sup_{|z|\leq
|x-y|}(1+2^{k}|x-z-x_{Q_{k,m}}|)^{-M}\text{\quad if\quad }|\gamma |=N.  \notag
\end{align}
A collection $\{b_{Q_{k,m}}\}_{k\in \mathbb{Z},m\in \mathbb{Z}^{n}}$ is
called a family of homogeneous smooth analysis molecules for $\dot{K}_{p}^{\alpha,q}\dot{A}_{\beta}(\mathbb{R}^{n},\{t_{k}\}_{k\in \mathbb{Z}})$, if each $b_{Q_{k,m}}$, $k\in \mathbb{Z},m\in 
\mathbb{Z}^{n}$, is an homogeneous smooth synthesis molecule for $\dot{K}_{p}^{\alpha,q}\dot{A}_{\beta}(\mathbb{R}^{n},\{t_{k}\}_{k\in \mathbb{Z}})$ supported near $Q_{k,m}$.
\end{defn}

We will use the notation $\{b_{k,m}\}_{k\in \mathbb{Z},m\in \mathbb{Z}^{n}}$
instead of $\{b_{Q_{k,m}}\}_{k\in \mathbb{Z},m\in \mathbb{Z}^{n}}$. The
proof of the following lemma is given in {\cite{D20}.}

\begin{lem}
\label{matrix-est}Let\ $\upsilon _{1},\upsilon _{2},\alpha, J,M,N,\delta ,\kappa ,p,q\ $
and $\beta$ be as in Definition {\ref{Atom-Def}}. Let $\{t_{k}\}_{k\in \mathbb{Z}}$ be a $p$
-admissible weight sequence. Suppose $\{\varrho _{v,h}\}_{v\in \mathbb{Z}
,h\in \mathbb{Z}^{n}}$ is a family of smooth synthesis molecules for $\dot{K}_{p}^{\alpha,q}\dot{A}_{\beta}(\mathbb{R}^{n},\{t_{k}\}_{k\in \mathbb{Z}})$ and $\{b_{k,m}\}_{k\in \mathbb{Z},m\in 
\mathbb{Z}^{n}}$\ is a family of homogeneous smooth analysis molecules for\ $%
\dot{K}_{p}^{\alpha,q}\dot{A}_{\beta}(\mathbb{R}^{n},\{t_{k}\}_{k\in \mathbb{Z}})$. Then there exist a positive real
number $\varepsilon _{1}$ and a positive constant $c$ such that 
\begin{equation*}
\left\vert \langle \varrho _{v,h},b_{k,m}\rangle \right\vert \leq c\text{ }
\omega _{Q_{k,m}P_{v,h}}(\varepsilon ),\quad k,v\in \mathbb{Z},h,m\in 
\mathbb{Z}^{n}
\end{equation*}
if $\varepsilon \leq \varepsilon _{1}$.
\end{lem}

As an immediate consequence, we have the following analogues of the
corresponding results on \cite[Corollary\ B.3]{FJ90}.

\begin{cor}
Let\ $\upsilon _{1},\upsilon _{2},\alpha, J,M,N,\delta ,\kappa ,p,\beta $ and $q$ be as in
Definition {\ref{Atom-Def}}. Let $\{t_{k}\}_{k\in \mathbb{Z}}$ be a $p$-admissible weight
sequence. Let $\Phi $ and $\varphi $ satisfy, respectively $\mathrm{\ %
\eqref{Ass1}}$ and $\mathrm{\eqref{Ass2}}$.\newline
$\mathrm{(i)}$\ If $\{\varrho _{k,m}\}_{k\in \mathbb{Z},m\in \mathbb{Z}^{n}}$
is a family of homogeneous synthesis molecules for  $\dot{K}_{p}^{\alpha,q}\dot{A}_{\beta}(\mathbb{R}^{n},\{t_{k}\}_{k\in \mathbb{Z}})$, then the operator $A$ with
matrix\ $a_{Q_{k,m}P_{v,h}}=\langle \varrho _{v,h},\varphi _{k,m}\rangle $, $%
k,v\in \mathbb{Z},m,h\in \mathbb{Z}^{n}$, is almost diagonal.\newline
$\mathrm{(ii)}$ If $\{b_{k,m}\}_{k\in \mathbb{Z},m\in \mathbb{Z}^{n}}$ is a
family of homogeneous smooth analysis molecules for the 
spaces $\dot{K}_{p}^{\alpha,q}\dot{A}_{\beta}(\mathbb{R}^{n},\{t_{k}\}_{k\in \mathbb{Z}})$, then \ the operator\ $A$,
with matrix $a_{Q_{k,m}P_{v,h}}=\langle \varphi _{v,h},b_{Q_{k,m}}\rangle $, 
$k,v\in \mathbb{Z},m,h\in \mathbb{Z}^{n}$, is almost diagonal.
\end{cor}

Let $f\in \dot{K}_{p}^{\alpha,q}\dot{A}_{\beta}(\mathbb{R}^{n},\{t_{k}\}_{k\in \mathbb{Z}}) $and $\{b_{k,m}\}_{k\in 
\mathbb{Z},m\in \mathbb{Z}^{n}}$ be a family of homogeneous\ smooth analysis
molecules. To prove that $\langle f,b_{Q_{k,m}}\rangle $, $k\in \mathbb{Z}%
,m\in \mathbb{Z}^{n}$, is well defined for all homogeneous smooth analysis
molecules for $\dot{K}_{p}^{\alpha,q}\dot{A}_{\beta}(\mathbb{R}^{n},\{t_{k}\}_{k\in \mathbb{Z}})$, we need the
following result, which proved in {\cite[Lemma 5.4]{BoHo06}}. Suppose that $%
\Phi $ is a smooth analysis (or synthesis) molecule supported near $Q\in 
\mathcal{Q}$ . Then there exists a sequence $\{\varphi _{k}\}_{k\in \mathbb{N%
}}\subset \mathcal{S(}\mathbb{R}^{n}\mathcal{)}$ and $c>0$ such that $%
c\varphi _{k}$ is a smooth analysis (or synthesis) molecule supported near $%
Q $ for every $k$,and $\varphi _{k}(x)\rightarrow \Phi (x)$ uniformly on $%
\mathbb{R}^{n}$ as $k\rightarrow \infty $.

Now we have the following smooth molecular characterization of  $\dot{K}_{p}^{\alpha,q}\dot{A}_{\beta}(\mathbb{R}^{n},\{t_{k}\}_{k\in \mathbb{Z}})$. We refer the reader to {\cite{D20}}
{for the corresponding result for Besov spaces.}

\begin{thm}
\label{molecules-dec}Let $\upsilon _{1}$, $\upsilon _{2}\in \mathbb{R},0<\theta
\leq p<\infty $ and $0<q,\beta\leq\infty $. Let $\{t_{k}\}_{k\in \mathbb{Z}}\in \dot{X}_{\upsilon
,\sigma ,p}$ be a $p$-admissible weight sequence with $\sigma _{1}=\theta
\left( \frac{p}{\theta} \right) ^{\prime }$ and\ $\sigma _{2}\geq p$. Let $%
J,M,N,\delta $ and $\kappa $ be as in Definition {\ref{Atom-Def}}. Let $\delta_{\{t_{k}\}_{k\in \mathbb{Z}}} $ be as in Lemma \ref{Ap-Property} and $\alpha>-\frac{
	n\delta_{\{t_{k}\}_{k\in \mathbb{Z}}} }{p}$. \newline
$\mathrm{(i)}$\ If $f=\sum_{v=-\infty }^{\infty }\sum_{h\in \mathbb{Z}
^{n}}\varrho _{v,h}\lambda _{v,h}$, where $\{\varrho _{v,h}\}_{v\in \mathbb{%
Z },h\in \mathbb{Z}^{n}}$ is a family of homogeneous smooth synthesis
molecules for $\dot{K}_{p}^{\alpha,q}\dot{A}_{\beta}(\mathbb{R}^{n},\{t_{k}\}_{k\in \mathbb{Z}})$, then for all $%
\lambda \in \dot{K}_{p}^{\alpha,q}\dot{a}_{\beta}(\mathbb{R}^{n},\{t_{k}\}_{k\in \mathbb{Z}})$ 
\begin{equation*}
{{\big\|}f{\big\|_{\dot{K}_{p}^{\alpha,q}\dot{A}_{\beta}(\mathbb{R}^{n},\{t_{k}\}_{k\in \mathbb{Z}})}}}\lesssim {\big\|
\lambda{\big\|_{\dot{K}_{p}^{\alpha,q}\dot{a}_{\beta}(\mathbb{R}^{n},\{t_{k}\}_{k\in \mathbb{Z}})}}.}
\end{equation*}
$\mathrm{(ii)}$\ Let $\{b_{k,m}\}_{k\in \mathbb{Z},m\in \mathbb{Z}^{n}}$ be
a family of homogeneous\ smooth analysis molecules.\ Then for all\ $f\in 
\dot{K}_{p}^{\alpha,q}\dot{A}_{\beta}(\mathbb{R}^{n},\{t_{k}\}_{k\in \mathbb{Z}})$ 
\begin{equation*}
{{\big\|}\{\langle f,b_{k,m}\rangle \}_{k\in \mathbb{Z},m\in \mathbb{Z}%
^{n}}{\big\|_{\dot{K}_{p}^{\alpha,q}\dot{A}_{\beta}(\mathbb{R}^{n},\{t_{k}\}_{k\in \mathbb{Z}})}}\lesssim {\big\|}f\big\|_{\dot{K}_{p}^{\alpha,q}\dot{A}_{\beta}(\mathbb{R}^{n},\{t_{k}\}_{k\in \mathbb{Z}})}.}
\end{equation*}
\end{thm}

\begin{proof}
The proof is a slight variant of \cite{FJ90} and \cite{D20}{. }For the
convenience of the reader, we give some details.

\textit{Step 1. Proof of }\textrm{(i)}. By\ \cite[Lemma 2.1]{YY2}. we can
write 
\begin{equation*}
\varrho _{v,h}=\sum_{k=-\infty }^{\infty }2^{-kn}\sum_{m\in \mathbb{Z}^{n}} 
\widetilde{\varphi }_{k}\ast \varrho _{v,h}(2^{-k}m)\psi _{k}(\cdot -2^{-k}m)
\end{equation*}
for\ any\ $v\in \mathbb{Z},h\in \mathbb{Z}^{n}$.\ Therefore, 
\begin{equation*}
f=\sum_{k=-\infty }^{\infty }\sum_{m\in \mathbb{Z}^{n}}S_{k,m}\psi
_{k,m}=T_{\psi }S,
\end{equation*}
where $S=\{S_{k,m}\}_{k\in \mathbb{Z},m\in \mathbb{Z}^{n}}$, with 
\begin{equation*}
S_{k,m}=2^{-kn/2}\sum_{v=-\infty }^{\infty }\sum_{h\in \mathbb{Z}^{n}} 
\widetilde{\varphi }_{k}\ast \varrho _{v,h}(2^{-k}m)\lambda _{v,h}.
\end{equation*}
From Theorem {\ref{phi-tran}, }we have 
\begin{equation*}
{{\big\|}f{\big\|_{\dot{K}_{p}^{\alpha,q}\dot{A}_{\beta}(\mathbb{R}^{n},\{t_{k}\}_{k\in \mathbb{Z}})}}={\big\|}T_{\psi
}S\big\|_{\dot{K}_{p}^{\alpha,q}\dot{A}_{\beta}(\mathbb{R}^{n},\{t_{k}\}_{k\in \mathbb{Z}})}\lesssim {\big\|}S|\big\|_{\dot{K}_{p}^{\alpha,q}\dot{a}_{\beta}(\mathbb{R}^{n},\{t_{k}\}_{k\in \mathbb{Z}})}.}
\end{equation*}
But 
\begin{equation*}
S_{k,m}=\sum_{v=-\infty }^{\infty }\sum_{h\in \mathbb{Z}
^{n}}a_{Q_{k,m}P_{v,h}}\lambda _{v,h},
\end{equation*}
with 
\begin{equation*}
a_{Q_{k,m}P_{v,h}}=\langle \varrho _{v,h},\widetilde{\varphi }_{k,m}\rangle
,\quad k,v\in \mathbb{Z},m,h\in \mathbb{Z}^{n}.
\end{equation*}
Applying Lemma {\ref{matrix-est} and Theorem \ref{almost-diag-est} we find
that} 
\begin{equation*}
{{\big\|}S{\big\|_{\dot{K}_{p}^{\alpha,q}\dot{a}_{\beta}(\mathbb{R}^{n},\{t_{k}\}_{k\in \mathbb{Z}})}}\lesssim {\big\|}
\lambda {\big\|_{\dot{K}_{p}^{\alpha,q}\dot{a}_{\beta}(\mathbb{R}^{n},\{t_{k}\}_{k\in \mathbb{Z}})}}.}
\end{equation*}

\textit{Step 2. Proof of }$\mathit{\mathrm{(ii)}}$\textit{.} We have 
\begin{align*}
\langle f,b_{k,m}\rangle &=\sum_{v=-\infty }^{\infty }2^{-vn}\sum_{m\in 
\mathbb{Z}^{n}}\langle \psi _{v}(\cdot -2^{-v}h),b_{k,m}\rangle \widetilde{
\varphi }_{v}\ast f(2^{-v}h) \\
&=\sum_{v=-\infty }^{\infty }\sum_{m\in \mathbb{Z}^{n}}\langle \psi
_{v,h},b_{k,m}\rangle \lambda _{v,h} \\
&=\sum_{v=-\infty }^{\infty }\sum_{h\in \mathbb{Z}^{n}}a_{Q_{k,m}P_{v,h}}
\lambda _{v,h},
\end{align*}
{\ where } 
\begin{equation*}
a_{Q_{k,m}P_{v,h}}=\langle \psi _{v,h},b_{k,m}\rangle ,\quad \lambda
_{v,h}=2^{-vn/2}\widetilde{\varphi }_{v}\ast f(2^{-v}h).
\end{equation*}
Again by Lemma {\ref{matrix-est} and Theorem \ref{almost-diag-est} we find
that} 
\begin{align*}
{\big\|\{\langle f,b_{k,m}\rangle \}_{k\in \mathbb{Z},m\in \mathbb{Z}^{n}}\big\|_{\dot{K}_{p}^{\alpha,q}\dot{a}_{\beta}(\mathbb{R}^{n},\{t_{k}\}_{k\in \mathbb{Z}})}} &{\lesssim }{\big\|\{\lambda
_{v,h}\}_{v\in \mathbb{Z},h\in \mathbb{Z}^{n}}\big\|_{\dot{K}_{p}^{\alpha,q}\dot{a}_{\beta}(\mathbb{R}^{n},\{t_{k}\}_{k\in \mathbb{Z}})}} \\
&=c{\big\|\{(S_{\varphi })_{v,h}\}_{v\in \mathbb{Z},h\in \mathbb{Z}^{n}}\big\|_{\dot{K}_{p}^{\alpha,q}\dot{a}_{\beta}(\mathbb{R}^{n},\{t_{k}\}_{k\in \mathbb{Z}})}.}
\end{align*}
Applying Theorem {\ref{phi-tran} we find that} 
\begin{equation*}
{{\big\|}\{\langle f,b_{k,m}\rangle \}_{k\in \mathbb{Z},m\in \mathbb{Z}%
^{n}}{\big\|_{\dot{K}_{p}^{\alpha,q}\dot{a}_{\beta}(\mathbb{R}^{n},\{t_{k}\}_{k\in \mathbb{Z}})}}\lesssim {\big\|}f{\big\|_{\dot{K}_{p}^{\alpha,q}\dot{A}_{\beta}(\mathbb{R}^{n},\{t_{k}\}_{k\in \mathbb{Z}})}}.}
\end{equation*}
The proof is complete.
\end{proof}

\begin{defn}
Let $\upsilon _{1},\upsilon _{2},\alpha\in \mathbb{R},0<p<\infty ,0<q,\beta\leq\infty $ and\ $%
N=\max \{J-n-\upsilon _{1},-1\}$. Let $\{t_{k}\}_{k\in \mathbb{Z}}$ be a $p$-admissible weight
sequence. A function $a_{Q_{v,m}}$ is called an homogeneous smooth atom for $%
\dot{K}_{p}^{\alpha,q}\dot{A}_{\beta}(\mathbb{R}^{n},\{t_{k}\}_{k\in \mathbb{Z}})$ supported near $Q_{k,m}$, $k\in 
\mathbb{Z}$ and $m\in \mathbb{Z}^{n}$, if 
\begin{equation}
\supp a_{Q_{k,m}}\subseteq 3Q_{k,m}  \label{supp-cond}
\end{equation}

\begin{equation}
|\partial ^{\gamma}a_{Q_{k,m}}(x)|\leq 2^{kn(|\gamma |+1/2)}\text{\quad
if\quad }0\leq |\gamma|\leq \max (0,1+\lfloor \upsilon _{2}\rfloor ),\quad
x\in \mathbb{R}^{n}  \label{diff-cond}
\end{equation}
and if 
\begin{equation}
\int_{\mathbb{R}^{n}}x^{\gamma }a_{Q_{k,m}}(x)dx=0\text{\quad if\quad }0\leq
|\gamma|\leq N\text{ and }k\in \mathbb{Z}.  \label{mom-cond1}
\end{equation}
\end{defn}

A collection $\{a_{Q_{k,m}}\}_{k\in \mathbb{Z},m\in \mathbb{Z}^{n}}$\ is
called a family of homogeneous smooth atoms for $\dot{K}_{p}^{\alpha,q}\dot{A}_{\beta}(\mathbb{R}^{n},\{t_{k}\}_{k\in \mathbb{Z}})$ is an homogeneous smooth atom for $\dot{K}_{p}^{\alpha,q}\dot{A}_{\beta}(\mathbb{R}^{n},\{t_{k}\}_{k\in \mathbb{Z}})$ supported near $Q_{v,m}$. We point
out that in the moment condition $\mathrm{\eqref{mom-cond1}}$ can be
strengthened into that%
\begin{equation*}
\int_{\mathbb{R}^{n}}x^{\gamma }a_{Q_{k,m}}(x)dx=0\text{\quad if\quad }0\leq
|\gamma |\leq \tilde{N}\text{ and }k\in \mathbb{Z}
\end{equation*}%
and the regularity condition $\mathrm{\eqref{diff-cond}}$ can be
strengthened into that%
\begin{equation*}
|\partial ^{\gamma }a_{Q_{k,m}}(x)|\leq 2^{kn(|\gamma |+1/2)}\text{\quad
if\quad }0\leq |\gamma|\leq \tilde{K},\quad x\in \mathbb{R}^{n},
\end{equation*}%
where $\tilde{K}$ and $\tilde{N}$ are arbitrary fixed integer satisfying $%
\tilde{K}\geq \max (0,1+\lfloor \upsilon _{2}\rfloor )$ and $\tilde{N}\geq
\max \{J-n-\upsilon _{1},-1\}$. If an atom $a$ is supported near $Q_{v,m}$,
then we denote it by $a_{v,m}$.\vskip5pt

Now we come to the atomic decomposition theorem, see {\cite{D20}} {for Besov
spaces and the same arguments are true for Triebel-Lizorkin spaces.}

\begin{thm}
\label{atomic-dec}Let $\upsilon _{1}$, $\upsilon _{2}\in \mathbb{R}$, $0<\theta
\leq p<\infty $, $0<q,\beta\leq\infty $. Let $\{t_{k}\}_{k\in \mathbb{Z}}\in \dot{X}_{\upsilon
,\sigma ,p}$ be a $p$-admissible weight sequence with $\sigma _{1}=\theta
\left( p/\theta \right) ^{\prime }$ and\ $\sigma _{2}\geq p$. Let $\delta_{\{t_{k}\}_{k\in \mathbb{Z}}} $ be as in Lemma \ref{Ap-Property} and $\alpha>-\frac{
	n\delta_{\{t_{k}\}_{k\in \mathbb{Z}}} }{p}$. Then for each $%
f\in \dot{K}_{p}^{\alpha,q}\dot{A}_{\beta}(\mathbb{R}^{n},\{t_{k}\}_{k\in \mathbb{Z}})$, there exist a family\ $%
\{\varrho _{k,m}\}_{k\in \mathbb{Z},m\in \mathbb{Z}^{n}}$ of homogeneous
smooth atoms for $\dot{K}_{p}^{\alpha,q}\dot{A}_{\beta}(\mathbb{R}^{n},\{t_{k}\}_{k\in \mathbb{Z}})$ and $\lambda
=\{\lambda _{k,m}\}_{k\in \mathbb{Z},m\in \mathbb{Z}^{n}}\in \dot{K}_{p}^{\alpha,q}\dot{a}_{\beta}(\mathbb{R}^{n},\{t_{k}\}_{k\in \mathbb{Z}})$ such that 
\begin{equation*}
f=\sum\limits_{k=-\infty }^{\infty }\sum\limits_{m\in \mathbb{Z}^{n}}\lambda
_{k,m}\varrho _{k,m},\text{\quad converging in }\mathcal{S}_{\infty
}^{\prime }(\mathbb{R}^{n})
\end{equation*}
and 
\begin{equation*}
{{\big\|}\{\lambda _{k,m}\}_{k\in \mathbb{Z},m\in \mathbb{Z}^{n}}|{\big\|_{\dot{K}_{p}^{\alpha,q}\dot{a}_{\beta}(\mathbb{R}^{n},\{t_{k}\}_{k\in \mathbb{Z}})}}\lesssim {\big\|}f\big\|_{\dot{K}_{p}^{\alpha,q}\dot{A}_{\beta}(\mathbb{R}^{n},\{t_{k}\}_{k\in \mathbb{Z}})}.}
\end{equation*}
Conversely, for any family of homogeneous smooth atoms for $\dot{K}_{p}^{\alpha,q}\dot{A}_{\beta}(\mathbb{R}^{n},\{t_{k}\}_{k\in \mathbb{Z}})$ and $\lambda =\{\lambda _{k,m}\}_{k\in \mathbb{Z}
,m\in \mathbb{Z}^{n}}\in \dot{K}_{p}^{\alpha,q}\dot{a}_{\beta}(\mathbb{R}^{n},\{t_{k}\}_{k\in \mathbb{Z}})$ 
\begin{equation*}
{{\big\|}\sum\limits_{k=-\infty }^{\infty }\sum\limits_{m\in \mathbb{Z}
^{n}}\lambda _{k,m}\varrho _{k,m} %
\big\|_{\dot{K}_{p}^{\alpha,q}\dot{A}_{\beta}(\mathbb{R}^{n},\{t_{k}\}_{k\in \mathbb{Z}})}\lesssim {\big\|}\{\lambda _{k,m}\}_{k\in \mathbb{Z},m\in \mathbb{Z}
^{n}}{\big\|_{\dot{K}_{p}^{\alpha,q}\dot{a}_{\beta}(\mathbb{R}^{n},\{t_{k}\}_{k\in \mathbb{Z}})}}.}
\end{equation*}
\end{thm}

\section{The non-homogeneous space\textbf{\ }$\dot{K}_{p}^{\alpha,q}A_{\beta}(\mathbb{R}
^{n},\{t_{k}\}_{k\in \mathbb{N}_{0}})$}

In this section, we present the inhomogeneous version of our results given
above. Let $\Phi ,\psi ,\varphi $ and $\Psi $ satisfy 
\begin{equation}
\Phi ,\Psi ,\varphi ,\psi \in \mathcal{S}(\mathbb{R}^{n})  \label{Ass1.1}
\end{equation}%
\begin{equation}
\supp\mathcal{F}\Phi \cup \supp\mathcal{F}\Psi \subset 
\overline{B(0,2)},\text{\quad }|\mathcal{F}(\Phi )(\xi )|,|\mathcal{F}(\Psi
)(\xi )|\geq c,  \label{Ass2.1}
\end{equation}%
if $|\xi |\leq 5/3$ and

\begin{equation}
\supp\mathcal{F}\varphi \cup \supp\mathcal{F}\psi
\subset \overline{B(0,2)}\backslash B(0,1/2),\text{\quad }|\mathcal{F}%
(\varphi )(\xi )|,|\mathcal{F}(\psi )(\xi )|\geq c,  \label{Ass3.1}
\end{equation}%
if $3/5\leq |\xi |\leq 5/3$, such that%
\begin{equation}
\overline{\mathcal{F}(\Phi )(\xi )}\mathcal{F}(\Psi )(\xi
)+\sum_{k=1}^{\infty }\overline{\mathcal{F(}\varphi )(2^{-k}\xi )}\mathcal{F}%
(\psi )(2^{-k}\xi )=1,\quad \xi \in \mathbb{R}^{n},  \label{Ass4.1}
\end{equation}%
where $c>0$. Let $\Phi ,\varphi \in \mathcal{S}(\mathbb{R}^{n})$\ satisfy,
respectively, $\mathrm{\eqref{Ass2.1}}$ and $\mathrm{\eqref{Ass3.1}}$. We
recall that by \cite[pp. 130--131]{FJ90} or \cite[Lemma 6.9]{FrJaWe01},
there exist functions $\Psi \in \mathcal{S}(\mathbb{R}^{n})$ satisfying $%
\mathrm{\eqref{Ass2.1}}$ and $\psi \in \mathcal{S}(\mathbb{R}^{n})$
satisfying $\mathrm{\eqref{Ass3.1}}$ such that $\mathrm{\eqref{Ass4.1}}$
holds.

The $\varphi $-transform $S_{\varphi }$ is defined by setting 
\begin{equation*}
(S_{\varphi}f)_{0,m}=\langle f,\Psi _{m}\rangle ,
\end{equation*}
with $\Psi _{m}(x)=\Psi (x-m)$ and 
\begin{equation*}
(S_{\varphi }f)_{k,m}=\langle f,\varphi _{k,m}\rangle,
\end{equation*}
where $\varphi _{k,m}(x)=2^{kn/2}\varphi (2^{k}x-m)$, $k\in \mathbb{N}$ and $%
m\in \mathbb{Z}^{n}$. The inverse $\varphi $-transform $T_{\psi }$ is
defined by 
\begin{equation*}
T_{\psi }\lambda =\sum_{m\in \mathbb{Z}^{n}}\lambda _{0,m}\Psi
_{m}+\sum_{k=1}^{\infty }\sum_{m\in \mathbb{Z}^{n}}\lambda _{k,m}\psi _{k,m},
\end{equation*}%
with $\lambda =\{\lambda _{k,m}\}_{k\in \mathbb{N}_{0},m\in \mathbb{Z}
^{n}}\subset \mathbb{C}$, see again \cite{FJ90}.

Now we present the inhomogenous version of Definition \ref{Tyulenev-class}.

\begin{defn}
\label{Tyulenev-class-inho}Let $\upsilon _{1}$, $\upsilon _{2}\in \mathbb{R}$, $%
\sigma _{1}$, $\sigma _{2}$ $\in (0,+\infty ]$, $\upsilon =(\upsilon _{1},\upsilon
_{2})$\ and let $\sigma =(\sigma _{1},\sigma _{2})$. We let $X_{\upsilon
,\sigma ,p}=X_{\upsilon ,\sigma ,p}(\mathbb{R}^{n})$ denote the set of $p$
-admissible weight sequences $\{t_{k}\}_{k\in \mathbb{N}_{0}}$ satisfying $%
\mathrm{\eqref{Asum1}}$ and $\mathrm{\eqref{Asum2}}$ for any $0\leq k\leq j$
, with constants $C_{1},C_{2}>0$ are independent of both the indexes $k$ and 
$j$.
\end{defn}

\begin{ex}
A sequence $\{\gamma _{j}\}_{j\in \mathbb{N}_{0}}$ of positive real numbers
is said to be admissible if there exist two positive constants $d_{0}$ and $%
d_{1}$ such that 
\begin{equation}
d_{0}\gamma _{j}\leq \gamma _{j+1}\leq d_{1}\gamma _{j},\quad j\in \mathbb{N}
_{0}.  \label{Farkas-Leop}
\end{equation}
For an admissible sequence $\{\gamma _{j}\}_{j\in \mathbb{N}_{0}}$, let 
\begin{equation*}
\underline{\gamma }_{j}=\inf_{k\geq 0}\frac{\gamma _{j+k}}{\gamma _{k}}\quad 
\text{and}\quad \overline{\gamma }_{j}=\sup_{k\geq 0}\frac{\gamma _{j+k}}{
\gamma _{k}},\quad j\in \mathbb{N}_{0}.
\end{equation*}
Let 
\begin{equation*}
\alpha _{\gamma }=\lim_{j\longrightarrow \infty }\frac{\log \overline{\gamma 
}_{j}}{j}\quad \text{and}\quad \beta _{\gamma }=\lim_{j\longrightarrow
\infty }\frac{\log \underline{\gamma }_{j}}{j},
\end{equation*}
be the upper and lower Boyd index of the given sequence $\{\gamma
_{j}\}_{j\in \mathbb{N}_{0}}$, respectively. Then 
\begin{equation*}
\underline{\gamma }_{j}\gamma _{k}\leq \gamma _{j+k}\leq \overline{\gamma }
_{j}\gamma _{k},\quad j,k\in \mathbb{N}_{0}
\end{equation*}
and for each $\varepsilon >0$, 
\begin{equation*}
c_{1}2^{(\beta _{\gamma }-\varepsilon )j}\leq \underline{\gamma }_{j}\leq 
\overline{\gamma }_{j}\leq c_{2}2^{(\alpha _{\gamma }+\varepsilon )j},\quad
j\in \mathbb{N}_{0}
\end{equation*}
for some constants $c_{1}=c_{1}(\varepsilon )>0$ and $c_{2}=c_{2}(
\varepsilon )>0$. Also, $\underline{\gamma }_{1}$ and $\overline{\gamma }
_{1} $ are the best possible constants $d_{0}$ and $d_{1}$ in %
\eqref{Farkas-Leop}, respectively. \newline
Clearly the sequence $\{\gamma _{j}\}_{j\in \mathbb{N}_{0}}$ lies in $%
X_{\upsilon ,\sigma ,p}$ for$\ \upsilon _{1}=\beta _{\gamma }-\varepsilon
,\upsilon _{2}=\upsilon _{\gamma }+\varepsilon $ and $0<p,\sigma _{1},\sigma
_{2}\leq \infty $.\newline
These type of admissible sequences are used in \cite{FL06} to study Besov
and Triebel-Lizorkin spaces in terms of a generalized smoothness, see also 
\cite{HaS08}.\newline
Let us consider some examples of admissible sequences. The sequence $%
\{\gamma _{j}\}_{j\in \mathbb{N}_{0}}$, 
\begin{equation*}
\gamma _{j}=2^{sj}(1+j)^{b}(1+\log (1+j))^{c},\quad j\in \mathbb{N}_{0}
\end{equation*}
with arbitrary fixed real numbers $s,b$ and $c$ is a an admissible sequence
with 
\begin{equation*}
\beta _{\gamma }=\alpha _{\gamma }=s.
\end{equation*}
\end{ex}

\begin{ex}
\label{Example1 copy(1)}Let $0<r<p<\infty $, a weight $\omega ^{p}\in A_{ 
\frac{p}{r}}(\mathbb{R}^{n})$ and 
\begin{equation*}
\{s_{k}\}=\{2^{ks}\omega ^{p}(2^{-k})\}_{k\in \mathbb{N}_{0}},\quad s\in 
\mathbb{R}.
\end{equation*}
Obviously, $\{s_{k}\}_{k\in \mathbb{N}_{0}}$ lies in $X_{\upsilon ,\sigma ,p}$
for $\upsilon _{1}=\upsilon _{2}=s$, $\sigma =(r(p/r)^{\prime },p)$.\ 
\end{ex}

Now, we define the spaces under consideration.

\begin{defn}
\label{B-F-def-inh}Let $0<p< \infty $ and $0<q,\beta \leq \infty $. Let $%
\{t_{k}\}_{k\in \mathbb{N}_{0}}$ be a $p$-admissible weight sequence. Let $%
\Phi ,\varphi \in \mathcal{S}(\mathbb{R}^{n})$\ satisfy $\mathrm{\ %
\eqref{Ass2.1}}$ and $\mathrm{\eqref{Ass3.1}}$, respectively, and we put $%
\varphi _{k}=2^{kn}\varphi(2^{k}\cdot) ,k\in \mathbb{N}_{0}$.
\newline
$\mathrm{(i)}$ The Herz-type Besov space $\dot{K}^{\alpha,q}_{p}B_{\beta}(\mathbb{R}^{n},\{t_{k}\}_{k\in \mathbb{{N}}_{0}})$\ is
the collection of all $f\in \mathcal{S}_{\infty }^{\prime }(\mathbb{R}^{n})$
\ such that 
\begin{equation*}
\big\|f\big\|_{\dot{K}^{\alpha,q}_{p}B_{\beta}(\mathbb{R}^{n},\{t_{k}\}_{k\in \mathbb{{N}}_{0}})}:=\Big(%
\sum\limits_{k=0 }^{\infty }\big\|t_{k}(\varphi _{k}\ast f)\big\|_{\dot{K}^{\alpha,q}_{p}}^{\beta}\Big)^{\frac{1}{\beta}}<\infty
\end{equation*}%
with the usual modification if $\beta=\infty $.\newline
 \newline
$\mathrm{(ii)}$ The Herz-type
Triebel-Lizorkin space $\dot{K}_{p}^{\alpha,q}F_{\beta}(\mathbb{R}
^{n},\{t_{k}\}_{k\in \mathbb{N}%
_{0}})$ is the collection of all $f\in \mathcal{S}^{\prime }(\mathbb{R}%
^{n}) $\ such that 
\begin{equation*}
\big\|f\big\|_{\dot{K}_{p}^{\alpha,q}F_{\beta}(\mathbb{R}
	^{n},\{t_{k}\}_{k\in \mathbb{N}%
		_{0}})}=\Big\| %
\Big(\sum\limits_{k=0}^{\infty }t_{k}^{\beta}|\varphi _{k}\ast f|^{\beta}\Big) %
^{1/q}\Big\|_{\dot{K}_{p}^{\alpha,q}}<\infty,
\end{equation*}
with the usual modifications if $\beta=\infty $, where $\varphi _{0}$ is
replaced by $\Phi $.
\end{defn}

Now we introduce the inhomogeneous sequence spaces $\dot{K}_{p}^{\alpha,q}a_{\beta}(\mathbb{R}^{n},\{t_{k}\}_{k\in \mathbb{N}_{0}})$. Let $0<p< \infty $ and $0<q\leq
\infty $. Let $\{t_{k}\}_{k\in \mathbb{N}_{0}}$ be a $p$-admissible weight
sequence. Then for all complex valued sequences $\lambda =\{\lambda
_{k,m}\}_{k\in \mathbb{N}_{0},m\in \mathbb{Z}^{n}}\subset \mathbb{C}$ we
define%
\begin{equation*}
\dot{K}_{p}^{\alpha,q}b_{\beta}(\mathbb{R}^{n},\{t_{k}\}_{k\in \mathbb{N}_{0}})=\Big\{\lambda :\big\|%
\lambda \big\|_{\dot{K}_{p}^{\alpha,q}b_{\beta}(\mathbb{R}^{n},\{t_{k}\}_{k\in \mathbb{N}_{0}})}%
<\infty \Big\},
\end{equation*}%
where
\begin{equation*}
\big\|%
\lambda \big\|_{\dot{K}_{p}^{\alpha,q}b_{\beta}(\mathbb{R}^{n},\{t_{k}\}_{k\in \mathbb{N}_{0}})}=\Big(%
\sum_{k=0 }^{\infty }2^{\frac{kn\beta}{2}}\Big\| \sum\limits_{m\in \mathbb{Z}%
	^{n}}|\lambda _{k,m}|t_{k,m}\chi _{k,m}\Big\|_{\dot{K}_{p}^{\alpha,q}}^{\beta}\Big)^{\frac{1}{\beta}},
\end{equation*}
with the usual modifications if $\beta=\infty $ and
\begin{equation*}
\dot{K}_{p}^{\alpha,q}f_{\beta}(\mathbb{R}^{n},\{t_{k}\}_{k\in \mathbb{N}_{0}})=\Big\{\lambda :\big\|%
\lambda \big\|_{\dot{K}_{p}^{\alpha,q}f_{\beta}(\mathbb{R}^{n},\{t_{k}\}_{k\in \mathbb{N}_{0}})}%
<\infty \Big\},
\end{equation*}%
where%
\begin{equation*}
\big\|%
\lambda \big\|_{\dot{K}_{p}^{\alpha,q}f_{\beta}(\mathbb{R}^{n},\{t_{k}\}_{k\in \mathbb{N}_{0}})}=%
\Big\|\Big(\sum_{k=0}^{\infty }2^{\frac{kn\beta}{2}}\sum\limits_{m\in \mathbb{Z}
^{n}}t_{k}^{\beta}|\lambda _{k,m}|^{\beta}\chi _{k,m}\Big)^{\frac{1}{\beta}}\big\|_{\dot{K}_{p}^{\alpha,q}}.
\end{equation*}%
We have the following analogue of Theorem \ref{phi-tran}.

\begin{thm}
\label{phi-tran-inho}Let $\upsilon =(\upsilon _{1},\upsilon _{2})\in \mathbb{R}
^{2},0<\theta \leq p<\infty $ and$\ 0<q,\beta\leq \infty $. Let $\{t_{k}\}_{k\in 
\mathbb{N}_{0}}\in X_{\upsilon ,\sigma ,p}$ be a $p$-admissible weight
sequence with $\sigma =(\sigma _{1}=\theta \left( \frac{p}{\theta}\right) ^{\prime
},\sigma _{2}\geq p)$.\ Let $\varphi $, $\psi $ satisfying $\mathrm{\ %
\eqref{Ass1.1}}$\ through\ $\mathrm{\eqref{Ass4.1}}$. Let $\delta_{\{t_{k}\}_{k\in \mathbb{Z}}} $ be as in Lemma \ref{Ap-Property} and $\alpha>-\frac{
n\delta_{\{t_{k}\}_{k\in \mathbb{Z}}} }{p}$. The operators 
\begin{equation*}
S_{\varphi }:\dot{K}_{p}^{\alpha,q}A_{\beta}(\mathbb{R}^{n},\{t_{k}\}_{k\in \mathbb{N}_{0}})\rightarrow \dot{K}_{p}^{\alpha,q}a_{\beta}(\mathbb{R}^{n},\{t_{k}\}_{k\in \mathbb{N}_{0}})
\end{equation*}
and 
\begin{equation*}
T_{\psi }:\dot{K}_{p}^{\alpha,q}a_{\beta}(\mathbb{R}^{n},\{t_{k}\}_{k\in \mathbb{N}_{0}})\rightarrow
\dot{K}_{p}^{\alpha,q}A_{\beta}(\mathbb{R}^{n},\{t_{k}\}_{k\in \mathbb{N}_{0}})
\end{equation*}
are bounded. Furthermore, $T_{\psi }\circ S_{\varphi }$ is the identity on $%
\dot{K}_{p}^{\alpha,q}A_{\beta}(\mathbb{R}^{n},\{t_{k}\}_{k\in \mathbb{N}_{0}})$.
\end{thm}

As a consequence the analogues of Corollary \ref{Indpendent} are now clear.
We obtain the following useful properties of these function spaces.

\begin{thm}
Let $\upsilon =(\upsilon _{1},\upsilon _{2})\in \mathbb{R}^{2},0<\theta \leq
p<\infty $ and $0<q,\beta\leq \infty $. Let $\{t_{k}\}_{k\in \mathbb{N}_{0}}\in
X_{\upsilon ,\sigma ,p}$ be a $p$-admissible weight sequence with $\sigma
=(\sigma _{1}=\theta \left( \frac{p}{\theta} \right) ^{\prime },\sigma _{2}\geq p)$.  Let $\delta_{\{t_{k}\}_{k\in \mathbb{Z}}} $ be as in Lemma \ref{Ap-Property} and $\alpha>-\frac{
	n\delta_{\{t_{k}\}_{k\in \mathbb{Z}}} }{p}$.
$F_{p,q}(\mathbb{R}^{n},\{t_{k}\}_{k\in \mathbb{N}_{0}})$ are quasi-Banach
spaces. They are Banach spaces if $1\leq p<\infty $ and $1\leq q,\beta \leq\infty $.
\end{thm}

\begin{rem} Let $0<\theta \leq p<\infty $ and $0<q,\beta \leq \infty $. Let\ $\{t_{k}\}_{k\in 
\mathbb{N}_{0}}\in X_{\upsilon ,\sigma ,p}$ be a $p$-admissible weight
sequence with $\sigma =(\sigma _{1}=\theta \left( \frac{p}{\theta} \right) ^{\prime
},\sigma _{2}\geq p)$ and $\upsilon =(\upsilon _{1},\upsilon _{2})\in \mathbb{R}%
^{2}$. Let $\delta_{\{t_{k}\}_{k\in \mathbb{Z}}} $ be as in Lemma \ref{Ap-Property} and $\alpha>-\frac{
	n\delta_{\{t_{k}\}_{k\in \mathbb{Z}}} }{p}$. As in Theorem \ref{embeddings-S-inf} we have the embedding%
\begin{equation*}
\mathcal{S}(\mathbb{R}^{n})\hookrightarrow \dot{K}_{p}^{\alpha,q}A_{\beta}(\mathbb{R}^{n},\{t_{k}\}_{k\in \mathbb{N}_{0}}).
\end{equation*}%
In addition $\mathcal{S}(\mathbb{R}^{n})$ is dense in $\dot{K}_{p}^{\alpha,q}A_{\beta}(\mathbb{R}^{n},\{t_{k}\}_{k\in \mathbb{N}_{0}})$. Also if $0<\theta \leq p<\infty $
and $0<q,\beta \leq\infty $, then 
\begin{equation*}
\dot{K}_{p}^{\alpha,q}A_{\beta}(\mathbb{R}^{n},\{t_{k}\}_{k\in \mathbb{N}_{0}})\hookrightarrow 
\mathcal{S}^{\prime }(\mathbb{R}^{n}).
\end{equation*}%
All the results in Sections 5 and 6 can be extended to the inhomogeneous case.
\end{rem}

\end{document}